\pdfoutput=1
\documentclass[final,pdftex]{siamart171218}

\usepackage{amsmath,amssymb,amscd}
\usepackage{amsfonts}
\usepackage{graphicx}
\usepackage{seqsplit}
\usepackage[hang,centerlast]{subfigure}

\DeclareMathOperator{\lift}{\mathcal{L}}
\newcommand{\tskip}{{t_\mathrm{skip}}}

\newcommand{\deltat}{\delta}

\DeclareMathOperator{\restrict}{\mathcal{R}}

\DeclareMathOperator{\dom}{dom}

\DeclareMathOperator{\diag}{diag}
\DeclareMathOperator{\linspan}{span}

\DeclareMathOperator{\var}{var}

\newcommand{\R}{\mathbb{R}}

\newcommand{\Lint}{\mathbb{L}}
\newcommand{\Hint}{\mathbb{H}}
\newcommand{\id}{I}

\DeclareMathOperator{\rg}{rg}
\DeclareMathOperator{\dist}{dist}
\DeclareMathOperator{\rank}{rank}

\newcommand{\e}{{\mathrm{e}}}
\newcommand{\llin}{\lift_{\mathrm{lin}}}
\newcommand{\lin}{{\mathrm{lin}}}
\newcommand{\lgauss}{\lift_{\mathrm{Gauss}}}
\newcommand{\gauss}{{\mathrm{Gauss}}}

\renewcommand{\d}{\mathop{}\!\mathrm{d}}

\renewcommand{\phi}{\varphi}

\renewcommand{\epsilon}{\varepsilon}
\newcommand{\dtan}{d_\mathrm{tan}}
\newcommand{\dtr}{d_\mathrm{tr}}

\DeclareMathOperator{\spec}{spec}

\newcommand{\order}[1]{\mathcal{O}(#1)}
\DeclareMathOperator{\res}{res}

\newcommand{\code}[1]{{\ttfamily\seqsplit{#1}}}
\newcommand{\scripts}[1]{}
\newsiamthm{assumption}{Assumption}

\title{Convergence of equation-free methods in the case of finite time
  scale separation with application to deterministic and stochastic systems}

\author{Jan Sieber\thanks{College of Engineering, Mathematics and
    Physical Sciences, University of Exeter, North Park Road, Exeter
    (Devon) EX4 4QF, United Kingdom (\url{j.sieber@exeter.ac.uk}).} %
  \and Christian Marschler\thanks{Department of Applied Mathematics
    and Computer Science, Technical University of Denmark,
    Matematiktorvet 303B, DK-2800 Kgs. Lyngby, Denmark}. %
  \and Jens Starke\thanks{Institute of Mathematics, University of
    Rostock, Ulmenstra{\ss}e 69, 18057 Rostock, Germany
    (\url{jens.starke@uni-rostock.de}).}  }

\headers{Convergence of equation-free methods}{Jan Sieber, Christian Marschler, Jens Starke}
\begin{document}

\maketitle

\newcommand{\slugmaster}{%
}

\begin{abstract}
A common approach to studying high-dimensional systems with emergent low-dimensional behavior is based on lift-evolve-restrict maps (called equation-free methods): first, a user-defined lifting operator maps a set of low-dimensional coordinates into the high-dimensional phase space, then the high-dimensional (microscopic) evolution is applied for some time, and finally a user-defined restriction operator maps down into a low-dimensional space again. We prove convergence of equation-free methods for finite time-scale separation with respect to a method parameter, the so-called healing time. Our convergence result justifies equation-free methods as a tool for performing high-level tasks such as bifurcation analysis on high-dimensional systems.

More precisely, if the high-dimensional system has an attracting invariant manifold with smaller expansion and attraction rates in the tangential direction than in the transversal direction (normal hyperbolicity), and restriction and lifting satisfy some generic transversality conditions, then an implicit formulation of the lift-evolve-restrict procedure generates an approximate map that converges to the flow on the invariant manifold for healing time going to infinity. In contrast to all previous results, our result does not require the time scale separation to be large. A demonstration with Michaelis-Menten kinetics shows that the error estimates of our theorem are sharp.

The ability to achieve convergence even for finite time scale separation is especially important for applications involving stochastic systems, where the evolution occurs at the level of distributions, governed by the Fokker-Planck equation. In these applications the spectral gap is typically finite. We investigate a low-dimensional stochastic differential equation where the ratio between the decay rates of fast and slow variables is $2$.
\end{abstract}




\section{Introduction}
\label{sec:introduction}

High-dimensional dynamical systems with time scale separation have
under certain assumptions the potential to be studied and understood
through a reduction to low-dimensional systems. In most cases these
reduction methods are applied directly to the high-dimensional systems
of equations \cite{roberts2014model}. The most common approaches are
referred to as averaging and mean-field approximation
\cite{Stanley1987}, the slaving principle or adiabatic elimination
\cite{Haken1983} in the physics literature.  The aim of
these methods is to reduce the complexity of a high-dimensional (here
also called \emph{microscopic}) system to a relatively simple
low-dimensional (here also called \emph{macroscopic}) system. After
reduction, the long-term dynamics of the system can be analyzed by
studying the low-dimensional macroscopic system, using techniques that
may only be available for low-dimensional deterministic systems (e.g.,
detailed bifurcation analysis). The underlying assumption is that a
trajectory of the microscopic system will rapidly relax onto a
low-dimensional manifold, which it will then track on a longer time
scale, following the slower macroscopic equations. Thus, one speaks of
\emph{slow} variables, which are the coordinates on the slow
low-dimensional manifold, and \emph{fast} variables transversal to the
slow manifold. The notion that the fast variables are ``slaved'' by
the slow variables describes that over long time the microscopic
trajectories track the slow manifold.

The justification for this reduction is simplest and strongest if the
underlying microscopic dynamical system possesses a low-dimensional
attracting invariant manifold. In these cases mathematical theorems on
persistence of invariant manifolds can be applied. Proofs were given
by Fenichel \cite{Fenichel1979} and Hirsch \emph{et al}. \cite{HPS77}
for finite-dimensional smooth dynamical systems such as ordinary
differential equation (ODEs) and maps and by Bates \emph{et
  al}. \cite{BLZ99,BLZ00} for general semiflows (covering certain
classes of partial differential equations). Certain cases of averaging
(such as periodic and quasi-periodic forcing) may also be reduced to
invariant manifold persistence.

The case for model reduction is more subtle if the microscopic system
is stochastic (or more generally, ergodic), for example, if the model
is given by a multi-particle or agent-based simulation. The time-scale
separation for these systems occurs if, for example, the number of
particles is large. A model case for stochastic systems is the
reduction of a high-dimensional system of stochastic differential
equations (SDE) to a low-dimensional SDE acting on a slower time scale
(smaller drift terms and smaller noise amplitude than the microscopic
system). In this case, the arguments for model reduction look formally
similar to the case of attracting manifolds in deterministic systems
\cite{PS08,givon2004}. However, the underlying mathematical
convergence results are not as strong. Two aspects in which the
deterministic results are stronger than the stochastic results will
have implications on convergence results for computational methods:

\emph{Validity for finite time-scale separation}: For invariant
manifolds in a deterministic ODE the time scale separation (let us
call it $\epsilon$) is measured as the ratio between the rate of
attraction along directions tangential to the manifold ($\dtan$) and
transversal to the manifold ($\dtr$, so $\epsilon=\dtan/\dtr$). As
long as this ratio $\epsilon$ is less than unity, the manifold
persists. Let us call this persistent low-dimensional manifold
${\cal C}$. Persistence implies that, even for a finite $\epsilon$, a
reduced model on this manifold ${\cal C}$ exists, describing some
trajectories of the microscopic system with perfect accuracy (those
that lie on ${\cal C}$). In practice the dynamics on the slow manifold
${\cal C}$ is often approximated by an expansion in $\epsilon$.
  
\emph{Shadowing}: Even more, every point $u$ from an open neighborhood
of ${\cal C}$ has a shadowing point $g(u)\in {\cal C}$. The difference
between the trajectories starting from $u$ and $g(u)$ goes to zero in
time with a rate close to $\dtr$. This means that the reduced model
describes \emph{all} nearby trajectories even for positive $\epsilon$
with perfect accuracy except for rapidly decaying terms. The nonlinear
projection $u\mapsto g(u)\in {\cal C}$ is called the \emph{stable
  fiber projection}.

Compared to the above, the precise mathematical convergence statements
in \cite{PS08,givon2004} for stochastic systems with time scale
separation are weaker. They are concerned with the limit
$\epsilon\to0$ and prove that moments of the slow coordinates of the
microscopic trajectory and of the trajectories of the slow model,
derived by a formal expansion in $\epsilon$, converge to each other
for $\epsilon\to0$ \cite{PS08}.
\begin{figure}[ht]
  \centering
  \includegraphics[scale=0.8]{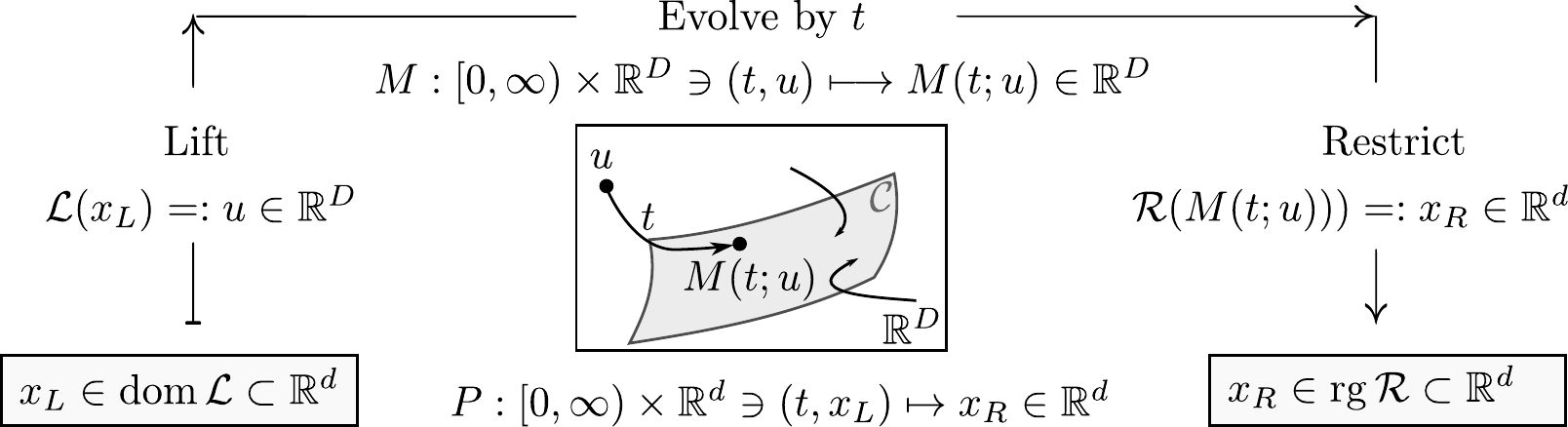}
  \caption{Lift-evolve-restrict map $P(t;\cdot)$ in low-dimensional space
    $\R^d$ used in equation-free computations.}
  \label{fig:lift-evolve-restrict}
\end{figure}
\paragraph{On-demand computation of slow flow --- Equation-free
  framework}
The above mathematical theorems underpin the derivation of approximate
low-dimensional models for large-dimensional systems. However, they
also provide guidance for the convergence analysis of computational
methods that avoid the explicit derivation of a low-dimensional model,
but merely assume its existence. A general framework for analysing
slow-time scale behaviour of systems with time scale separation was
proposed by Kevrekidis \emph{et al}.  under the name ``equation-free
computations'' \cite{kevrekidisgear2003,Gear2005,Kevrekidis2010}.  The
assumption behind equation-free computations is the existence of a
slow low-dimensional description (in $\R^d$) for some macroscopic
quantities of the high-dimensional microscopic system (which is
defined in $\R^D$) that contains a $d$-dimensional invariant manifold
${\cal C}$, which we will call the slow manifold. We do not append a
subscript $\epsilon$ to ${\cal C}$, since our main result will only
assume existence and smoothness of the invariant manifold ${\cal C}$
and its stable fiber projection, but not consider the limit of time
scale separation $\epsilon\to0$.  The framework, illustrated in
\cref{fig:lift-evolve-restrict}, only relies on the availability of a
microscopic time stepper (a map $M(t;\cdot):\R^D\mapsto\R^D$ for
$t\geq0$) that can be called at selected microscopic initial values
$u\in\R^D$. The goal is to compose a macroscopic time stepper
$\Phi_*(\delta;\cdot): \R^d\mapsto\R^d$ for $\delta\in\R$ (possibly
including $\delta<0$) in some coordinates for the slow manifold ${\cal
  C}$, which is then amenable to higher-level tasks such as
bifurcation analysis.

For equation-free computations the user also has to choose two
operators, the \emph{lifting} $\lift:\R^d\mapsto\R^D$ and the
\emph{restriction} $\restrict:\R^D\mapsto\R^d$, which are maps between
the original high-dimensional ($\R^D$) microscopic level and the
low-dimensional ($\R^d$) macroscopic level. The user-defined
  lifting ${\cal L}$ and restriction ${\cal R}$, together with the
  time stepper $M(t;\cdot)$ define the central building block of
  equation-free methodology, the ``lift-evolve-restrict'' map,
  \begin{displaymath}
    P:[0,\infty)\times\dom\lift\ni (t,x_L)\mapsto 
  \restrict(M(t;\lift(x_L)))\in\rg \restrict\mbox{}
\end{displaymath}
(see \cref{fig:lift-evolve-restrict}). For a given
value $x_L\in\R^d$ of macroscopic quantities, one first applies the
lifting $\lift$ to $x_L$ getting a microscopic state $u$, then one
runs the microscopic simulation for time $t$ starting from $u$
(applying the microscopic evolution $M(t;u)$), and finally one applies
the restriction $\restrict$ to the result $M(t;u)$.

The use of the lift-evolve-restrict map $P$ assumes that the
trajectory $t\mapsto M(t;u)$ of the time stepper will be close to the
slow manifold ${\cal C}$ most of the time. Assuming this,
equation-free methods aim to extract information about the slow flow
along ${\cal C}$ by calling the lift-evolve-restrict map $P$
judiciously. The simplest approach would be to use $P(t;\cdot)$ as an
approximation for $M(t;\cdot)$ restricted to ${\cal C}$ (called
\emph{explicit} equation-free computation in \cite{Marschler2014}).

When using equation-free methods one faces several challenges, both
analytical and in terms of implementation.  First, as the slow
manifold ${\cal C}$ cannot be assumed to be known to the user, the
method cannot assume that the user provides a lifting operator $\lift$
that maps onto ${\cal C}$. This leads to initial fast transients in
the trajectory that will also change the supposedly slow variables,
unless the stable fiber projection $g:\R^D\to{\cal C}$ keeps the
restriction constant (the criterion would be
$\restrict\circ g\circ \lift\approx \id$). Since the projection $g$
cannot be assumed to be known either, this implies that an unknown
nonlinear transformation is applied to the variables in $\dom\lift$
before the slow dynamics start. A detailed illustration of this
problem is given in \cref{fig:slowfast-parallel} and its
description in \cref{sec:state}.

Second, the justification for equation-free methods relies on the
stronger results for classical attracting invariant manifolds of
deterministic systems (including persistence of the slow manifold for
finite time-scale separation and its shadowing properties via stable
fiber projection). However, the methods are commonly applied to
stochastic or deterministic chaotic systems with time scale
separation, for which convergence results are weaker. In the
stochastic case the microscopic time stepper $M(t;\cdot)$ applies to
densities not single trajectories.  Finally, for applications with
stochastic microscopic systems the additional difficulty of low
computational accuracy in the evaluation of $M(t;\cdot)$ and possibly
$\lift$ and $\restrict$ may impose practical limitations.

This paper addresses the first challenge, the unknown slow manifold
and fiber projection. It proves convergence of the \emph{implicit}
approximation $y=\Phi_\tskip(\deltat;x)$ for the slow flow
$\Phi_*(\deltat;x)$, given by the solution $y$ of the $d$-dimensional
nonlinear system
\begin{displaymath}
  P(\tskip;y)=P(\tskip+\deltat;x)
\end{displaymath}
for sufficiently large \emph{healing time} $\tskip$ and a fixed finite
time scale separation for the scenario of an attracting
$d$-dimensional invariant manifold in $\R^D$ (strong reduction results
are available in this scenario). We also give a demonstration how the
implicit equation-free formulation behaves when it is applied to
moments of distributions in a stochastic system. Starting from this
demonstration, we outline in our subsequent discussion how convergence
statements for stochastic systems may have to be formulated.

\paragraph{Applications and recent practical improvements}
A motivation for using the equa\-tion-free framework is that it extends
methods which are otherwise only applicable to low-dimensional
dynamical systems directly to simulations of high-dimensional complex
systems. Classical applications of equation-free methods were
\emph{macroscopic bifurcation analysis} for microscopic simulations in
chemical engineering (see \cite{KS09} for a review). Recently similar
analysis was performed on stochastic network models of neurons
\cite{avitabile2016,MarschlerFaust2014} or disease spread \cite{GK08},
or on agent-based models in ecology \cite{TLS16} and social sciences
(for example, for consumer lock-in \cite{avitabile2014}, for
pedestrian flow \cite{Marschler2014,Marschler2015}, or for trading
\cite{SGK12}). Another example for a high-level task accessible via
equation-free methods is control design \cite{SMK04,SGK12}.

Recent modifications and improvements to equation-free methods in
multi-particle or agent-based simulations are variance reduction
\cite{papavasiliou2007variance,avitabile2014}, restriction of
computations to patches in space
\cite{SKR06,SRK10,liu2015acceleration} (for which a-priori error
estimates can be proven \cite{SKR06,SRK10}), and data-driven selection
of the slow variables using diffusion maps
\cite{coifman2008,Marschler2014}. Debrabant \emph{et al}.
\cite{debrabant2017} construct an acceleration scheme for Monte-Carlo
simulations of high-dimensional SDEs based on moments of densities
(the macroscopic variables), and prove its convergence as the number
of moments goes to infinity.

\section{Current state of analysis}
\label{sec:state}
\begin{figure}[ht]
  \centering
  \includegraphics[scale=0.7]{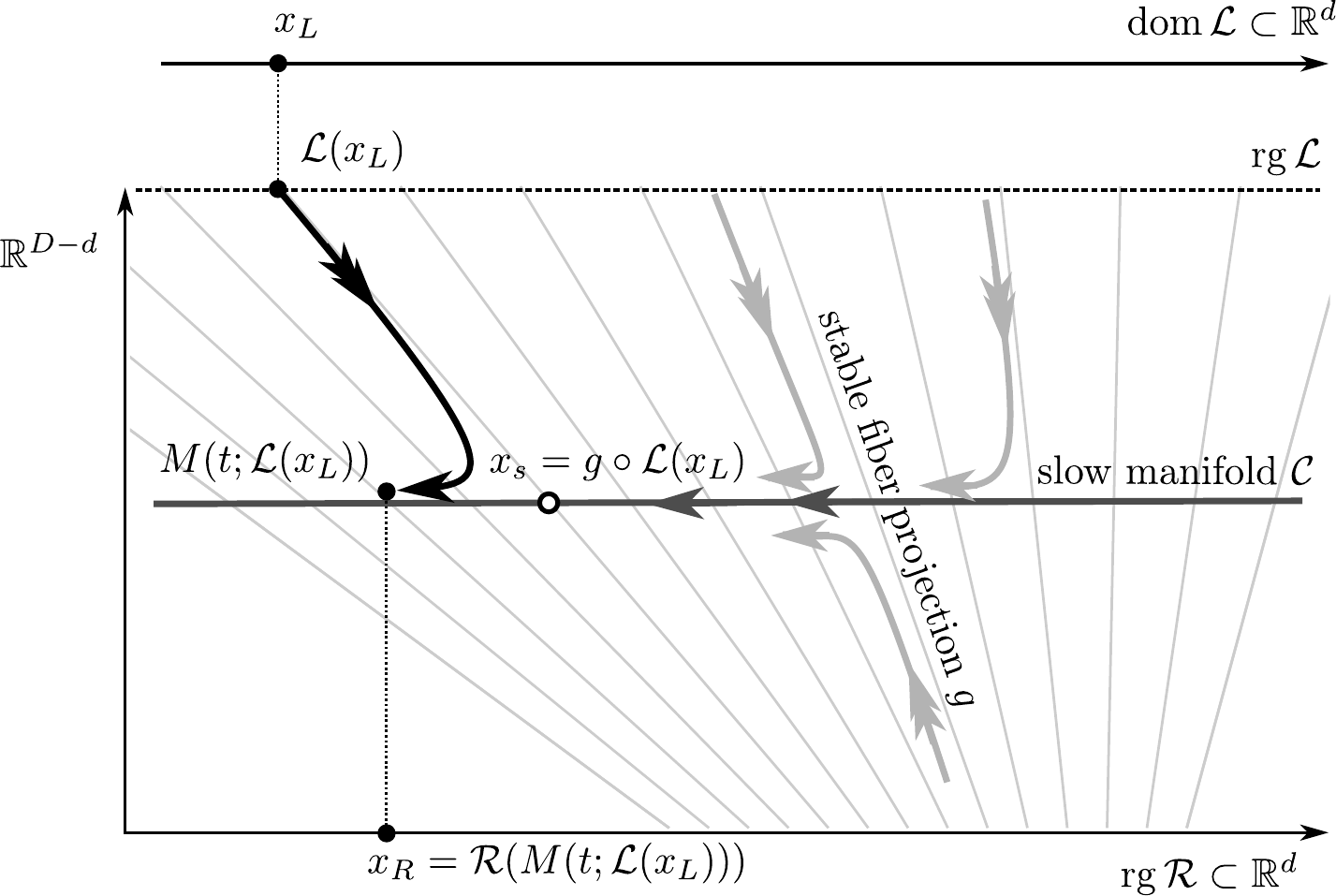}
  \caption{Geometry of lift-evolve-restrict map near slow
    manifold\textup{:} macroscopic value $x_L\in\R^d$ gets lifted to
    $\lift(x_L)$, then evolved to $M(t;\lift(x_L))$, then restricted
    to $x_R=\restrict(M(t;\lift(x_L))\in\R^d$. The aim is to
    approximate the true slow flow
    $x_L\mapsto (g\circ\lift)^{-1}\circ M(t;g(\lift(x_L)))$, which
    involves the unknown fiber projection $g$, using the map
    $\restrict\circ M(t;\cdot)\circ\lift$, and assuming invertibility
    of $g\circ\lift:\R^d\mapsto{\cal C}$.  For this sketch, $d=1$ and
    $D=2$.}
  \label{fig:slowfast-parallel}
\end{figure}

\paragraph{Geometry of the idealized case of an attracting slow manifold}
Analysis of the equation-free framework (based on
lift-evolve-restrict) is still ongoing. Convergence analysis with
general a-priori error estimates has been performed mostly for the
idealized case where the $D$-dimensional microscopic problem has a
$d$-dimensional attracting invariant slow manifold ${\cal C}$, which
is rarely encountered in the practical applications listed
above. Exceptions are, for example, a study of bursting neurons
\cite{BK16} and the application of implicit equation-free computations
to generalize an algorithm for growing stable manifolds of fixed
points of two-dimensional maps a delay-differential equation with an unknown
two-dimensional slow manifold  \cite{quinn2018effects}. Even for this
idealized case one faces the geometric difficulty illustrated in
\cref{fig:slowfast-parallel}.
The geometry shows an example scenario where the microscopic system is
two-dimensional, and the slow manifold ${\cal C}$ is horizontal (and,
thus, the slow motion is purely horizontal, drifting to the
left). Here we choose a lifting $\lift$ that maps also onto a
horizontal line $\rg \lift$.  However, $\rg\lift$ is at a distance to
${\cal C}$, because the precise location of ${\cal C}$ is in practice
unknown. The restriction $\restrict$ is the horizontal component of
any point $u\in\R^2$. The spaces $\dom \lift$ and $\rg\restrict$ (both
one-dimensional) are drawn separately for clarity in
\cref{fig:slowfast-parallel}, but they may be identical in
examples. The fast motion of $M(t;\cdot)$ is not perfectly vertical,
but has a significant horizontal component.
\Cref{fig:slowfast-parallel} also shows how the map $P$ acts on
a typical point $x_L$, showing its image $\lift(x_L)$, the result of
the evolution, $M(t;\lift(x_L))$, and the result of the restriction
$P(t;x_L)=x_R=\restrict(M(t;\lift(x_L)))$.

The point $x_s\in{\cal C}$ in \cref{fig:slowfast-parallel} is
defined as the unique point $x_s$ on ${\cal C}$ such that
$M(t;\lift(x_L))-M(t;x_s)$ converges at an exponential rate $\dtr$
that is larger than the maximal rate of contraction $\dtan$ tangential
to ${\cal C}$ (which is horizontal). As mentioned in the introduction
as shadowing, this mapping is defined for every point $u$ in the
neighborhood of ${\cal C}$: for every $u$ near ${\cal C}$ there exists
a point $g(u)\in{\cal C}$ such that
$M(t;u)-M(t;g(u))\sim\exp(-\dtr t)$ (in the illustration
$u=\lift(x_L)$, $g(u)=x_s$).  This point $g(u)$ is called the stable
fiber projection of $u$. The map $g$ is known to have the same
regularity as ${\cal C}$ \cite{Fenichel1979,HPS77}. For $\dtan\ll\dtr$
the map can be expanded in orders of $\epsilon=\dtan/\dtr$. The thin
grey lines (called stable fibers or isochrones) in
\cref{fig:slowfast-parallel} indicate how points in the plane are
projected onto ${\cal C}$ under the nonlinear projection $g$ for the
illustrative example.

\Cref{fig:slowfast-parallel} makes clear that the dynamics of
the map $P(t;x)$ is qualitatively different from the dynamics of
$M(t;\cdot)$ restricted to ${\cal C}$, $M(t;\cdot)\vert_{\cal C}$. For
the particular geometry shown in \cref{fig:slowfast-parallel}
$P(t;\cdot)$ has a unique stable fixed point if the horizontal
attraction/expansion rate $\dtan$ of $M(t;\cdot)$ on ${\cal C}$ is
sufficiently small compared to the attraction rate $\dtr$ transversal
to ${\cal C}$. This fixed point is nearly independent of the dynamics
of $M(t;\cdot)$ on ${\cal C}$.

More generally, if the lifting operator $\lift$ does not map $x_L$
into the low-dimensional slow manifold ${\cal C}$ then the initial
part of the trajectory $t\mapsto M(t;\lift(x_L))$, which is computed
as part of the lift-evolve-restrict map $P$, is a rapidly changing
transient toward the slow manifold ${\cal C}$, which will generically
also change the resulting $x_R$.


In the limit of infinite time-scale separation (that is, the
derivative of $M$ with respect to time, $\partial_1M(t;,u)$, goes to
$0$ for $u\in{\cal C}$) the dynamics of the lift-evolve-restrict map
$P$ is a small perturbation of the map $\restrict\circ g\circ \lift$.
Unless this limit map equals the identity, $P(t;\cdot)$ cannot be a
good approximation of the slow flow along the manifold ${\cal C}$. 
Using $x$ in the domain of the lifting $\lift$ and the
map $g\circ \lift:\dom\lift\mapsto{\cal C}$ onto the manifold ${\cal
  C}$ as the coordinate map, the slow flow $\Phi_*$ has the form
\begin{equation}\label{intro:phistar}
  \begin{aligned}
    \Phi_*(\deltat;\cdot)&=(g\circ \lift)^{-1} \circ M(\deltat;\cdot)
    \circ (g \circ \lift)\mbox{,}&&\mbox{or implicitly defined by}\\
    \Phi_*:&\ \R\times\dom\lift\ni (\deltat,x)\mapsto y_*\in\dom\lift\mbox{,}
    &&\mbox{where $y_*$ solves}\\
    \restrict(g(\lift(y_*)))&=\restrict(M(\deltat;g(\lift(x))))
  \end{aligned}
\end{equation}
(using the notation $(\cdot)^{-1}$ for the inverse map).  This
definition is not directly computable since the nonlinear
projection $g$ is unknown in general.

Feasible approaches to construct an accurate approximation of
$M(t;\cdot)$ restricted to ${\cal C}$ are \emph{constrained runs}, as
discussed by Gear, Zagaris \emph{et al}.
\cite{Gear2005,Zagaris2009,Zagaris2012}, or the introduction of a
healing time $\tskip$. The latter approach is studied in this
paper.


\subparagraph{Constrained runs} The approach of
\cite{Gear2005,Zagaris2009,Zagaris2012} to ensuring that
$\restrict\circ g\circ \lift$ is close to the identity is to enforce
that the lifting $\lift$ maps onto the manifold ${\cal C}$ with
sufficient accuracy for all $x$ in its domain. Usually, this requires
an additional scheme involving the iterative application of $\lift$
and $M$; see \cite{Gear2005,Zagaris2009,Zagaris2012}. The a-priori
error estimates prove that the lift-evolve-restrict scheme with these
additional iterations has an error of order $(\dtan/\dtr)^m$ if the
constrained runs scheme is of order $m$, where $\dtan$ is the
attraction/repulsion time scale tangential to the slow manifold
${\cal C}$ and $\dtr$ is the transversal attraction rate. The ratio
$\dtan/\dtr$ measures the time scale separation. It is assumed to be
small when applying constrained runs (and called $\epsilon$), and
$O(\epsilon^m)$ convergence is proven in
\cite{Gear2005,Zagaris2009,Zagaris2012} in the limit $\epsilon\to0$.
This limit will not be required in our proof, later on.

\subparagraph{Implicit formulation with healing time} A second,
alternative, approach is to introduce a \emph{healing time} $\tskip$,
exploiting that $M$ attracts along the fibers
\cite{Kevrekidis2010,Barkley2006}.
\begin{figure}[ht]
  \centering
  \includegraphics[scale=0.7]{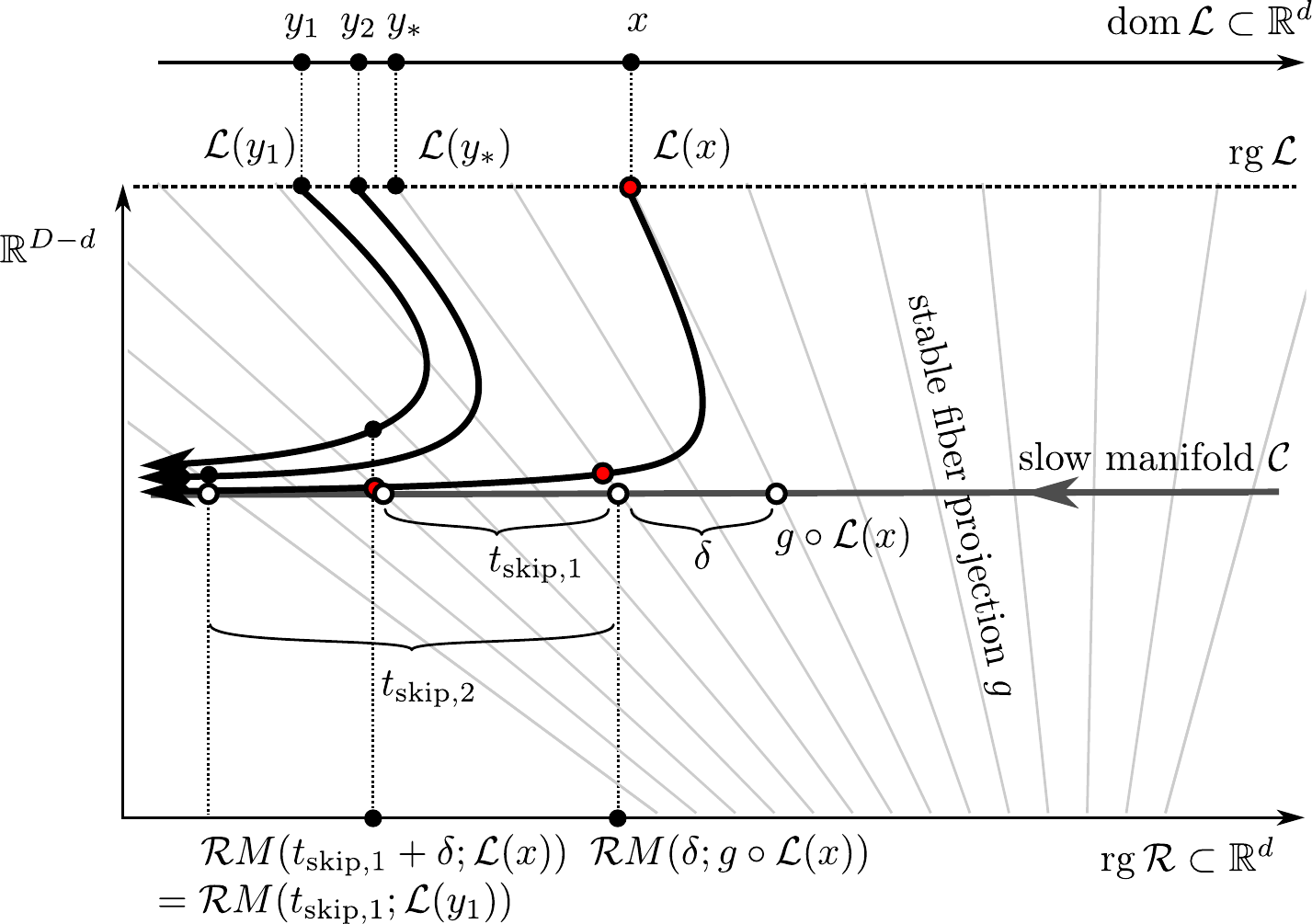}
  \caption{Trajectories involved in the implicit definition of $y_*$
    and $y$ using the scenario of
    \cref{fig:slowfast-parallel}\textup{:} $y_*$ is the image of $x$
    under the true slow flow, $y_*=\Phi_*(\deltat;x)$; $y_j$ are
    approximations for different $\tskip$,
    $y_j=\Phi_{\tskip,j}(\deltat;x)$.}
  \label{fig:slowfast-convsketch}
\end{figure}
Marschler \emph{et al}.
\cite{Marschler2014b} show that the healing time $\tskip$ can be
motivated by introducing an additional shift $M(\tskip;\cdot)$ and its
inverse into \cref{intro:phistar} (note that $M(\tskip;\cdot)$ is
invertible on the slow manifold ${\cal C})$:
\begin{equation}
  \label{eq:phistar:tskip}
    y_*=\Phi_*(\deltat;x)=(g\circ \lift)^{-1} 
    \circ M(\tskip;\cdot)^{-1}\circ M(\deltat+\tskip;\cdot) \circ g
  \circ \lift(x)\mbox{.}
\end{equation} 
Removing the inverses in \cref{eq:phistar:tskip} leads to an implicit
equation for $y_*=\Phi_*(\delta;x)$ with the healing time $\tskip$ as
an additional parameter:
\begin{equation}
  \label{eq:phistar:tskip:impl}
  \restrict\circ M(\tskip;\cdot)\circ g\circ\lift(y_*)=\restrict\circ 
  M(\deltat+\tskip;\cdot)\circ g\circ \lift(x)
\end{equation}
In \cref{eq:phistar:tskip:impl} the parameter $\tskip$ has no effect
since $M(\tskip;\cdot)$ is invertible on the slow manifold.  However,
the difference $M(\tskip;\cdot)\circ g-M(\tskip;\cdot)$ decreases with
$\tskip$ (at rate $\sim\exp(-\dtr\tskip)$).  In
\cref{fig:slowfast-convsketch} the distance between points along the
trajectory starting from $\lift(x)$ (in red) and their projections
$g\circ \lift(x)$ (white) illustrates this convergence. Thus, we may
approximate $M(\tskip;\cdot)\circ g$ by $M(\tskip;\cdot)$ in
\cref{eq:phistar:tskip:impl}. This results in a computable
approximation $y_\tskip=\Phi_\tskip(\deltat;x)$ of $y_*$, given
implicitly by the equation
\begin{equation}\label{eq:intro:phidef}
  \restrict(M(\tskip;\lift(y_\tskip)))=
  \restrict(M(\deltat+\tskip;\lift(x)))\mbox{.}
\end{equation}
\Cref{fig:slowfast-convsketch} illustrates the effect of
increasing healing time $\tskip$ in the scenario introduced in
\cref{fig:slowfast-parallel}. The points $y_1$ and $y_2$ are the
solutions of \cref{eq:intro:phidef} for two different healing times
$t_{\mathrm{skip},1}<t_{\mathrm{skip},2}$.
Equation~\eqref{eq:intro:phidef} means that the points $y_j$ are
defined as those elements of $\dom\lift$ for which the trajectory
starting from $\lift(y_j)$ has the same horizontal component
(restriction $\restrict$) as $M(t_{\mathrm{skip},j}+\deltat;\lift(x))$
after time $t_{\mathrm{skip},j}$.

The implicit approach was analyzed and
illustrated in a traffic model in \cite{Marschler2014b} and will also
be studied in this paper. Vandekerckhove \emph{et al}.
\cite{Vandekerckhove2011} proposed and demonstrated a similar approach, but applied
the healing time backward in time by fixing the image of the
restriction: they solve $x=\restrict(M(\tskip;\lift(x_b)))$ for $x_b$
first and then set $y=\restrict(M(\deltat+\tskip;\lift(x_b)))$. This
gives an (approximate) representation $\Phi_*^{\restrict}$ of the slow
flow in the coordinates on the image of the restriction $\restrict$:
$\Phi_*^{\restrict}(\deltat;x)=\restrict\circ M(\deltat;\cdot)\circ
[\restrict\vert_{\cal C}]^{-1}(x)$. 

The coordinates for the flow on the slow manifold ${\cal C}$ are
somewhat arbitrary as the difference between the expressions used by
Vandekerckhove \emph{et al}.  \cite{Vandekerckhove2011} and the
implicit expression \cref{eq:intro:phidef} for $\Phi_\tskip$
shows. For the coordinates in the space $\dom\lift$ the diffeomorphism
between $\dom\lift$ and ${\cal C}$ is $g\circ\lift$, where $g$ is the
stable fiber projection, as implied by \cref{eq:phistar:tskip}. The
diffeomorphism can be approximately computed by solving
$\restrict(M(2\tskip;\lift (x_g)))=\restrict(M(\tskip;\lift (x)))$ for
$x_g$ and then using $M(\tskip;\lift(x_g))$ as the approximation for
$[g\circ\lift](x)$. The approximate diffeomorphism for the expression
of Vandekerckhove \emph{et al}.  \cite{Vandekerckhove2011} is
$\left[\restrict\vert_{\cal C}\right]^{-1}:x\mapsto
M(\tskip;\lift(x_b))$.

Marschler \emph{et al}.
\cite{Marschler2014b} proved that the approximation $y_\tskip$ is
exponentially accurate if $\dtan/\dtr\to0$: $\|y_\tskip-y_*\|\sim
\exp(-K\dtr/\dtan)$ (for some constant $K$ depending on $\tskip$). The
error estimates in \cite{Marschler2014b} require that
$\tskip\dtan/\dtr$ and $(\tskip+\deltat)\dtan/\dtr$ stay bounded from
above such that the convergence result is valid in the limit of
infinite time scale separation $\dtan/\dtr\to0$.
%
%
This means that the assumptions of \cite{Marschler2014b} are similar
to those required by schemes involving constrained runs
\cite{Gear2005,Zagaris2009,Zagaris2012}.  The analysis left open if
the error goes to zero for $\tskip\to\infty$ but the time scale
separation stays finite: $\dtan/\dtr\in(0,1)$.

Our paper will prove the general a-priori error
estimate that $\|y_\tskip-y_*\|\sim \exp((\dtan-\dtr)\tskip)$ for
$\tskip\to\infty$ and fixed $\dtan<\dtr$ under some genericity
conditions on $\restrict$ and $\lift$. It will also give a convergence
result for the derivatives of $y_\tskip$ with respect to its argument
$x$: $\|\partial^jy_\tskip-\partial^jy_*\|\sim
\exp(((2j+1)\dtan-\dtr)\tskip)$ if $(2j+1)\dtan<\dtr$.

\paragraph{Analysis beyond attracting manifolds in slow-fast systems}
As mentioned above, equation-free analysis based on
lift-evolve-restrict maps is more commonly applied to problems that
are assumed to have a fast subsystem, where the fast time scale
converges only in a statistical sense to a stationary measure
conditioned on the slow variables. In these cases the microscopic time
stepper $M(\delta;\cdot)$ operates on measures (or densities). It may
be approximated by Monte Carlo simulations on ensembles of initial
conditions. Barkley \emph{et al}. \cite{Barkley2006} investigated the
behaviour of the lift-evolve-restrict map
$P(\deltat;\cdot)=\restrict\circ M(\deltat;\cdot)\circ \lift$ where
the slow variables were leading moments (thus, $P$ was called moment
map in \cite{Barkley2006}) on prototype examples from the class of
stochastic problems. The simplest example from \cite{Barkley2006} is a
scalar stochastic differential equation (SDE), for which the evolution
of the probability distribution is governed by a (linear)
Fokker-Planck equation (FPE). Hence, the measure of time-scale
separation is the size of the spectral gap in the right-hand side of
the FPE. The analysis in \cite{Barkley2006} found that the dynamics of
the map $P$ was qualitatively different from the dynamics of the
underlying linear FPE. For example, $P$ was nonlinear and had several
coexisting fixed points for certain choices of time $\deltat$.

Our paper will demonstrate for two different lifting operators $\lift$
that the approximation $y_\tskip$, defined by \cref{eq:intro:phidef},
behaves exactly as predicted by our convergence theorem. In
particular, it preserves the metastability features and the linearity
of the flow generated by the FPE, thus, addressing the problems
highlighted in \cite{Barkley2006}.

\subsection{Outline of results}
\Cref{sec:slowfast} states the precise assumptions (time scale
separation for decay rates tangential and transversal to
the invariant manifold ${\cal C}$ ($\dtan<\dtr$) and transversality of
$\restrict$ and $\lift$) for exponential convergence:
\begin{equation}\label{eq:outline:conv}
  \partial^jy_\tskip-\partial^jy_*\sim\exp(((2j+1)\dtan-\dtr)\tskip)
  \quad\mbox{for $\tskip\to\infty$}
\end{equation}
(using the convention that $\partial^0y=y$ and assuming that the
derivatives up to order $j+1$ exist). Estimate~\cref{eq:outline:conv}
predicts that convergence in $\tskip$ is slower for derivatives of
higher order.  \Cref{sec:exampl-mich-ment} demonstrates the
convergence rates in $\tskip$ for $y_\tskip$ and its first two
derivatives with respect to $x$ for a singularly perturbed ODE
modelling the Michaelis-Menten kinetics (which was also used by
\cite{Gear2005,Zagaris2009,Zagaris2012} for
illustration). \Cref{sec:sde-dynam} studies the evolution of
densities under a scalar SDE with a double-well potential drift term
also considered by Barkley \emph{et al}. \cite{Barkley2006}.  We
demonstrate global convergence of implicit equation-free methods for a
linear lifting $\llin$. We also demonstrate local convergence for the
nonlinear lifting $\lgauss$ used in \cite{Barkley2006}.

\Cref{sec:discussion} discusses differences between
observations of the behaviour in the SDE and the predictions from the
theoretical result. These are caused by the numerical errors in the
evaluations of lifting, evolution and restriction and their growth
along trajectories.

We conclude with an outlook on possible consequences of the results on
application of equation-free methods to Monte-Carlo simulations of
multi-particle or agent-based systems. One important observation is
that in some cases increasing the number of agents or particles does
\emph{not} increase the spectral gap (and, thus, the time scale
separation). Didactic examples where the finiteness of the spectral
gap is apparent are the dynamic networks as considered by Gross and
Kevrekidis \cite{GK08}. The slow system is an ODE derived from the
pair-wise interaction approximation, cutting off an infinite series of
ODEs of higher-order interaction terms. The spectral gap between
pair-wise interaction terms and triplet interaction terms is finite
even in the limit of infinitely large networks.

Thus, the results from \Cref{sec:slowfast} are potentially
applicable to equation-free analysis of stochastic multi-particle
systems, where distributions of microscopic initializations are
studied. This is in contrast to previous convergence results on
constrained runs \cite{Gear2005,Zagaris2009,Zagaris2012} and implicit
lifting \cite{Marschler2014b}, which only apply in the limit of
infinite time scale separation.

\section{Convergence in the case of finite time-scale
  separation}
\label{sec:slowfast}
We consider a smooth dynamical system
\begin{equation}
\label{eq:f}
  \dot u=f(u), \quad u\in\R^D\mbox{,}
\end{equation}
where $D$ is large. We assume that the flow $M$ generated by \cref{eq:f},
\begin{displaymath}
M:\R\times \R^D\to \R^D\mbox{,} \qquad (t;u)\mapsto M(t;u)
\end{displaymath}
has a $d$-dimensional compact relatively invariant manifold ${\cal C}$
(possibly with boundary). That is, trajectories $M(t;u)$ starting in
$u\in{\cal C}$ either stay in ${\cal C}$ for all times $t\in\R$, or
they stay in ${\cal C}$ until they cross the boundary $\partial{\cal C}$ of
${\cal C}$. We assume that ${\cal C}$ is at
least $k_{\max}$ times differentiable. For a point $u\in{\cal C}$, let
us denote by ${\cal N}(u)$ the $d$-dimensional tangent space to ${\cal
  C}$.  The following assumption states that attraction transversal to
the manifold ${\cal C}$ is faster than attraction or expansion
tangential to ${\cal C}$.
\begin{assumption}[Hyperbolicity --- Separation of time scales and transversal
  stability]\label{ass:timescales}
  There exists an open neighborhood ${\cal U}$ of the manifold ${\cal
    C}$, a \textup{(}possibly nonlinear\textup{)} projection $g:{\cal
    U}\mapsto{\cal C}$ \textup{(}the so-called \emph{stable fiber
    projection}\textup{)}, a pair of constants \textup{(}decay
  rates\textup{)} $0<\dtan<\dtr$, and a bound $C$ such that the
  following two conditions hold.
  \begin{enumerate}\item \textup{(tangential expansion/attraction rate)}\label{ass:hyptan}
    For all points $u\in{\cal C}$ on the manifold, all tangent vectors
    $v_1,\ldots,v_{k_{\max}}\in{\cal N}(u)$ and all $t\in\R$ with
    $M([0,t];u)\subset{\cal C}$
    \begin{align}
      \label{eq:mcebound}
      \|\partial_2^jM(t;u)[v_1,\ldots,v_{k_j}]\| &\leq
      C\exp(\dtan|t|)\|v_1\|\cdot\ldots\cdot\|v_j\|
    \end{align}
    for all $j\in\{1,\ldots,k_{\max}\}$.
  \item \textup{(Stability along transversal fiber projections)} For
    all $u\in{\cal U}$ and all $t>0$ with $M([0,t];g(u))\in{\cal C}$
    \begin{align}
      \label{eq:gepsconv}
      \|\partial_2^jM(t;u)-\partial_2^jM(t;g(u))\|&\leq 
      C\exp(-t\dtr)\|u-g(u)\|
    \end{align}
    for all $j\in\{0,\ldots,k_{\max}\}$.
\end{enumerate}
\end{assumption}
\noindent
In \cref{eq:mcebound} and \cref{eq:gepsconv} we use the convention
that $\partial^j_kM$ is the $j$th-order partial derivative of $M$ with
respect to its $k$th argument, and that $\partial_2^0M$ ($j=0$) is the
flow $M$ itself. The norm on the left side of \cref{eq:gepsconv} is
the usual operator norm for the multi-linear operators
$\partial_2^jM(t,\cdot)$. The constants $C$, $\dtr$ and $\dtan$ are
assumed to be independent of the point $u$ and the time
$t$. Assumption \cref{eq:mcebound} is also made for negative times
$t$ (using the convention that $M([0,t];u)$ means $M([t,0];u)$ for
$t<0$) such that it is also an assumption about the inverse of $M$,
when restricted to ${\cal C}$: $M(-t,\cdot)=M^{-1}(t,\cdot)$. The
constant $\dtr$ is the decay rate toward the manifold ${\cal C}$, the
constant $\dtan$ is the rate of attraction and expansion along the
flow restricted to ${\cal C}$. The main requirement of
\cref{ass:timescales} is that $\dtr>\dtan$.

\paragraph{Transversality of restriction and lifting}
Second, we assume basic compatibility between 
\begin{equation}
\begin{aligned}
\label{eq:lift_restrict}
  \restrict&: {\cal U}\subset\R^D\to \R^d &&\mbox{the \emph{restriction} operator,}\\
  \lift&: \dom\lift\subset\R^d\to\R^D &&\mbox{the \emph{lifting} operator,}
\end{aligned}
\end{equation}
and the invariant manifold ${\cal C}$: the lifting $\lift$ should map
into the neighborhood ${\cal U}$ of ${\cal C}$ in which the stable
fiber projection $g$ is defined, and the restriction $\restrict$
should be defined on the projection $g$ of the image of $\lift$ along
the stable fibers:
\begin{align*}
  \lift(\dom\lift)&\subset{\cal U}\mbox{,} &
  g(\lift(\dom\lift))&\subset \dom{\cal R}\cap{\cal C}\mbox{,}
\end{align*}
In addition to these compatibility conditions, we impose the
following two transversality conditions on lifting $\lift$ and
restriction $\restrict$.
\begin{assumption}[Transversality of $\restrict$ and
  $\lift$]\label{ass:transversality}
  \begin{enumerate}
  \item\label{ass:ltrans}the projection $g$ is a diffeomorphism
    between $\rg\lift=\lift(\dom\lift)$ and ${\cal C}$. In particular, for all
    $x\in\dom\lift\subset\R^d$
    \begin{displaymath}
      \rank \frac{\partial}{\partial x}\left[g(\lift(x))\right]=
      \rank\left[\partial g(\lift(x))\circ \partial \lift(x)\right]=d\mbox{.}
    \end{displaymath}
  \item\label{ass:rtrans} The map $\restrict$, restricted to ${\cal
      C}$, is a diffeomorphism between ${\cal C}$ and $\R^d$. In
    particular, for all $u\in{\cal C}$ \textup{(}${\cal N}(u)$ is the tangent
    space to ${\cal C}$ in $u$\textup{)}
    \begin{displaymath}
      \dim\partial\restrict(u){\cal N}(u)=d\mbox{.}
    \end{displaymath}
  \end{enumerate}
\end{assumption}

\paragraph{Coordinates on the slow manifold ${\cal C}$}
The maps $\restrict$ and $\lift$ create two natural ways to define
local coordinate representations on the invariant manifold ${\cal C}$,
one by a map from $\dom\lift$ to ${\cal C}$, one by a map from ${\cal
  C}$ to $\rg\restrict$.  For our presentation we choose the
representation in $\dom\lift$ coordinates:
\begin{align*}
  g\circ \lift&:\dom\lift\subset\R^d\mapsto{\cal  C}\subset \R^D\mbox{,}&
  x\mapsto g(\lift(x))
  \mbox{.}
\end{align*}
The inverse of $g\circ\lift$ is defined implicitly.  Assume that
$u_0=g(\lift(x_0))$ for some $x_0\in\dom\lift$. Then for $u\in{\cal
  C}$ near $u_0$ the pre-image $x=(g\circ\lift)^{-1}(u)$ is found by
solving $\restrict(u)=\restrict(g(\lift(x)))$ for $x\approx x_0$,
which has a locally unique solution by
\cref{ass:transversality}. 

We can represent the flow $M$ on ${\cal C}$ as a flow in $\dom\lift$,
denoting it by $\Phi_*$:
\begin{align}\label{eq:phistardef1}
  \Phi_*&:\R\times\dom\lift\mapsto\dom\lift\mbox{,}
  &\Phi_*(\delta;x)_{\phantom{r}}&=[(g\circ\lift)^{-1}\circ M(\deltat;\cdot)\circ
  g\circ\lift](x):=y\mbox{}
\end{align}
(for $\delta\in\R$ and $x\in\dom\lift$), where
$y$ is given implicitly as solution of a
$d$-dimensional system of nonlinear equations
\begin{align}
  \restrict(g(\lift(y)))&=
  \restrict(M(\delta;g(\lift(x))))\mbox{.}
  \label{eq:phistarts0}
\end{align}
\Cref{ass:transversality} on transversality implies that
$\Phi_*$ is well defined for small $\delta$ (since $y=x$ is a regular
solution of \cref{eq:phistarts0} at $\delta=0$). For larger $\delta$,
one can break down the solution into smaller steps by increasing
$\delta$ gradually from $0$ and tracking the curve $y(\delta)$ of
solutions of \cref{eq:phistarts0}, which is well parametrized by
$\delta$ in every point by \cref{ass:transversality}. If
$\dom\lift$ is simply connected then this continuation approach makes
the implicit solution $y$ used in the definition of $\Phi_*$
unique. 
Let us define the map
\begin{align}\label{eq:pstardef}
  P_*&:\R\times\dom\lift\ni (t,x)\mapsto\restrict(M(t;g(\lift(x))))\in\R^d\mbox{.}
\end{align}
This map $P_*$ is well defined and invertible for all $t\in\R$ and
$x\in\dom\lift$ for which the trajectory $s\mapsto M(s;g(\lift(x)))$
stays in ${\cal C}$ for all $s$ between $0$ and $t$. The implicit
definition \cref{eq:phistardef1} of the flow $\Phi_*$ on ${\cal C}$
has the following form when expressed with the help of this map $P_*$
on $\dom\lift$:
 \begin{align}
  \label{eq:phistar0p}
  y=\Phi_*(\delta;x) \mbox{\quad if\quad } P_*(0;y)=P_*(\delta;x)\mbox{.}
\end{align}
Since the flow $M(\delta;\cdot)$ is a diffeomorphism on ${\cal C}$, we
can replace the times $0$ and $\delta$ in the above implicit
definition with $\tskip$ and $\tskip+\delta$ for an arbitrary
so-called \emph{healing time} $\tskip\in\R$ (as long as
$M([0,\tskip];g(\lift(x)))\subset{\cal C}$). So, equivalent to
\cref{eq:phistar0p}, we have for $\tskip>0$ with $M([0,\tskip];g(\lift(x)))\subset{\cal C}$
 \begin{align}
  \label{eq:phistarp}
  y=\Phi_*(\delta;x)\mbox{\quad if\quad } P_*(\tskip;y)=P_*(\tskip+\delta;x)\mbox{.}
\end{align}
\paragraph{Convergence Theorem for implicit equation-free
  computations with finite time-scale separation}
The stable fiber projection $g$ (which is part of the definition of
$P_*$) is not known in most practical applications. Thus, implicit
equation-free computations use the explicit macroscopic time-$t$ map
$P$ instead of $P_*$ in the equation defining $y$
in \cref{eq:phistarp}:
 \begin{align}
  \label{eq:pdef}
  P&:[0,\infty)\times \dom\lift\ni(t,x)\mapsto
  \restrict(M(t;\lift(x)))\in\rg\restrict
\end{align}
such that we may define the approximate flow map
 \begin{align}
  \label{eq:phip}
  \Phi_\tskip:\R\times\dom\lift\ni(\deltat,x)\mapsto y\in\dom\lift
  \mbox{,\ where $y$ solves\ }P(\tskip;y)=P(\tskip+\delta;x)\mbox{}
\end{align}
implicitly in a similar way to \cref{eq:phistarp}.
Our general convergence theorem, the following \cref{thm:conv},
states that $\Phi_\tskip$ is well defined for large $\tskip$ (that is,
the equation in \cref{eq:phip}, defining $\Phi_\tskip$ implicitly,
has a locally unique solution), and that $\partial^j\Phi_\tskip$ is an
approximation of $\partial^j\Phi_*$ of order
$\exp(((2j+1)\dtan-\dtr)\tskip)$ (including $j=0$ for the map
$\Phi_\tskip$).
\begin{theorem}[Convergence of approximate flow map at finite
  time-scale separation]\label{thm:conv}\ \\
  Let us assume that the microscopic flow $M$ satisfies
  \cref{ass:timescales} on time-scale separation, and that
  the maps $\restrict$ and $\lift$ satisfy
  \cref{ass:transversality} on transversality.  

  Let $\delta_{\max}>0$ and $x\in\dom\lift$ be arbitrary. Let us also
  assume that $x\in\dom\lift$ maps to a point under $g\circ\lift$
  that keeps a positive distance from the boundary $\partial {\cal C}$
  of ${\cal C}$ for all times $t\geq-\delta_{\max}$ under $M$. That
  is,
  \begin{equation}\label{eq:doml:inflow}
    \dist(M(t;g(\lift(x))),\partial {\cal C})\geq c_\partial
    \mbox{\quad for all $t\geq-\delta_{\max}$ and some given $c_\partial>0$.}
  \end{equation}
     Then there exists a $t_0\geq \delta_{\max}$ such that
  $y=\Phi_\tskip(\deltat;x)$ is well defined by
 \cref{eq:phip} for all
  $\deltat\in[-\delta_{\max},\delta_{\max}]$ and $\tskip>t_0$. The estimate
  \begin{equation}
    \label{eq:thm:accuracy}
    \|\partial_2^j\Phi_\tskip(\deltat;x)-\partial_2^j\Phi_*(\deltat;x)\|\leq 
    C\exp(((2j+1)\dtan-\dtr)\tskip)
  \end{equation}
  holds for all orders $j\in\{0,\ldots,k_{\max}-1\}$ satisfying
  $(2j+1)\dtan<\dtr$. The constant $C$ depends on
   $\delta_{\max}$ and $x$, but not on $\tskip$.
\end{theorem}
Assumption \cref{eq:doml:inflow} in \cref{thm:conv} is made to
permit arbitrarily large $\tskip$ while still having
\cref{ass:timescales} and \cref{ass:transversality} uniformly satisfied. If one
considers $x\in\dom\lift$ for which the trajectory $t\mapsto
M(t;g(\lift(x)))$ leaves ${\cal C}$ (by crossing the boundary
$\partial {\cal C}$) then one has to put restrictions on $\deltat$ and
$\tskip$ to avoid crossing $\partial {\cal C}$. The theorem permits
negative integration times $\deltat$ shorter than $t_0\leq \tskip$ and
positive integration times larger than $\tskip$ as long as the factor
$\exp(\dtan|\deltat|)$ is of order $1$. Since $1/\dtan$ and $1/\dtr$
are the time scales of the dynamics inside the invariant manifold
${\cal C}$ and transversal to it, the theorem covers time steps of
length $\deltat$ of order $1$ in the slow ($1/\dtan$) time scale. The
statement of \cref{thm:conv} does not require that the
time-scale separation $\dtan/\dtr$ goes to zero for convergence of the
approximate map. It only requires that $\dtan<\dtr$, where $\dtr$ is
the attraction rate along fibers (see \cref{eq:gepsconv}) and $\dtan$
is the attraction and expansion rate tangential to the invariant
manifold ${\cal C}$.
Since the constant $C$ in \cref{eq:thm:accuracy} is independent of
$\tskip$ it can be chosen uniformly for compact domains $\dom\lift$.
\paragraph{Outline of proof of \cref{thm:conv}}
\emph{Existence and
error of $\Phi_\tskip$}: For the proof of \cref{thm:conv} we
have to analyze the difference between the two defining equations for
the approximate solution $y$ and the true solution
$y_*=\Phi_*(\deltat;x)$, both depending on $x$ as a parameter:
\begin{align}\label{eq:out:perturbed}
\restrict(M(\tskip;\lift(y)))&=
  \restrict(M(\tskip+\deltat;\lift(x)))\mbox{,}\\
    \label{eq:out:exact}
    \restrict(M(\tskip;g(\lift(y_*))))&=
    \restrict(M(\tskip+\deltat;g(\lift(x))))
    \mbox{.}
\end{align}
Rearranging the difference between \cref{eq:out:perturbed} and
\cref{eq:out:exact}, we obtain an implicit fixed-point
problem for $y$ (recall that $P_*(t;\cdot)=\restrict\circ
M(t;\cdot)\circ g\circ \lift$ and $P(t;\cdot)=\restrict\circ
M(t;\cdot)\circ \lift$):
\begin{equation}\label{eq:bfix}
  P_{*}(\tskip;y)=P_{*}(\tskip;y_*)+[P_{*}(\tskip;y)-P(\tskip;y)]+
  [P(\tskip+\deltat;x)-P_{*}(\tskip+\deltat;x)]\mbox{.}
\end{equation}
The norms of both terms in square brackets, $P_*(\tskip;y)-P(\tskip;y)$ and
$P(\tskip+\deltat;x)-P_*(\tskip+\deltat;x)$, are of order
$\exp(-\dtr\tskip)$ by \cref{ass:timescales},
equation~\cref{eq:gepsconv} (transversal stability with rate $\dtr$
of ${\cal C}$). For the same reason, the Lipschitz constant of
$P_*(\tskip;y)-P(\tskip;y)$ is of order $\exp(-\dtr\tskip)$, too. By
\cref{ass:timescales}, equation~\cref{eq:mcebound}
(tangential decay rate inside the manifold is less than $\dtan$), and
\cref{ass:transversality} on transversality of $\lift$ and
$\restrict$, the inverse of $P_*(\tskip;\cdot)$ has a local Lipschitz
constant of order $\exp(\dtan\tskip)$ near $y_*$. These two facts
enable us to apply the Banach Contraction Mapping Principle to
\cref{eq:bfix} to obtain a unique solution $y\approx y_*$ for large
$\tskip$. More precisely, $y-y_*$ is of order
$\exp((\dtan-\dtr)\tskip)$.

\emph{Inductive proof of error estimate for derivatives}: We differentiate
\cref{eq:bfix} with respect to $x$ in its fixed point $y(\tskip;x)$ up to
$j$ times and then re-arrange the resulting equation for the
$j$th-order derivatives of $y$ and $y_*$ into the form
\begin{equation}
  \label{eq:error:recursion}
  \partial P_*(y)\left[\partial^jy-\partial^jy_*\right]=
  \left[\partial P_*(y)-\partial P_*(y_*)\right]\partial^jy_*+r\mbox{.}
\end{equation}
In \cref{eq:error:recursion} we abbreviated
$\partial_2P_*(\tskip;\cdot)=\partial P_*(\cdot)$,
$\partial^jy=\partial^j_2y(\tskip;x)$ and dropped the argument $x$ from
$y_*$ and the arguments $\tskip$ and $x$ from $y$. The remainder $r$
is less than $C\exp(((2j-1)\dtan-\dtr)\tskip)$ for some constant $C$
by induction hypothesis.  The implicit expression
\cref{eq:error:recursion} for $\partial^jy-\partial^jy_*$ shows why
errors in derivatives of the solution can grow for increasing $\tskip$
and insufficient time scale separation: the norms of $\partial
P_*(y)-\partial P_*(y_*)$ and of
$\left[\partial P_*(y)\right]^{-1}$ are of order
$\exp(\dtan\tskip)$ due to \cref{eq:mcebound}.  
The details of the proof are given in
\cref{sec:app:conv}.

\section{Example: Michaelis-Menten kinetics}
\label{sec:exampl-mich-ment}
To illustrate the consequences of error estimate \cref{eq:thm:accuracy},
we look at a model for Michaelis-Menten kinetics with explicit time
scale separation as studied in
\cite{OMalley1991,Gear2005,Zagaris2009,Zagaris2012}.  The system is
given in $\R^D$ with $D=2$ as
\begin{equation}
  \begin{aligned}
    \label{eq:mmdyn}
    \dot x           &= \epsilon\left[-x + (x+\kappa - \lambda)y\right]\mbox{,} & 
    \dot y &= x  - (x+\kappa)y\mbox{,}
  \end{aligned}
\end{equation}
where $x\in\R$ is the slow variable, $y\in\R$ is the fast variable,
and $\epsilon$ measures the time scale separation. The parameters
$\kappa=1, \lambda=0.5$ and $\epsilon=0.01$ are kept fixed throughout
this section.  

\begin{figure}[t]
  \centering
  \newlength{\mmwidth}
  \setlength{\mmwidth}{0.42\textwidth}
 \subfigure[Phase space geometry of \cref{eq:mmdyn}]{\label{fig:mm_dyn_ph_s}\includegraphics[width=\mmwidth]{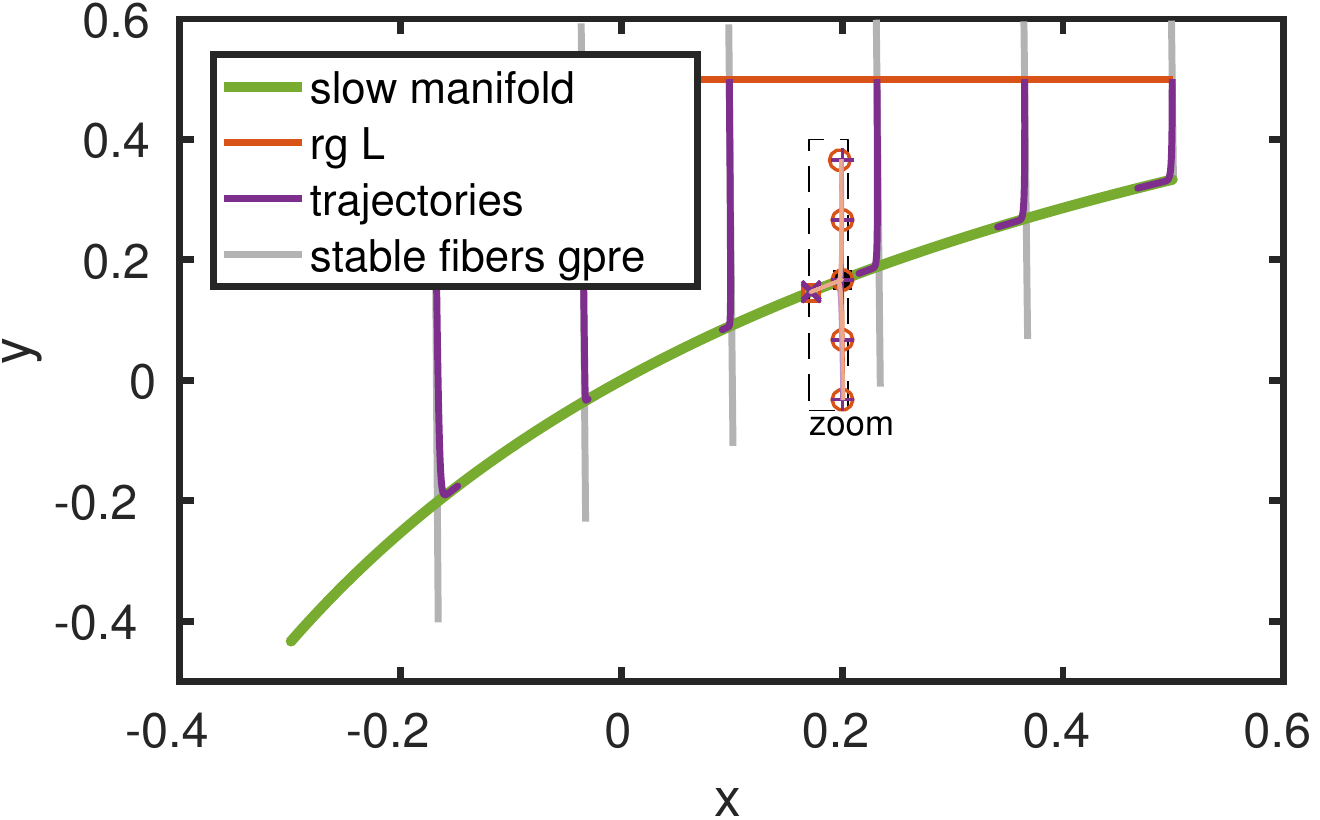}} 
 \subfigure[Geometry in $(v,w)$ coordinates]{\label{fig:mm_dyn_ph_srot}\includegraphics[width=\mmwidth]{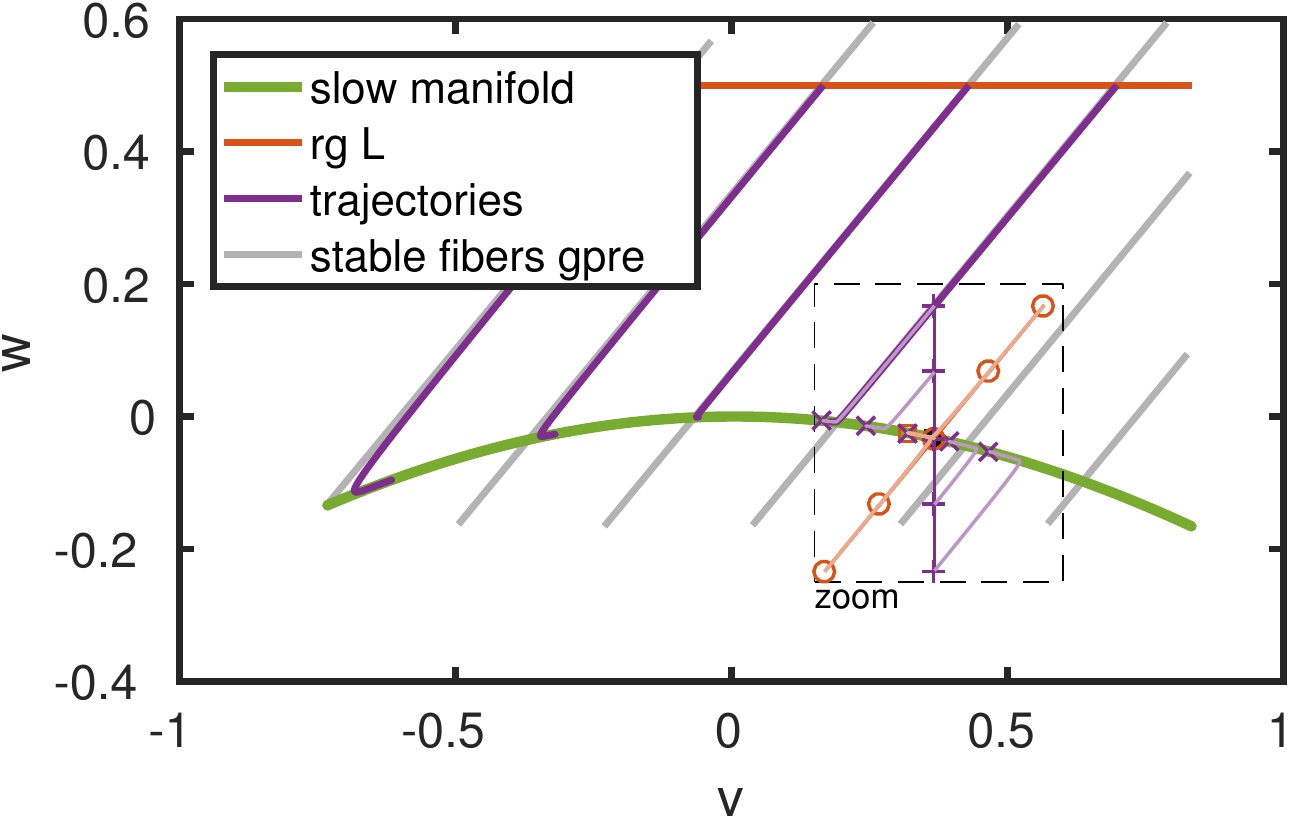}} 
 \subfigure[Zoom to fiber for  \cref{eq:mmdyn}]{\label{fig:mm_dyn_ph_zoom}\includegraphics[width=\mmwidth]{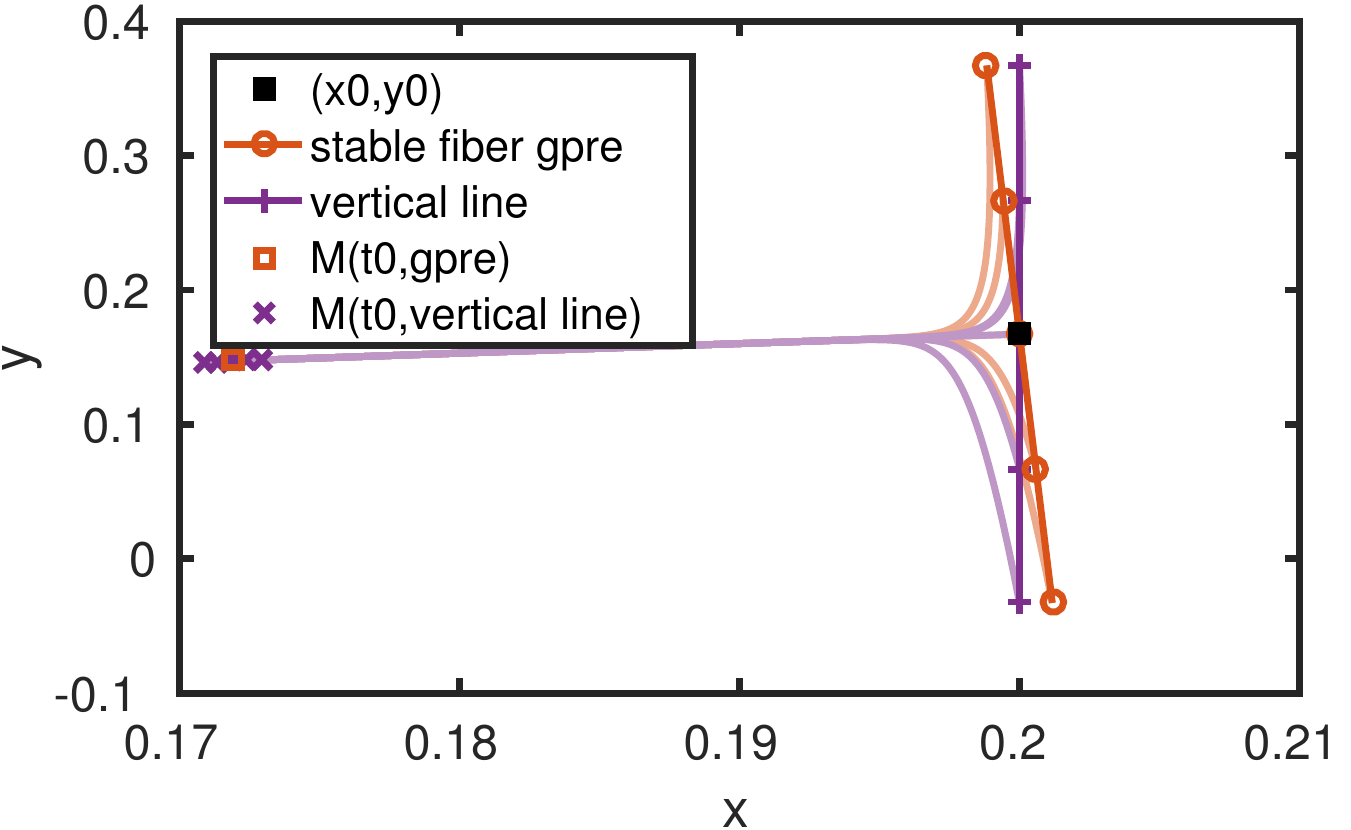}} 
 \subfigure[Zoom to fiber in $(v,w)$ coordinates]{\label{fig:mm_dyn_phr_zoom}\includegraphics[width=\mmwidth]{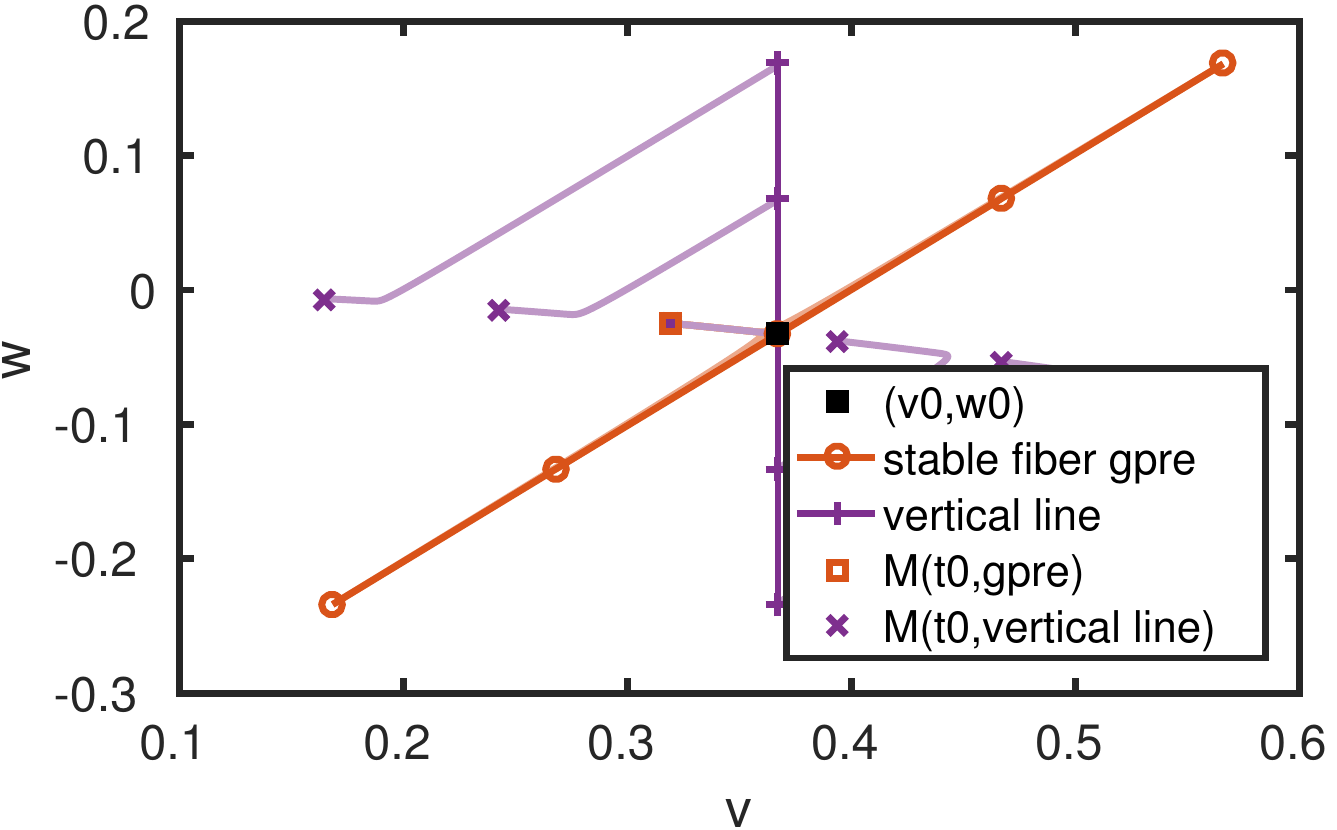}} 
 \subfigure[Error of $\Phi_{\tskip}$ in \cref{eq:mmdyn}]{\label{fig:mm_dyn_e}\includegraphics[width=\mmwidth]{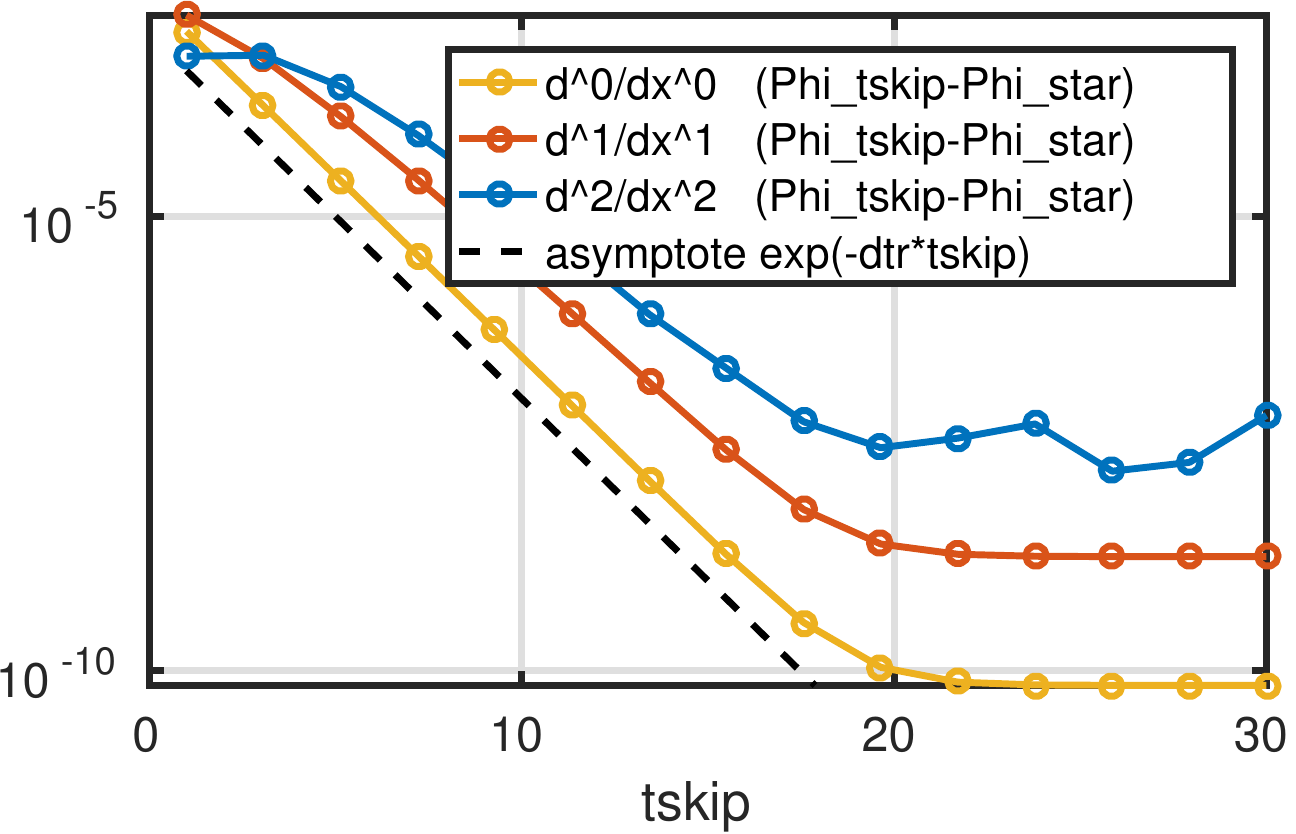}} 
 \subfigure[Error  of $\Phi_{\tskip}$ in $(v,w)$ coordinates]{\label{fig:mm_dyn_e_rot}\includegraphics[width=\mmwidth]{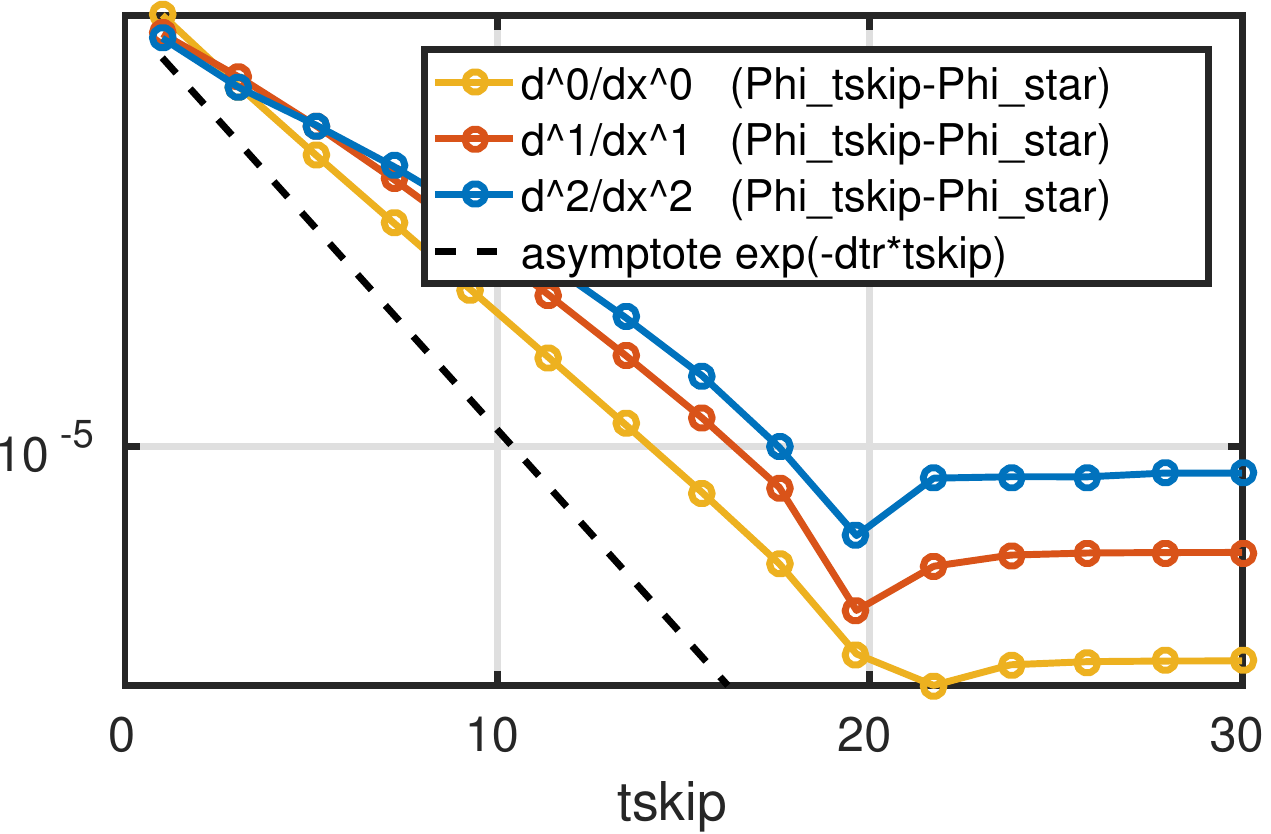}}
 \caption{Michaelis-Menten dynamics: panels \textup{(a,\,c,\,e)} for
   \cref{eq:mmdyn}, panels \textup{(b,\,d,\,f)} for rotated
   coordinates~\cref{eq:rot_dyn}. \textup{(a,\,b):} Location of slow
   manifold ${\cal C}_\epsilon=\{(x,y):y=h_\epsilon(x)\}$ and stable
   fibers in phase space. \textup{(c,\,d):} Trajectories
   $M([0,\mathtt{t0}];U_{\mathrm{ini}})$ \textup{($\mathtt{t0}=20$)}
   starting from two sets $U_{\mathrm{ini}}$ of initial conditions:
   once from a subset of a vertical line, and once from a pre-image of
   a point under the stable fiber projection (called \textup{\textsf{gpre}} in
   the legend). \textup{(e,\,f):} Difference of $\Phi_\tskip(\delta;\cdot)$ to
   third-order asymptotic expansion of $\Phi_*(\delta;\cdot)$ as a function of
   $\tskip$. Parameters: $\lambda=0.5, \kappa = 1, \epsilon = 0.01$ in
   \cref{eq:mmdyn}, $x=-0.1\in\dom\lift$ and $\delta=25$ for
   \textup{(e,f)}.}
  \label{fig:mm_dyn}
\end{figure}

In the singular case $\epsilon = 0$, system \cref{eq:mmdyn} has a
critical manifold $\mathcal{C}_0$ of equilibria, given by the graph
${\cal C}_0=\{(x,y): y=x/(x+\kappa)\}$. For positive $\epsilon$, the
system has a transversally stable invariant manifold, which can be
represented as a graph ${\cal C}_\epsilon=\{(x,y): y=h_\epsilon(x)\}$,
such that $d=1$. In this section we put a subscript $\epsilon$ on
quantities to indicate their dependence on the parameter $\epsilon$
(so, writing, for example, ${\cal C}_\epsilon$ instead of ${\cal
  C}$). The graph $h_\epsilon$ can be expanded in $\epsilon$ for small
$\epsilon>0$:
\begin{equation}
  \label{eq:slow_man_expand}
  y = h_\epsilon(x)=\frac{x}{x+\kappa} + \frac{\kappa \lambda x}{(x+\kappa)^4}
  \epsilon + \frac{\kappa \lambda x (2\kappa \lambda - 3\lambda x -
    \kappa x - \kappa^2)}{(x+\kappa)^7} \epsilon^2 + \order{\epsilon^3}.
\end{equation}
We plan to compare an equation-free approximate flow
$\Phi_\tskip(\delta;\cdot)$, which is constructed below, to the true flow
$\Phi_*(\delta;\cdot)$. For this simple example we may approximate the
true flow $\Phi_*(\delta;\cdot)$ by obtaining an approximation of the
stable fiber projection $g_\epsilon$.  For $\epsilon\to0$ in
\cref{eq:mmdyn}, $g_\epsilon$ has the limit
$g_\epsilon(x,y)\to(x,h_0(x))$. Thus, every point in phase space is
approximately projected along vertical lines, as shown in
\cref{fig:mm_dyn_ph_s}. For positive $\epsilon$ this stable fiber
projection persists and is perturbed by terms of order $\epsilon$. A
general approximation algorithm for stable fibers in slow-fast systems
was provided by Kristiansen \emph{et al}.
\cite{kristiansen2014iterative}. However, we need the stable fibers
only to a degree of accuracy that permits us to compare
$\Phi_*(\delta;\cdot)$ to $\Phi_\tskip(\delta;\cdot)$ for
demonstration purposes. Thus, we expand
$g_\epsilon=(g_{x,\epsilon}(x,y),g_{y,\epsilon}(x,y))$ in
$\epsilon$. Since $g_\epsilon$ projects onto the manifold ${\cal
  C}_\epsilon$, we know that
$g_{y,\epsilon}(x,y)=h_\epsilon(g_{x,\epsilon}(x,y))$.  The
first-order expansion of $g_{x,\epsilon}$ is
\begin{displaymath}
  g_{x,\epsilon}(x,y)=
  x+\frac{(x+\kappa-\lambda)(y-1)x+\kappa y}{x+\kappa}\epsilon+O(\epsilon^2)\mbox{.}
\end{displaymath}
\Cref{fig:mm_dyn} was produced using a third-order expansion of
$h_\epsilon$ and $g_{x,\epsilon}$. The supplementary material provides
Matlab code which reproduces the graphs in \cref{fig:mm_dyn} and
computes the expansion coefficients for $h_\epsilon$ and
$g_{x,\epsilon}$ to third order (see also \cref{sec:supp}).  \Cref{fig:mm_dyn_ph_s,fig:mm_dyn_ph_zoom} show the phase space geometry. The slow manifold
${\cal C}_\epsilon=\{(x,y): x=h_\epsilon(x)\}$ is shown in green, the
stable fibers (pre-images of selected points $(x_j,h_\epsilon(x_j))$
on the slow manifold under the projection $g_\epsilon$) are the almost
straight grey lines, and some sample trajectories of \cref{eq:mmdyn} are
shown in purple. After a rapid transient all trajectories approach the
slow manifold ${\cal C}_\epsilon$, given approximately by
\cref{eq:slow_man_expand}. Furthermore, \cref{fig:mm_dyn_ph_zoom} shows
in more detail how initial conditions on the same stable fiber,
defined as the pre-image of $x_0$ under $g_{x,\epsilon}$,
$G_\mathrm{pre}=\{(x,y): g_{x,\epsilon}(x,y)=x_0\}$, collapse onto the
same slow limiting trajectory (trajectories shown in red in
\cref{fig:mm_dyn_ph_zoom}). In contrast, initial conditions with the same
$y$-component do so only up to an error of order $\epsilon$
(trajectories shown in purple in \cref{fig:mm_dyn_ph_zoom}).

To define the approximate flow $\Phi_\tskip$, we specify the
restriction and lifting operators, $\restrict$ and $\lift$, for the
Michaelis-Menten system as
\begin{align}
  \label{eq:rest_mm}
  \restrict&:\R^2\mapsto\R\mbox{,} & \restrict
  \begin{pmatrix}
    x\\ y
  \end{pmatrix}
  &= x\mbox{,\quad and} &
  \lift&:\R\mapsto\R^2\mbox{,}& \lift (x)
  &= \begin{pmatrix}x\\0.5\end{pmatrix}\mbox{.}
\end{align}

The approximate time-$\delta$ map $\Phi_\tskip(\delta;\cdot)$ on the
slow manifold is determined by the root $z_\tskip$ of
\begin{equation}
  \label{eq:impl_ts}
  F:\R\ni z \mapsto \restrict(M(\tskip;\lift(z))) -\restrict(M(\tskip+\deltat;\lift(x)))
  \in\R
\end{equation}
and setting $\Phi_\tskip(\delta;x) := z_\tskip$ (cf.\ \cref{eq:phip}). We
compare this to the true solution (or, rather, the alternative
approximation by expansion) $\Phi_*(\delta;x)$, determined by the root $z_*$ of
\begin{equation}
  \label{eq:impl_mm_star}
  F_*:\R\ni z\mapsto\restrict(g_\epsilon(\lift(z)))-\restrict(M(\deltat;g_\epsilon(\lift(x))))\in\R\mbox{,}
\end{equation}
setting $\Phi_*(\delta;x):=z_*$.  Note, that $F$ depends on $\tskip$
and $\deltat$, and $F_*$ depends on $\deltat$, which are not included
in the list of arguments to simplify notation.  We solve
\cref{eq:impl_ts} and \cref{eq:impl_mm_star} using a Newton
iteration with tolerance $10^{-12}$, where we approximate $M$ with the
fifth-order component of the \texttt{DOPRI45} Runge-Kutta scheme with
fixed step size $0.1$ for the ODE~\cref{eq:mmdyn}. We approximate the
first two derivatives of $\Phi_*$ and $\Phi_\tskip$ (and the Jacobians
needed inside the Newton iteration) by central finite differences with
step size $\Delta_z=10^{-4}$. The supplementary material contains
  didactic implementations of $\Phi_*$ and $\Phi_\tskip$ for this
  example in the form of matlab code and its published output.  The
  error
\begin{equation}
  E^j(\tskip)=|\partial_2^j\Phi_*(\delta;x_0)-\partial_2^j\Phi_\tskip(\delta;x_0)|\label{eq:error_mm}
  \mbox{\quad ($\partial_2^0\Phi$ refers to $\Phi$)}
\end{equation}
for $x_0=-0.1$ and $\deltat=25$ is shown in \cref{fig:mm_dyn_e} for a
range of healing times $\tskip \in [0;30]$. Since $\epsilon=10^{-2}$
and $\dtan\sim\epsilon$, the quantity $\exp(\dtan\delta)$ is of order
$1$, as required by \cref{thm:conv}.

The error plot \cref{fig:mm_dyn_e} shows that $\Phi_\tskip$,
$\partial_2 \Phi_\tskip$ and $\partial_2^2\Phi_\tskip$ approach a
limit at an exponential rate in $\tskip$ up to an accuracy determined
by the accuracy of the asymptotic expansion of $\Phi_*$ ($\sim
\epsilon^4$), round-off errors in the finite difference approximations
of the derivatives ($\sim 10^{-4}$ for $\partial_2^2\Phi_\tskip$), and
the tolerance of the Newton iteration. Furthermore, the convergence
rate is indeed lower for higher orders of the derivative of
$\Phi_\tskip$ as the estimate \cref{eq:thm:accuracy} in
\cref{thm:conv} suggests. The slope of $\exp(-\dtr\tskip)$ is
included as a lower bound for the error for comparison.

We also observe that the error of the flow and its derivatives is
acceptably small ($\approx10^{-3}$) even for the minimal value $\tskip =
0$. Since $\restrict \circ \lift$ equals the identity (see
\cref{eq:rest_mm}), the implicit equation-free method turns into an
explicit formulation if $\tskip = 0$. 
However, the geometry of system \cref{eq:mmdyn} is not generic. The
system \cref{eq:mmdyn} is given in an explicit slow-fast form with one
fast and one slow variable. This leads with our choice of lifting and
restriction to the degenerate situation that the stable fiber
projection $g_\epsilon$ is aligned with lifting and restriction to
first order: $\restrict\circ g_0\circ\lift=\id$ such that
$\restrict\circ g_\epsilon\circ\lift$ is a small (order $\epsilon$)
perturbation of the identity. In this case the explicit equation-free
method without healing time ($y=\restrict (M(0;\lift
(y)))=\restrict(M(\deltat;\lift (x)))$, such that $\tskip=0$) is accurate up to
order $\epsilon$. To create a situation with a generic arrangement of
the stable fiber projection $g_\epsilon$, we study a rotated system of
the Michaelis-Menten dynamics (which was also used by
\cite{Gear2005}).

We apply the rotation matrix $R$
to the system in order to obtain the dynamics in the new coordinates
$(v,w)^T\in \R^2$ by
\begin{align}
  \label{eq:rot_dyn}
  \begin{pmatrix} v\\w \end{pmatrix} &= R \begin{pmatrix}
    x\\y \end{pmatrix}\mbox{,}& 
    R &= \begin{pmatrix} \phantom{-}1 & 1\\ -1 & 1\end{pmatrix}
  \mbox{, such that}&
  M_\mathrm{rot}\left(t;(v,w)^T\right)&= RM\left(t;R^{-1}(x,y)^T\right)
\mbox{}
\end{align}
is the microscopic simulator in the new $(v,w)^T$ coordinates. In this
rotated system the time scale separation is no longer visible between
$v$ and $w$, since the slow and fast variables are
mixed. \Cref{fig:mm_dyn_ph_srot} shows the phase space geometry with
slow manifold (green), stable fibers (grey) and sample trajectories.
The initial transients are no longer following a straight line
parallel to a coordinate axis such that both, $v$ and $w$, change
rapidly during transients. This situation is expected in a generic
situation when one applies equation-free methods without precise knowledge
about the slow and fast variables. We use the same restriction and
lifting operators as defined in \cref{eq:rest_mm} (but in the new
coordinates $(v,w)$: $\lift(x)=(x,0.5)^T$, $\restrict(v,w)=v$). All
parameter values are as in the unrotated system \cref{eq:mmdyn},
otherwise. The error $E^j$, defined in \cref{eq:error_mm}, is now
much larger: it is of order $1$ for $\tskip = 0$; see
\cref{fig:mm_dyn_e_rot}.  In the implicit framework the error decreases
with increasing healing time $\tskip$ down to $10^{-8}$ for $\tskip =
20$.  Note again that the slope of the curves is smaller for
higher-order derivatives as predicted by the estimate
\cref{eq:thm:accuracy} in \cref{thm:conv}. The error for the flow
is bounded from below again by the accuracy of the asymptotic
expansion for $\Phi_*$, the accuracy of the Newton iteration and
round-off error caused by the finite difference approximations of
$\partial_2^j\Phi_*$ and $\partial_2^j\Phi_\tskip$.

\section{Application: stochastic dynamics}
\label{sec:sde-dynam}
A common area where equation-free methods are applied are
multi-particle systems where slow dynamics emerges for macroscopic
(typically averaged) quantities, e.~g.\ 
\cite{Makeev2002,Barkley2006}. More precisely, the macroscopic
quantities are assumed to satisfy a low-dimensional stochastic
differential equation (SDE). For example, the SDE could be assumed to
be of the form $\d x=f(x) \d t+\sigma \d W_t$, where the noise term
$\sigma\d W_t$ approximates the microscopic fluctuation as white noise
and the deterministic part $f(x)$ is the systematic average drift of
the macroscopic quantities. Givon \emph{et al}. \cite{givon2004} review
rigorous results concerning dimension reduction of SDEs.

\begin{figure}[t]
  \centering
  \includegraphics[width=0.5\textwidth]{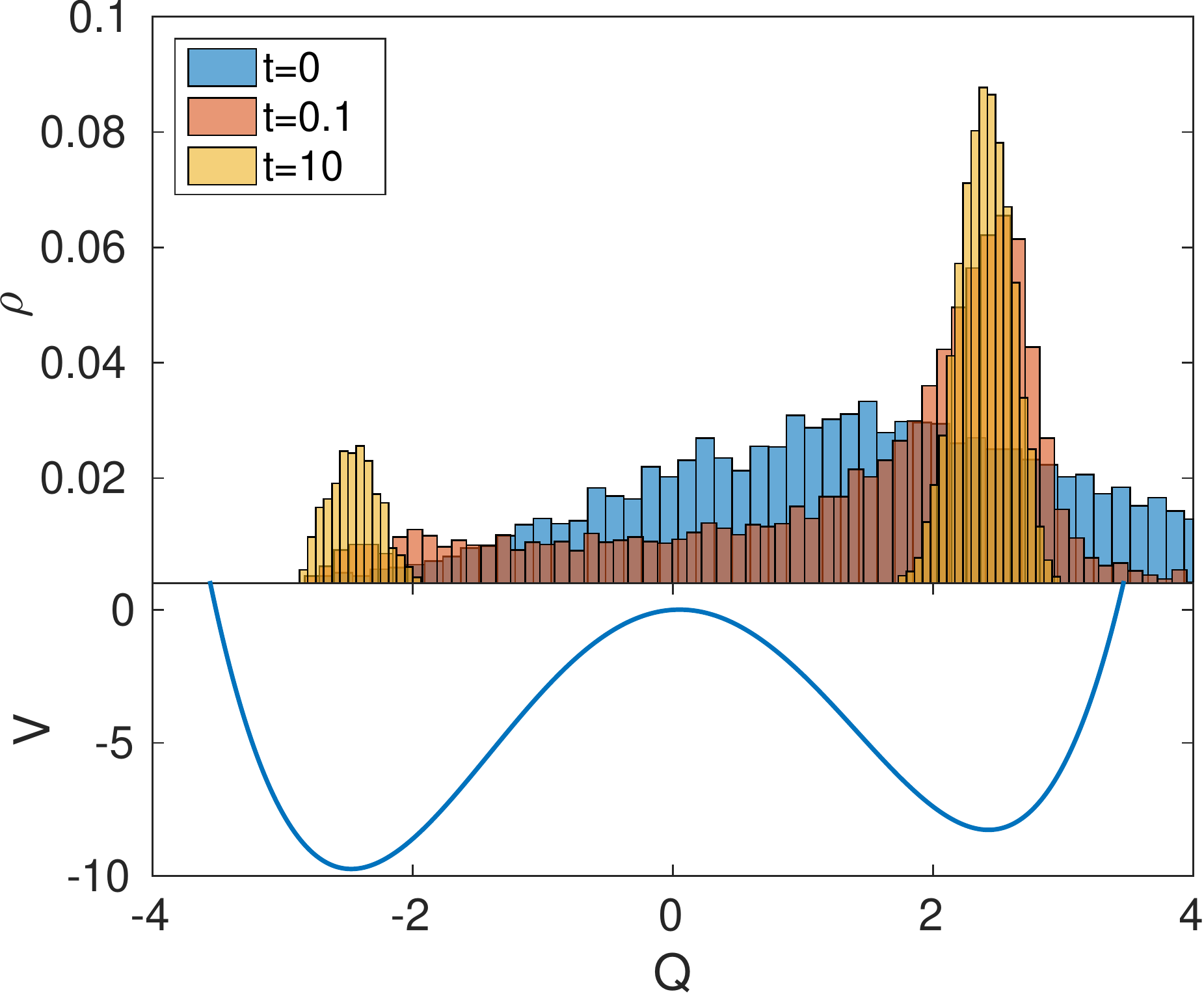}
  \caption{\textup{(}Top\textup{)} Dynamics of distributions for SDE
    model \cref{eq:sde_rescaled} for $\mu=6,\nu=0.3,\sigma = 1$
    \textup{(}sampled from $N=10000$ realizations\textup{)}. On the
    microscopic level of distributions, the Gaussian distributed
    initial condition \textup{($t=0$)} with mean $1.5$ and variance
    $3.5$ converges to a bimodal distribution by $t=10$. Afterwards,
    the transition to the stationary distribution happens
    \textup{(}see mode $1$ in \cref{fig:eigenfunctions}\textup{)} on a
    slow time scale corresponding to the second eigenvalue
    $\lambda_2\sim10^{-8}$ of the operator $L$ given in
    \cref{eq:fp_sde_stat}. \textup{(}Bottom\textup{)} Shape of
    potential well $V(Q)$, as given in \cref{eq:potential_hb}.
  }
  \label{fig:sde:dyn}
\end{figure}

Typically, a stochastic simulation is performed not just once, but for
an ensemble of initial conditions and realizations (as part of a Monte
Carlo simulation). At the level of an SDE, an ensemble of initial
conditions corresponds to (a sampling of) an initial distribution
density $\rho(x)$.  In this section we restrict ourselves to the study
of a scalar SDE of the form
\begin{equation}
  \label{eq:sde_rescaled}
  \d Q = -V'(Q) \d t + \sigma \d W_t,
\end{equation}
where $W_t$ is a Wiener process, an example for which explicit
equation-free methods have been thoroughly analyzed by Barkley
\emph{et al}. \cite{Barkley2006}. As in \cite{Barkley2006}, we set the
noise strength $\sigma$ equal to $1$ in \cref{eq:sde_rescaled}
without loss of generality.  The potential
\begin{equation}
  \label{eq:potential_hb}
  V(Q) = \frac{Q^4}{4} - \frac{\mu Q^2}{2} + \nu Q\mbox{,}
\end{equation}
forms for $\mu>0$ a double well with two local minima $Q_\pm$ and a
local maximum $Q_s$ (see lower panel of \cref{fig:sde:dyn} for a graph
of $V$). The parameters $\mu$ and $\nu$ determine the depth and the
asymmetry of the double-well potential, respectively. We use $\mu=6$,
$\nu=0.3$ such that $Q_-<Q_s<Q_+$ and the well around $Q_-$ is deeper
than the well around $Q_+$. The microscopic simulation is a Monte
Carlo simulation of \cref{eq:sde_rescaled} starting from initial
(ensemble) density $\rho_0(Q)$ of initial conditions. Thus, the phase
space is the space of possible initial distributions in $Q$, which has
dimension $D$ equal to infinity. Strictly, the infinite-dimensional
case is outside of the scope of \cref{thm:conv}. However, the
observations to follow agree with the convergence predicted by the
theorem for reasons that will be discussed after defining lifting,
evolution and restriction. We will make the connection to
multi-particle systems or high-dimensional SDEs in
\cref{sec:discussion}.

\subsection{Lifting, evolution and restriction for distributions}
\begin{figure}[t]
  \centering
  \includegraphics[width = 0.8\textwidth]{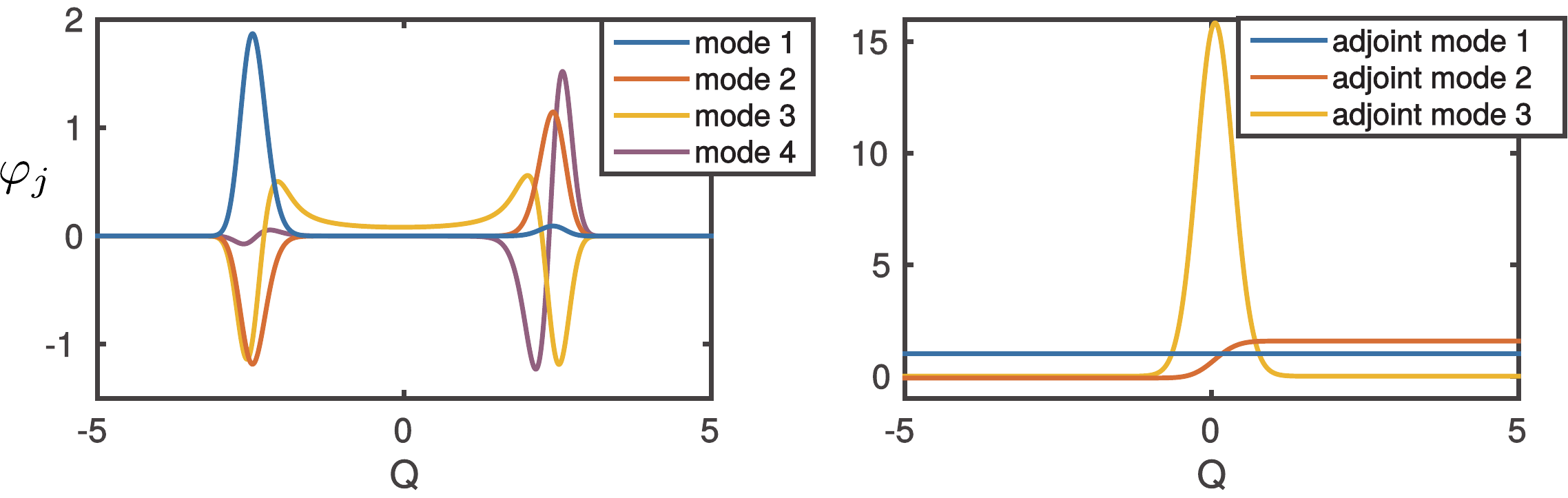}
  \caption{First four eigenfunctions (left panel) and first three
    $\Lint^2$-adjoint eigenfunctions (right panel) of the differential operator
    $L$, defined in \cref{eq:fp_sde_stat} and computed with
    \textup{\texttt{chebfun}} \textup{\cite{chebfun,Driscoll2014}}. The stationary solution is
    shown in blue for $\lambda_1 = 0$ as mode $1$ in the left
    panel. The asymmetric eigenfunction corresponding to $\lambda_2$
    \textup{(}mode $2$\textup{)} is responsible for the transportation
    of mass from one well to the other. Eigenvalues $\lambda_k$ are
    given in \cref{eq:spec:L}\textup{:}
    $(0, -8.63\cdot 10^{-8}, -5.71, -10.3)$. Parameters\textup{:}
    $\mu=6$, $\nu=0.3$, $\sigma=1$. The approximation of the density by Chebyshev
    polynomials is chosen automatically by \textup{\texttt{chebfun}}:
    the degree is $394$ for the left panel and $941$ in the right panel.
  }
  \label{fig:eigenfunctions}
\end{figure}
The evolution of the probability density function (pdf) $\rho(Q,t)$
for the realization of
\cref{eq:sde_rescaled} 
is determined by the Fokker-Planck equation with $\sigma=1$,
\begin{equation}
  \label{eq:fp_sde}
  \partial_t \rho = - \partial_Q (V'(Q)\rho) + \frac{1}{2} \partial_{QQ}\rho\mbox{.}
\end{equation}
The right-hand side of \cref{eq:fp_sde} is linear, of the form
\begin{equation}
  \label{eq:fp_sde_stat}
  L\rho = - \partial_Q (V'\rho) + \frac{1}{2} \partial^2_Q\rho\mbox{,}
\end{equation}
where the operator $L:\Hint^2_1(\R;\R)\mapsto\Lint^2_1(\R;\R)$ is
self-adjoint with respect to the scalar product
\begin{equation}
  \langle\rho_1,\rho_2\rangle_1=\int_\R\frac{\rho_1(Q)\rho_2(Q)}{\phi_1(Q)}\d
  Q\mbox{\quad where \quad}\phi_1(Q)=\frac{\exp(-2V(Q))}{\int_\R\exp(-2V(q))\d q}\mbox{.}\label{eq:scprod}
\end{equation}
The space $\Lint^2_1(\R;\R)$ is in our case the space of all
measurable functions $u:\R\to\R$ with $\int_\R u^2(x)/\phi_1(x) \d
x<\infty$ (a subset of $\Lint^2(\R;\R)$, which has the scalar product
$\langle \rho_1,\rho_2\rangle=\int_\R\rho_1(Q)\rho_2(Q)\d Q$).  The
space $\Hint^\ell_1(\R;\R)$ is the space of all $u\in\Lint^2_1(\R;\R)$
with $u^{(j)}\in\Lint^2_1(\R;\R)$ for all $j\leq \ell$.  The spectrum
of $L$ is real and consists of point spectrum only. It has the form
$0=\lambda_1>\lambda_2>\ldots$ with eigenvectors
$\phi_j(Q)\in\Hint^2_1(\R;\R)$ that can be orthonormalized with
respect to $\langle\cdot,\cdot\rangle_1$. The function $\phi_1$ is the
eigenvector for the trivial eigenvalue $\lambda_1=0$ (which is present
due to the preservation of total probability $\int_\R\rho(Q,t)\d Q$
along trajectories).  The spectrum and the corresponding
eigenfunctions $\phi_j$ are shown in \cref{fig:eigenfunctions} (left
panel), together with the $\Lint^2$-adjoint eigenfunctions
$\phi_j/\phi_1$ (right panel). A solution of the Fokker-Planck
equation~\cref{eq:fp_sde} can be expanded in the eigenfunctions of
$L$ with time-dependent coefficients $a_j(t)$:
\begin{equation}
  \label{eq:dist_expan}
  \rho(Q,t) = \sum_{j=1}^\infty a_j(t) \phi_j(Q).
\end{equation}
The coefficients satisfy $\dot a_j(t)=\lambda_j a_j(t)$ for all $j$,
and the series $\sum_{j=1}^\infty a_j^2$ converges for all $t>0$.  The
orthonormality of the basis $\{\phi_j:j\geq1\}$ with respect to
$\langle\cdot,\cdot\rangle_1$, defined in \cref{eq:scprod}, implies that
\begin{align}
  \label{eq:ef_props}
  \int_{-\infty}^\infty \phi_1(Q) \d Q &= 1\mbox{,}& 
  \int_{-\infty}^\infty
  \phi_j(Q) \d Q &= 0 \mbox{\quad for $j\geq2$.}
\end{align}
Since $\lambda_1=0$, $a_1(t)$ equals $a_1(0)$ for all times
$t\geq0$. One usually chooses $a_1(0)=1$ such that $\rho(Q,t)$
converges to the \emph{stationary density} $\phi_1(Q)$ for
$t\to\infty$.  While \cref{thm:conv} was only formulated
  for flows in $\R^D$, the linearity of $L$ implies that statements
  identical to \cref{thm:conv} can be made for the PDE
  \cref{eq:fp_sde}. Instead of Fenichel's Theorem on invariant
  manifolds in ODEs \cite{Fenichel1979} (persistence and regularity of
  invariant manifolds and fiber projections) we rely on the spectral
  mapping properties for the linear operator $L$.  For any chosen
  dimension $d$ of the slow variables, the slow manifold ${\cal C}$ is
  the subspace spanned by $\phi_1,\ldots,\phi_d$. Instead of the
  stable fiber projection in the finite-dimensional case, we have the
  linear spectral (for $L$) projection
  $g:\Lint^2_1(\R;\R)\mapsto{\cal C}$ ($M$ is also linear, such that we
  write $M(t)\rho$ and $g\rho$) which is explicitly known in terms of
  the eigenvectors of $L$:
\begin{equation}
  \label{eq:sde_g}
  g:\Lint^2_1(\R;\R)\ni \rho\mapsto g\rho=\sum_{j=1}^d
  \langle \phi_j,\rho\rangle_1
  \phi_j\in {\cal L}(\phi_1,\ldots,\phi_d)\subset\Lint^2_1(\R;\R)\mbox{.}
\end{equation}
With this definition of ${\cal C}$ and $g$ the decay and growth
properties of the evolution map $M$ replacing
\cref{ass:timescales} are
\begin{align}
        \label{eq:sdebound}
        \|M(t)\vert_{\cal C}\| &\leq C\exp(\lambda_d t)&&\mbox{for all
          $t\leq0$,} & \|M(t)\vert_{\cal C}\| &\leq C&&\mbox{for all
          $t\geq0$}\\
      \label{eq:sdeconv}
      \|M(t)-M(t)\circ g\|&\leq 
      C\exp(\lambda_{d+1}t)&&\mbox{for all $t\geq0$}
\end{align}
and some constant $C$, such that $\dtan=-\lambda_d$,
$\dtr=-\lambda_{d+1}$. The approximation statement of
\cref{thm:conv} then follows immediately from error estimates
for finite-dimensional matrices and will be derived after the
definition of the restriction and lifting operators. The lifting and
restriction operators are chosen to map from a macroscopic description
of $\rho$, for example, by moments, to the full density $\rho$ and
vice versa.  In particular, we will investigate the behaviour of
implicit equation-free methods for $d=3$ using the following
restriction and two different lifting operators:
\begin{align}
  \label{eq:sde_restrict}
  \restrict&:\Lint^2_1(\R;\R)\mapsto\R^d & \restrict \rho &=
  \left(\int_\R Q^{k-1}\rho(Q)\d Q\right)_{k=1,\ldots,d}\mbox{,}\\
  \label{eq:sde_llin}
  \llin&:\R^d\mapsto\Lint^2_1(\R;\R) &
  \llin(x)(Q) &=\sum_{j=1}^{d} x_j\rho_j(Q)\mbox{,} \\ 
  \label{eq:sde_lgauss}
  \lgauss&:\R^3\mapsto\Lint^2_1(\R;\R) &
  \lgauss(x)(Q)&=\frac{x_1}{\sqrt{2\pi x_3}}\exp\left(\frac{-(Q-x_2)^2}{2x_3}\right)\mbox{.}
\end{align}
Thus, $\restrict$ projects a density onto its first $d$ moments
(counting from the zeroth moment, which is preserved by $M$ since
$\lambda_1=0$). In a Monte-Carlo simulation the zeroth moment would
correspond to the (possibly scaled) number of realizations. The
functions $\rho_j$ in the definition \cref{eq:sde_llin} of $\llin$
are arbitrary in $\Lint^2_1(\R;\R)$ with $\int_\R\rho_j(Q)\d Q=1$, which
ensures that $\int_\R\llin(x)(Q)\d Q=\sum_{j=1}^d x_j$ is conserved
under $M(t)$.  For $\lgauss$, the $x_1$ component is preserved under
$M(t)$ and becomes the first component of $\restrict$ such
that always $[\restrict M(t)\lgauss x]_1=x_1$.

For the combination of $\llin$ and $\restrict$ all components of the
lift-evolve-restrict map $P(t;\cdot)=\restrict\circ
M(t;\cdot)\circ\lift$ and its exact counterpart $P_*(t;\cdot)=\restrict\circ
M(t;\cdot)\circ g\circ\lift$ from
\cref{sec:slowfast} are linear such that we can reduce the
study of convergence for arbitrary coordinates $x$ to convergence
estimates for matrices.

The combination of $\lgauss$ and $\restrict$ was studied in detail in
\cite{Barkley2006} for explicit equation-free methods, where the
authors observed that the nonlinearity of $\lgauss$ introduced a
nonlinearity in the moment map and that the resulting flow depended
qualitatively on the choice of the healing time $\tskip$. We will
demonstrate that for $\lgauss$ the implicitly defined flow
$\Phi_{\gauss,\tskip}$ converges to a nonlinear transformation of the
linear flow $M(t)\vert_{\cal C}$. Since the $x_1$ component does not
change under $P(t;\cdot)$ and $P_*(t;\cdot)$ for $\lift=\lgauss$, it
can be ignored, making the choice of $\lgauss$ and $\restrict$
identical to the situation studied in \cite{Barkley2006}.

We use the \textsc{MATLAB} \cite{matlab} package \texttt{chebfun}
\cite{chebfun,Driscoll2014} to numerically compute the spectrum and
eigenfunctions of $L$, the flow $M$, the projection $g$, restriction
and lifting for the example potential $V$ given in
\cref{eq:potential_hb}. The package \texttt{chebfun} uses Chebyshev
polynomials of adaptive degree to approximate arbitrary functions on
finite intervals to optimal precision. For a typical result, as shown
in \cref{fig:eigenfunctions} the degree is larger than $100$
($394$ for the left panel, $941$ for the right panel). The numerically
computed spectrum of $L$ is
\begin{align}\label{eq:spec:L}
  \spec(L) &= (\lambda_1,\lambda_2,\lambda_3,\lambda_4,\ldots) \approx (
  -2.37\cdot
  10^{-9}, -8.63\cdot 10^{-8}, -5.71, -10.3, \ldots)\\
  \intertext{for $V$ with the
    parameters} \mu&=6\mbox{,\qquad}\nu=0.3\mbox{.}\label{eq:Vpar}
\end{align}
Note that $\lambda_1 = 0$ is the correct value for the first
eigenvalue on an infinite domain. In numerical computations we choose
a bounded domain $[-10,10]$ with Dirichlet boundary conditions,
leading to a small probability of escape from the domain.  The
spectrum and the corresponding eigenfunctions $\phi_j$ are shown in
\cref{fig:eigenfunctions}. The eigenvector $\phi_1$ corresponds to the
stationary solution of the Fokker-Planck equation and $\phi_2$ is the
mode representing escape from one well to another.

\subsection{Convergence for the linear lifting operator \texorpdfstring{$\llin$}{Llin} with
  $d=3$}
\label{sec:llin} We express the maps $P_*(t;\cdot)$ and $P(t;\cdot)$
in terms of $M$, the eigenvectors $\phi_j$ and the scalar product
$\langle\cdot,\cdot\rangle_1$, initially for a general dimension
$d$. The exact macroscopic flow $\Phi_*$ is defined using the map
$P_*(t;\cdot)$ in \cref{eq:phistarp}, and the approximate macroscopic
flow $\Phi_\tskip$ is defined using the map $P(t;\cdot)$ in
\cref{eq:phip}. The definitions \cref{eq:sde_restrict} for
$\restrict$ and \cref{eq:sde_llin} for $\llin$ imply
\begin{align}
  \left[P_{\lin,*}(t)x\right]_k&=[\restrict M(t) g \llin
  x]_k=\sum_{\ell,\,j=1}^{d}\exp(\lambda_\ell t)\int_\R
  Q^{k-1}\phi_\ell(Q)\d Q \langle\phi_\ell,\rho_j\rangle_1x_j\label{eq:pslin}\\
  \left[P_\lin(t)x\right]_k&=[\restrict M(t) \llin
  x]_k=\sum_{j=0}^{d}\left[\int_\R Q^{k-1}[M(t)\rho_j](Q)\d Q\right]
  x_j\label{eq:plin}
\end{align}
where $k=1,\ldots,d$. Using the $d\times d$
matrices ($k,\ell,j=1,\ldots,d$)
\begin{align}\label{eq:sde:tlin:M:R}
  (T_{\lin})_{\ell,j}&=\langle\phi_\ell,\rho_j\rangle_1\mbox{,} &
  M_d(t)&=\diag\left[
    \exp(\lambda_\ell t)_{\ell=1,\ldots,d}\right]\mbox{,} &
  (R_d)_{k,\ell}&=\int_\R Q^{k-1}\phi_\ell(Q)\d Q\mbox{.}
\end{align}
we can express the map $P_{\lin,*}(t)$ and the exact slow flow
$\Phi_{\lin,*}(\delta)$ in the form
\begin{align}
  \label{eq:lin:pslinmat}
  P_{\lin,*}(t)x &=R_dM_d(t)T_\lin x
  \\
  \label{eq:lin:phis}
  \Phi_{\lin,*}(\delta) x &=P_{\lin,*}(\tskip)^{-1}P_{\lin,*}(\tskip+\delta)=T_\lin^{-1}M_d(\delta)T_\lin x \mbox{.}
\end{align}
General \cref{ass:transversality} on transversality of $\restrict$ and
$\lift$ for \cref{thm:conv}, when applied to the SDE
\cref{eq:sde_rescaled} and $\restrict$ and $\llin$, demands the
regularity of the matrices $R_d$ and $T_\lin$. If both matrices are
indeed regular then $P_{\lin,*}$ is invertible:
$P_{\lin,*}(t)^{-1}=T_\lin^{-1}M_d(-t)R_d^{-1}$. Thus, the claim of
\cref{thm:conv} can be simplified to a statement about perturbations
of matrices using the quantities
\begin{align}
  \label{eq:pslin:inv}
  n(t)&:=\|P_{\lin,*}(t)^{-1}\|\leq\|T_\lin^{-1}\|\|R_d^{-1}\|\exp(-\lambda_dt)
  &&\mbox{for $t\geq0$, and}\\
  r(t)&:=\|P_{\lin,*}(t)-P_\lin(t)\|\leq
  C\exp(\lambda_{d+1}t) && \mbox{for all $t\geq0$ and a fixed $C>0$.}\label{eq:plin-pslin}
\end{align}
The estimate for $r(t)$ follows from
\cref{eq:pslin,eq:plin}:
  \begin{align*}
    \left|\left[P_\lin(t)x\right]_k-\left[P_{\lin,*}(t)x\right]_k\right|&=\left|\sum_{j=0}^{d}\left[\int_\R
        Q^{k-1}\left[M(t)\left[\rho_j-\sum_{\ell=1}^d\langle
  \phi_\ell,\rho_j\rangle\phi_\ell\right]\right](Q)\d Q\right]
      x_j\right|\nonumber\\
    &\leq C\exp(\lambda_{d+1}t)|x|\mbox{.}
  \end{align*}
  The integrand contains the spectral projection
  $\rho\mapsto\rho-\sum_{\ell=0}^d\langle
  \phi_\ell,\rho_j\rangle\phi_\ell$ onto the complement of the space
  spanned by $\phi_1,\ldots,\phi_d$. On the complement of ${\cal
    L}(\{\phi_1,\ldots,\phi_d\})$ the evolution operator $M(t)$ decays
  exponentially with rate $\lambda_{d+1}$ in time. Together with the
  boundedness of $\restrict$ and the spectral projection, this decay
  of $M(t)$ implies estimate \cref{eq:plin-pslin}. Since
$P_\lin(t;\cdot)$ is linear, the approximate flow map
$\Phi_{\lin,\tskip}(\delta)$ is given by
\begin{equation}\label{eq:philin}
  \Phi_{\lin,\tskip}(\delta)=P_\lin(\tskip)^{-1}P_\lin(\tskip+\delta)
\end{equation}
assuming the inverse of $P_\lin(\tskip)$ exists (all involved matrices
have dimension $d\times d$).  The linear expressions for the exact flow $\Phi_{\lin,*}$,
\cref{eq:lin:phis}, and $\Phi_{\lin,\tskip}$, \cref{eq:philin}, imply
\begin{equation}
  \label{eq:lin:err:est}
  \|\Phi_{\lin,*}(\delta)-\Phi_{\lin,\tskip}(\delta)\|\leq
  \frac{n(\tskip)}{1-n(\tskip)r(\tskip)}
  \left[r\left(\tskip\right)\left\|\Phi_{\lin,*}(\delta)\right\|+r(\tskip+\delta)\right]\mbox{.}
\end{equation}
The exponential estimates $r(t)\leq C \exp(\lambda_{d+1} t)$ in
\cref{eq:plin-pslin} and $n(t)\leq \exp(-\lambda_d t)$ in
\cref{eq:pslin:inv} immediately imply the statement of \cref{thm:conv}
(given that $\Phi_{\lin,*}$ is globally bounded).
\begin{figure}[ht]
  \centering
  \subfigure[Convergence of $\Phi_{\lin,\tskip}$ to $\Phi_{\lin,*}$]{\label{fig:plin:convergence}\includegraphics[height=41ex]{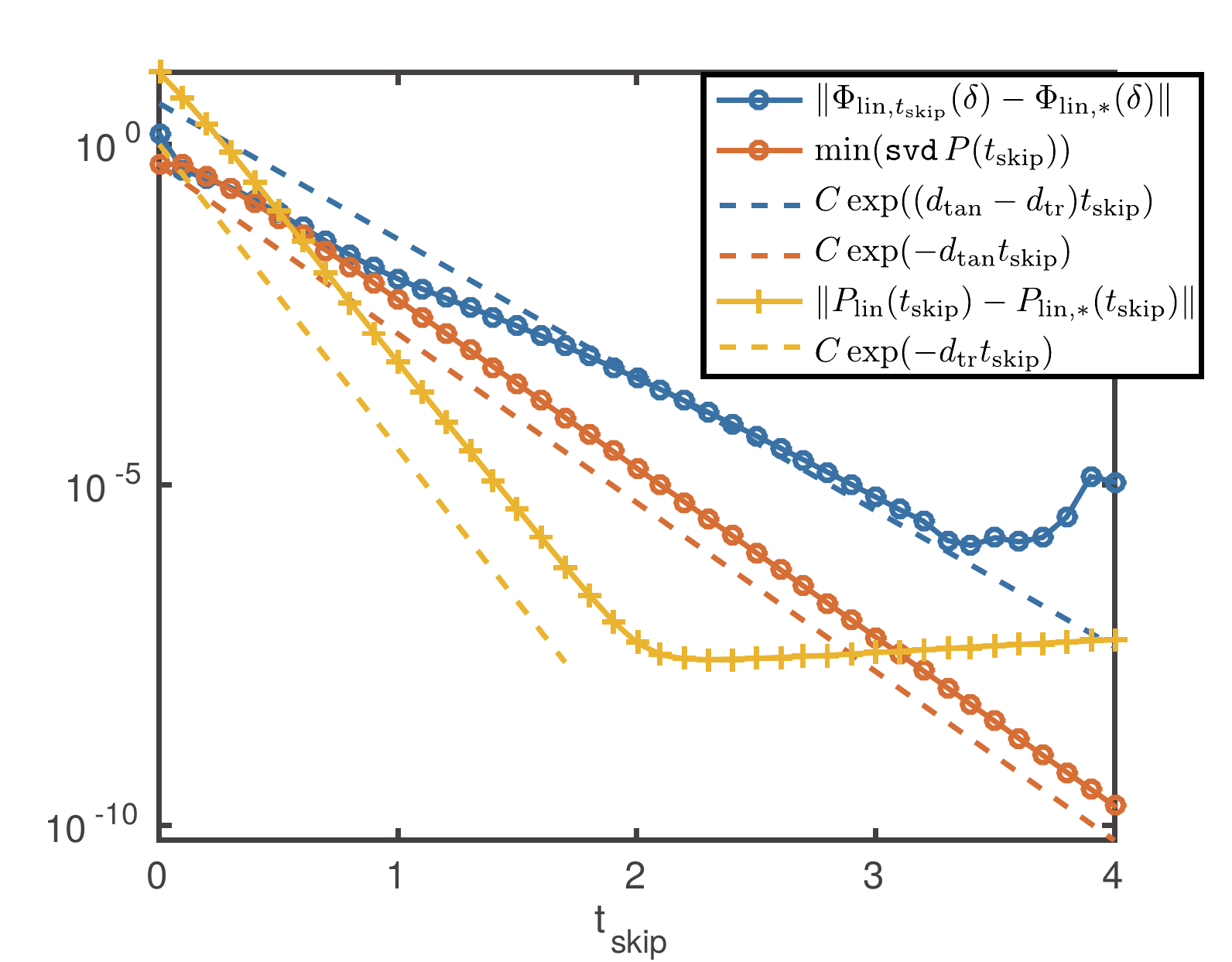}}\hfill
    \subfigure[Phase portrait of $\Phi_{\lin,*}$]{\label{fig:plin:phaseportrait}\includegraphics[height=39ex]{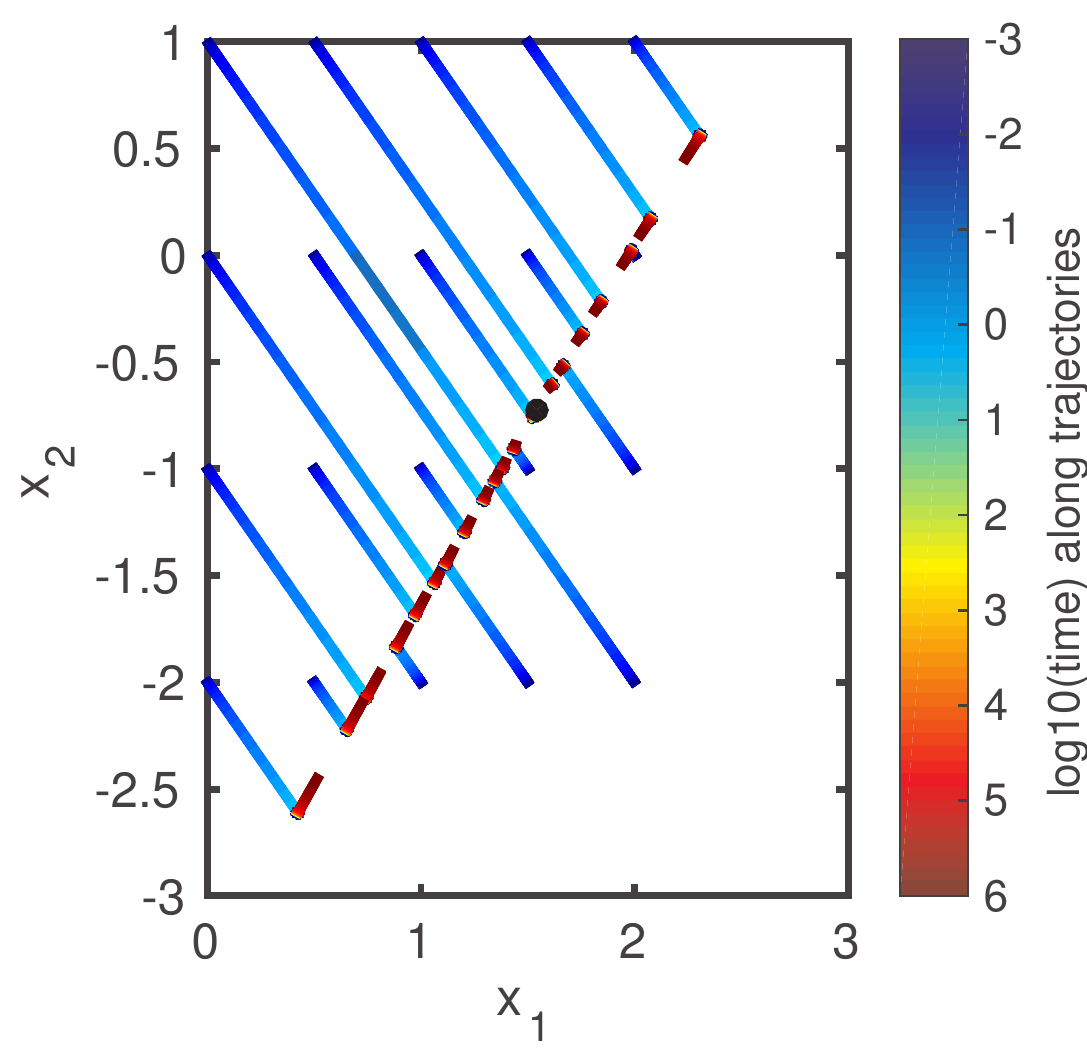}}
    \caption{Analysis for $\llin$ and $\Phi_{\lin,\tskip}$ with
      $d=3$\textup{.}  \textup{(}a\textup{)} Matrix norm of
      $\Phi_{\lin,\tskip}(\delta)-\Phi_{\lin,*}(\delta)$ and smallest
      singular value of $P(\tskip)$, compared to theoretical
      exponential estimates. \textup{(}b\textup{):} the plot shows
      trajectories of $\Phi_{\lin,*}$ in the $(x_1,x_2)$-plane
      starting on a grid of initial points. The coloring indicates the
      time along the trajectory \textup{(}in logarithmic
      scale\textup{)}. The black point is the fixed point (see
      text). Parameters\textup{:} $\delta=0.1$, $\rho_j$
      \textup{(}$j=1\ldots3$\textup{)} are Gaussians with means
      $-1.5$, $-0.5$ and $1$ and variance $1$, shape parameters of
      potential $V$ are $\mu=6$, $\nu=0.3$, see also
      \cref{eq:spec:L,eq:Vpar}.}
  \label{fig:plin}
\end{figure}

The semilogarithmic plot in \cref{fig:plin:convergence} shows
the difference  between $\Phi_{\lin,\tskip}(\delta)$ and
$\Phi_{\lin,*}(\delta)$ (blue line with circles) for $d=3$, for a
linear basis of three Gaussians $\rho_j$ (with variance $1$ and means
$-1.5$, $-0.5$ and $1$), $\delta=0.1$ and the double-well potential
well $V$ with parameters $\mu=6$, $\nu=0.3$. The decay rate inside the
slow manifold ${\cal C}=\linspan(\phi_1,\phi_2,\phi_3)$ is
$\dtan=-\lambda_3\approx5.71$ and the attraction rate toward ${\cal
  C}$ is $\dtr=-\lambda_4\approx10.3$.
\Cref{fig:plin:convergence} also shows the two components of the
error $\Phi_{\lin,\tskip}(\delta)-\Phi_{\lin,*}(\delta)$ and their
theoretical estimates:
\begin{enumerate}
\item (In yellow) The difference $r(\tskip)$ between
  $P_\lin(\tskip)=\restrict\circ M(\tskip)\circ\llin$ and
  $P_{\lin,*}(\tskip)=\restrict\circ M(\tskip)\circ
  g\circ\llin=R_dM_d(\tskip)T_\lin$, which decays according to the
  attraction toward ${\cal C}$ until it reaches the limits of
  numerical accuracy of \texttt{chebfun} ($\sim10^{-8}$):
  $\|P_\lin(\tskip)-P_{\lin,*}(\tskip)\|\sim\exp(-\dtr\tskip)$.
\item (In red) The norm $n(\tskip)$ of the inverse of $P_\lin(\tskip)$, which grows
  like $\exp(\dtan\tskip)$. \Cref{fig:plin:convergence} shows
  the inverse (the minimal singular value).
\end{enumerate}
The overall error \cref{eq:lin:err:est} is approximately the product
of these two components, which is proportional to
$\exp((\dtan-\dtr)\tskip)$ (shown as a blue dashed line in
\cref{fig:plin:convergence}). In particular, the combination of
$\|P_{\lin,*}(\tskip)^{-1}\|\sim\exp(\dtan\tskip)$,
$\|P_\lin(\tskip)-P_{\lin,*}(\tskip)\|\sim\exp(-\dtr\tskip)$ and
$\dtan<\dtr$ implies that $P_\lin(\tskip)$ is invertible for
sufficiently large $\tskip$ and that
$\|P_\lin(\tskip)^{-1}\|\sim\exp(\dtan\tskip)$.

\Cref{fig:plin:phaseportrait} shows a phase portrait of the
exact flow in the coordinates in $\dom\llin$. Since $\Phi_{\lin,*}$
and $\Phi_{\lin,\tskip}$ both preserve the quantity $\sum_{j=1}^dx_j$
(which corresponds to $\int_\R\llin(x)(Q)\d Q$), we set
$x_3=1-x_1-x_2$ in the initial values for the sample trajectories,
keeping $\sum_{j=1}^dx_j=1$ along trajectories without loss of
generality. This leads to an affine flow in the $(x_1,x_2)$-plane with
a non-trivial fixed point (shown in black in
\cref{fig:plin:phaseportrait}). The coloring along the sample
trajectories illustrates the extreme difference in the time scale
along the directions corresponding to $\lambda_2$ ($\approx -10^{-7}$;
escape between wells), mostly evolving on time scales $\gg10^4$
(dark red in \cref{fig:plin:phaseportrait}), and $\lambda_3$
($\approx-5.71$; relaxation into the nearest well), mostly decaying
on time scale of order $1$ and less (blue and light blue in
\cref{fig:plin:phaseportrait}).

\paragraph{ Remark: Densities with sign changes in
  \cref{sec:llin}} The phase portrait
\cref{fig:plin:phaseportrait} of the exact flow $\Phi_{\lin,*}$
includes coordinates $x=(x_1,x_2,x_3)$ where the lifted initial
density $\llin(x)$ has sign changes. This is not unphysical. If
one performs Monte Carlo simulations with ensembles on the example
with the lifting operator $\llin(x)$, one would run a Monte Carlo
simulation on an ensemble for each of the three initial densities
$\rho_j$. Then one would sum the densities at the end of the
simulation with the weights $x_j$ ($j=1,\ldots,3$). These weights can
be negative to get a combined density.


\subsection{Convergence for the nonlinear lifting operator \texorpdfstring{$\lgauss$}{Lgauss}}
\label{sec:lgauss}
\begin{figure}[ht]
  \centering
  \includegraphics[width=\textwidth]{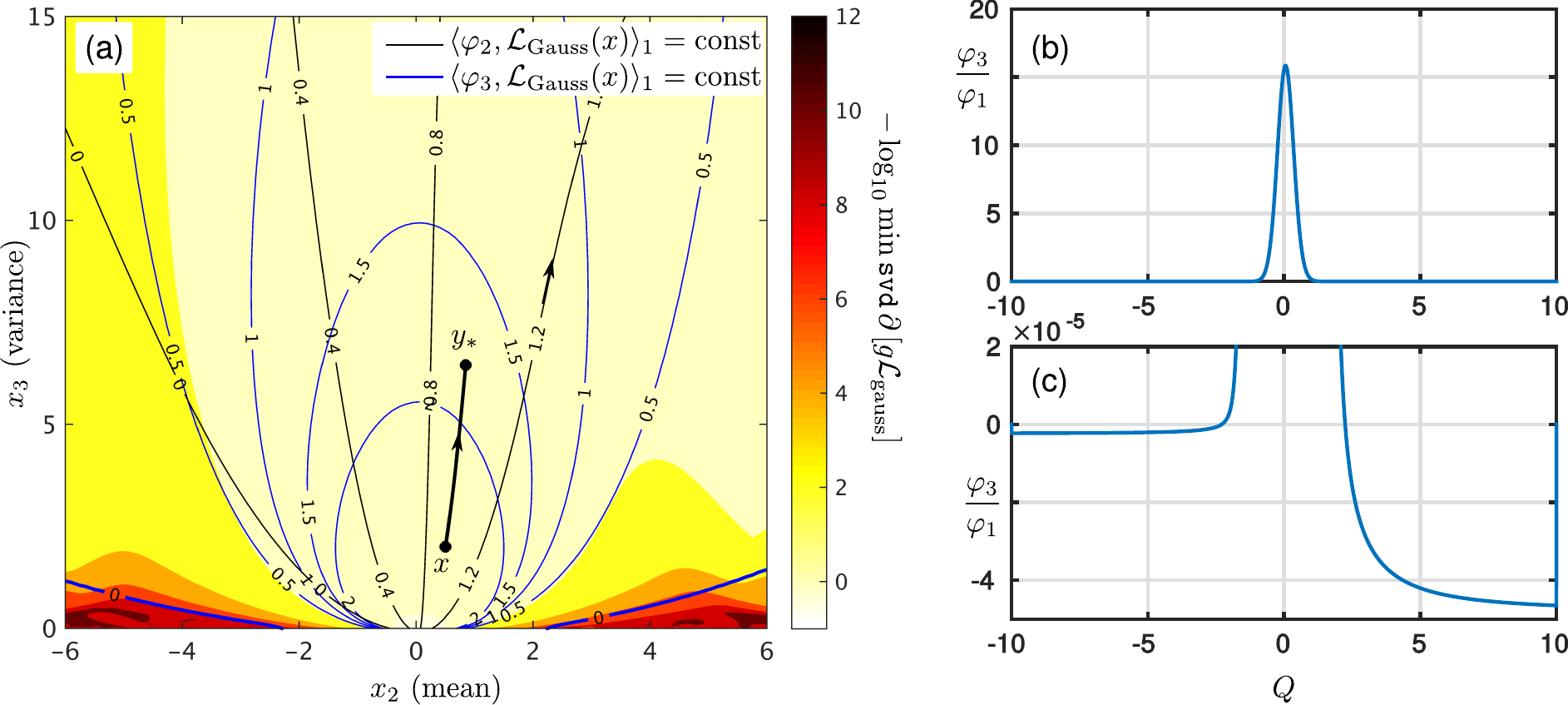}
  \caption{\textup{(}a\textup{)} Contour curves in coordinates in
    $\dom\lgauss$ and condition of $\partial
    T_\gauss=\partial[g\lgauss]$. The piece of trajectory from $x$ to
    $y_*=\Phi_{\gauss,*}(\deltat;x)$ with $\deltat=0.1$ is used for
    convergence analysis in
    \cref{sec:lgauss}. \textup{(}b\textup{)} Profile of
    $\phi_3(Q)/\phi_1(Q)$ in $Q$. Note that away from its peak around
    $0$, the profile is slightly negative \textup{(}see zoom in
    panel \textup{(}c\textup{))}. Parameters\textup{:} $\mu=6$, $\nu=0.3$.}
  \label{fig:lgauss:phaseportrait}
\end{figure}
The exact and approximate lift-evolve-restrict maps for lifting with a
Gaussian distribution of mass $x_1$, mean $x_2$ and variance $x_3$, of the form
$\lgauss(x)=q\mapsto x_1\exp(-(q-x_2)^2/(2x_3))/\sqrt{2\pi x_3}$, are given by
\begin{align}\nonumber
  P_{\gauss,*}(t;x)_k&=\left[\restrict M(t) g \lgauss(x)\right]_k\\
  &=
  \sum_{\ell=1}^3\exp(\lambda_\ell t)\int_\R Q^{k-1}\phi_\ell(Q)\d Q
  \langle\phi_\ell,\lgauss(x)\rangle_1\label{eq:psgauss}\\
  P_\gauss(t;x)_k&=\left[\restrict M(t)\lgauss(x)\right]_k=
  \int_\R Q^{k-1}M(t)\lgauss(x)\d  Q\mbox{,}\label{eq:pgauss}
\end{align}
where $k=1,\ldots,d$ ($d=3$).  The flow $M(t)$ preserves the integral
of the initial distribution such that $P_\gauss(t;x)_1=x_1$ and
$P_{\gauss,*}(t;x)_1=x_1$. Thus, we can fix $x_1=1$ without loss of
generality and focus on the dynamics in the $(x_2,x_3)$-plane in
$\dom\lgauss$.
\paragraph{Phase portrait of the exact flow $\Phi_{\gauss,*}$}
The exact flow
$\Phi_{\gauss,*}$ on ${\cal C}$ in the coordinates of $\dom\lgauss$ is
a nonlinear transformation of the linear map $M_d(t)=\diag\left[
  \exp(\lambda_\ell t)_{\ell=1}^{\ell=d}\right]:\R^d\mapsto\R^d$,
defined in \cref{eq:sde:tlin:M:R} (with $d=3$). We call the nonlinear
transformation
\begin{align}
  \label{eq:tgauss}
  T_\gauss&:\R^3\mapsto\R^3 &
  T_\gauss(x)_k&=\langle\phi_k,\lgauss(x)\rangle_1
  &&\mbox{($k=1,\ldots,3$).}
\end{align}
In particular $T_\gauss(x)_1=x_1$ by construction.  Using $T_\gauss$,
$M_d$ and the matrix $R_d$ (defined in \cref{eq:sde:tlin:M:R}), the map
$P_{\gauss,*}(t;x)$, and the exact flow $\Phi_{\gauss,*}$ are given by
(using the notation $T_\gauss^{-1}$ for the inverse of the nonlinear
map $T_\gauss$)
\begin{equation}\label{eq:phigauss}
  \begin{split}
    P_{\gauss,*}(t;x)&=R_dM_d(t)T_\gauss(x)\mbox{,}\\
    \Phi_{\gauss,*}(\deltat;x)&= T_\gauss^{-1}(M_d(\deltat)
    T_\gauss(x))\mbox{,}
  \end{split}
\end{equation}
where all involved quantities are maps from $\R^3$ to $\R^3$. Since the
map $T_\gauss$ is nonlinear, it is not clear if the inverse exists for
all $x\in\R^3$, or if it is unique where it
exists. \Cref{fig:lgauss:phaseportrait}(a) shows the contours of
$T_\gauss(x)_2$ (in black) and $T_\gauss(x)_3$ (in blue; remember that
$T_\gauss(x)_1=x_1$), and the norm of $[\partial T_\gauss(x)]^{-1}$ as
color shading (in logarithmic scale). Since the difference in time
scale between motion along $\phi_2$ and motion along $\phi_3$ is large
($0>\lambda_2\gg\lambda_3$), the flow $\Phi_{\gauss,*}$ follows the
black curves in the direction of the arrow until it reaches the
zero-level of $T_\gauss(x)_3$ (slightly wider blue curve, only visible
close to the bottom of \cref{fig:lgauss:phaseportrait}(a)).

\paragraph{Near-singularity of $T_\gauss$}
The zero curve $\{x:T_\gauss(x)_3=0\}$ in the $(x_2,x_3)$ plane
(wide blue in \cref{fig:lgauss:phaseportrait}(a)) is given by
$\int_\R\lgauss(x)(Q)\phi_3(Q)/\phi_1(Q)\d Q=0$, where $\lgauss(x)$ is
a Gaussian of mean $x_2$ and variance $x_3$ and $\phi_3(Q)/\phi_1(Q)$
is shown in \cref{fig:lgauss:phaseportrait}(b,c). From the profile of
$\phi_3/\phi_1$ it is clear that the zero-level forms a single curve
connecting the two pieces of the wide blue curve
$\{x:T_\gauss(x)_3=0\}$ visible in
\cref{fig:lgauss:phaseportrait}(a). However, this curve has a large
radius (passing through the region $x_3\gg1$). For example, there
exists a Gaussian $u=\lgauss(x)$ with mean $x_2=0$ and large variance
$x_3$ such that $T_\gauss(x)_3=0$, because $\phi_3/\phi_1$ is negative
everywhere outside its peak, but the negative values have small
modulus (note the scaling of the vertical axis in the zoom of
$\phi_3/\phi_1$ in \cref{fig:lgauss:phaseportrait}(c)). The fixed
point of $x\mapsto \Phi_{\gauss,*}(\deltat;x)$ (assuming $x_1=1$) is
the intersection of the two zero-level curves (not visible in
\cref{fig:lgauss:phaseportrait}(a) as it has large $x_3$). The color
shading in \cref{fig:lgauss:phaseportrait}(a) indicates that the
nonlinear transformation $T_\gauss$ is nearly singular close to the
line $x_3=0$, because the $\Lint^2$-adjoint modes $\phi_2/\phi_1$ and
$\phi_3/\phi_1$ are both nearly constant away from the region around
$Q\in[-2,2]$ (see \cref{fig:eigenfunctions}, right panel) such that,
when inverting $T_\gauss$, the mean $x_2$ is very sensitive for small
changes in the coefficients for the $\Lint^2$-adjoint modes
$\phi_2/\phi_1$ and $\phi_3/\phi_1$ .

\paragraph{Components of error $\Phi_{\gauss,\tskip}-\Phi_{\gauss,*}$}
We perform a detailed convergence analysis along the example
trajectory of the exact flow $\Phi_{\gauss,*}$ shown in
\cref{fig:lgauss:phaseportrait}(a): $y_*=\Phi_{\gauss,*}(\deltat;x)$,
where $x=(1,0.5,2)^T$ and $\deltat=0.1$ (thus, $y_*\approx(1,
0.8459, 6.4556)^T$).
\begin{figure}[ht]
  \centering
  \includegraphics[width=1\textwidth]{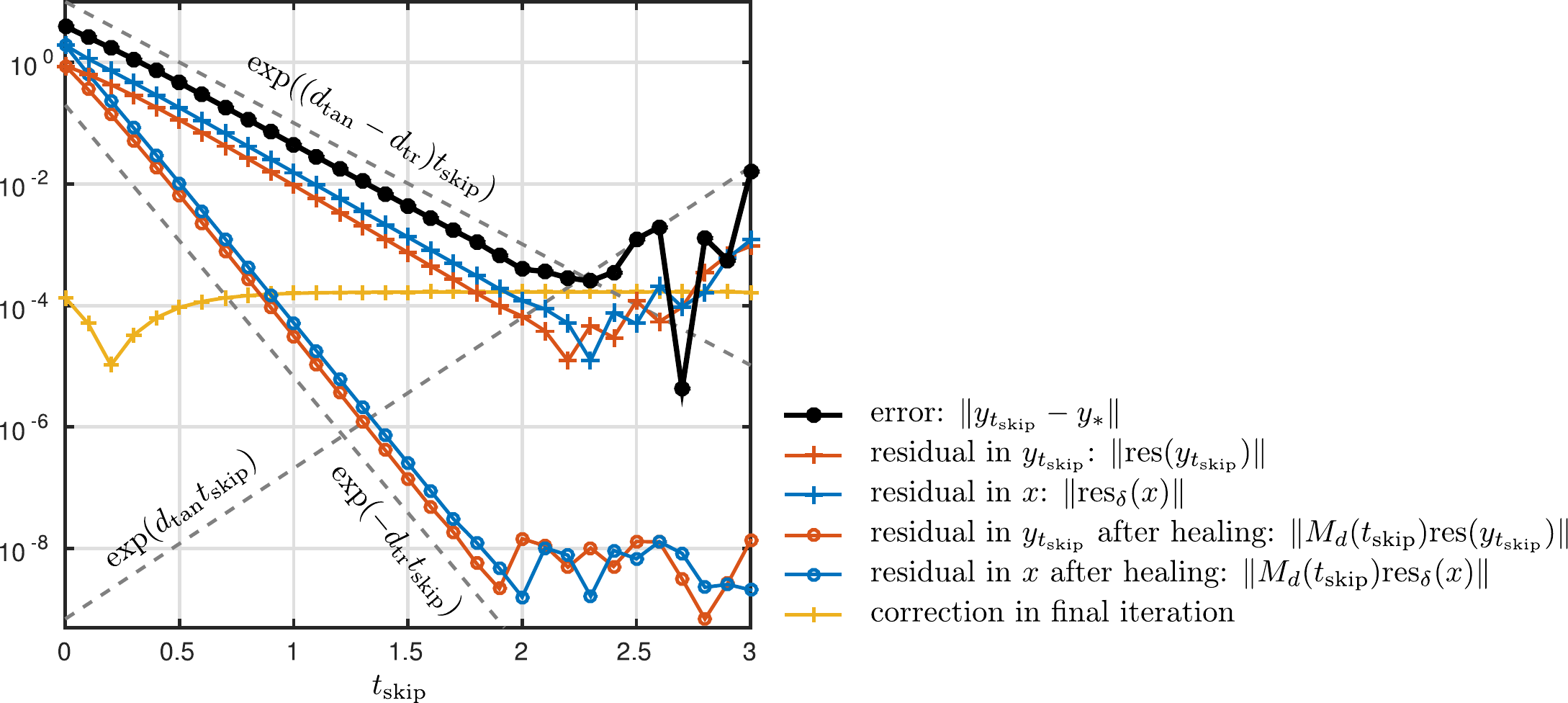}
  \caption{Convergence analysis along trajectory
    $y_*=\Phi_{\gauss,*}(\deltat;x)$, shown in
    \cref{fig:lgauss:phaseportrait}\,\textup{(}a\textup{)}. Parameters\textup{:}
    $\mu=6$, $\nu=0.3$, $\deltat=0.1$, $x=(1,0.5,2)^T$.}
  \label{fig:lgaussconv}
\end{figure}
To understand the factors entering the practically achievable lower
limit of the error
$\|y_\tskip-y_*\|=
\|\Phi_{\gauss,\tskip}(\deltat;x)-\Phi_{\gauss,*}(\deltat;x)\|$,
we consider again the identity \cref{eq:bfix} used in the proof of
\cref{thm:conv}:
\begin{multline}
  P_*(\tskip;y_\tskip)=P_*(\tskip;y_*)+\\
  [P_*(\tskip;y_\tskip)-P(\tskip;y_\tskip)]+
  [P(\tskip+\deltat;x)-P_*(\tskip+\deltat;x)]\mbox{,}  
\end{multline}
but re-arrange it using the concrete expressions for $P_{\gauss,*}$
and $P_\gauss$:
\begin{equation}\label{eq:lgauss:error:rec}
  \begin{aligned}
    R_d M_d(\tskip)T_\gauss(y_\tskip)
    =&\ R_d M_d(\tskip)T_\gauss(y_*)\ldots\\
    &+[R_dM_d(\tskip)T_\gauss(y_\tskip)-\restrict
    M(\tskip)\lgauss(y_\tskip)]\\
    &+[\restrict
    M(\tskip+\deltat)\lgauss(x)-R_dM_d(\tskip+\deltat)T_\gauss(x)]
  \end{aligned}
\end{equation}
Since the matrices $R_d$ and $M_d(\tskip)$ are invertible, we can apply
their inverses to both sides in \cref{eq:lgauss:error:rec}. For a
general distribution $\rho$, the composition of $R_d^{-1}$ and $\restrict$
\begin{align*}
  T_R \rho&:\Lint^2_1(\R;\R)\mapsto\R^3 &T_R\rho&=R_d^{-1}\restrict \rho
  =R_d^{-1}\left[\int_\R Q^{k-1}\rho(Q)\d Q\right]_{k=1,2,3}
\end{align*}
is a projection onto the slow manifold ${\cal C}$ in the coordinates
$(\phi_1,\phi_2,\phi_3)$. Furthermore, the nonlinear map $T_\gauss$ is
locally invertible in $y_*$ (and, hence, also in $y_\tskip$, if
$y_\tskip$ is near $y_*$). Its Jacobian is invertible in $y_*$ with a
moderate norm of its inverse $\|[\partial T_\gauss(y_*)]^{-1}\|\approx 10$ for
the chosen $y_*$. Hence, the identity \cref{eq:lgauss:error:rec} can
be written in the form
\begin{align}\nonumber
  T_\gauss(y_\tskip)=&\ T_\gauss(y_*)+\ldots\\
  \label{eq:lgauss:err}
  &\underbrace{T_\gauss(y_\tskip)-M_d(-\tskip)T_R M(\tskip)\lgauss(y_\tskip)}_{\displaystyle\res(y_\tskip)}+\ldots \\
  \nonumber
  &\underbrace{M_d(-\tskip)T_R
  M(\tskip+\deltat)\lgauss(x)-M_d(\deltat)T_\gauss(x)}_{\displaystyle\res_\deltat(x)}\mbox{.}
\end{align}
The two residual terms on the right-hand side, labelled
$\res(y_\tskip)$ and $\res_\deltat(x)$, are the two contributions to
the error, before it gets amplified by a moderate factor ($\|[\partial
T_\gauss(y_*)]^{-1}\|\approx 10$) when inverting $T_\gauss$. The
spectral properties of the flow $M$ ensure that
\begin{displaymath}
M_d(t)T_\gauss(\eta)-T_RM(t)\lgauss(\eta)\sim \exp(-\dtr t)\mbox{,}
\end{displaymath}
where $\dtr=-\lambda_4\approx10.3$.  Applying this estimate to
$\eta=y_\tskip$ and $t=\tskip$, and to $\eta=x$ and $t=\tskip+\delta$
gives the asymptotics $\sim\exp(-\dtr \tskip)$ in $\tskip$ for
$M_d(\tskip)\res(y_\tskip)$ and $M_d(\tskip)\res_\delta(x)$, shown in
\cref{fig:lgaussconv} (red and blue curves with circles). The healed
residuals $M_d(\tskip)\res(y_\tskip)$ and $M_d(\tskip)\res_\deltat(x)$
indeed decay with rate $\dtr$ until computational errors for computing
the distributions dominate (in this case $10^{-8}$). The matrix
$M_d(\tskip)^{-1}=M_d(-\tskip)$ has norm of order $\exp(\dtan\tskip)$
(where $\dtan=-\lambda_3\approx5.71$; see grey dashed line sloping
upward in \cref{fig:lgaussconv}) such that the residuals
$\res(y_\tskip)$ and $\res_\deltat(x)$ are of order
$\sim\exp((\dtan-\dtr)\tskip)$ (blue and red curves with $+$ marks in
\cref{fig:lgaussconv}). The residuals indeed decrease with rate
$\dtr-\dtan$ for increasing $\tskip$ until the amplification of the
computational errors by $\exp(\dtan\tskip)$ starts to dominate (at
$\tskip\approx2$). The true error $y_\tskip-y_*$ (shown in black in
\cref{fig:lgaussconv}) is then amplified approximately by the norm of
$\|[\partial T_\gauss(y_*)]^{-1}\|\approx 10$, because the residuals
$\res$ and $\res_\deltat$ occur on the manifold ${\cal C}$ (in the
coordinates $(\phi_1,\phi_2,\phi_3)$), while the error $y_\tskip-y_*$
is defined in $\dom\lift$. The relation between the error $y_\tskip-y_*$ and
the residual errors is independent of $\tskip$. Overall, the error
$y_\tskip-y_*$ decays with rate $\dtr-\dtan$ asymptotically for increasing
$\tskip$, but the computational error grows with rate $\dtan$. The
optimal healing time $\tskip$ is when both errors are of the same order
of magnitude.

The identity \cref{eq:lgauss:err} becomes a
nonlinear fixed-point problem after applying $T_\gauss^{-1}$, for
which the right-hand side is a contraction for sufficiently large
$\tskip$ (see the proof of \cref{thm:conv}). For
\cref{fig:lgaussconv} we applied this fixed-point iteration. The final
fixed-point iteration correction (shown as a yellow curve in
\cref{fig:lgaussconv} is always smaller than the
error $y_\tskip-y_*$.

\subsection{The size of 
  computational errors in ensemble computations}
The results shown in \cref{fig:plin:convergence,fig:lgaussconv} show
the qualitative behaviour of implicit lifting for increasing
$\tskip$. Two sources contribute to the overall error. One source is
the mismatch between the trajectory started from the lifted point and
the projected (along the stable fiber) trajectory on the slow
manifold. The size of this contribution is estimated in
\cref{thm:conv} as decaying with rate $\dtr-\dtan$ with
increasing $\tskip$ (also observed in
\cref{fig:plin:convergence,fig:lgaussconv}). The other source is the
limited accuracy in the computations of the lifting $\lift$, the
microscopic flow $M$ and the restriction $\restrict$. Errors
introduced from this limited accuracy grow with rate $\dtan$ for
increasing $\tskip$. The analysis in
\cref{fig:plin:convergence,fig:lgaussconv} illustrates the trade-off
between these two sources of error when the computational error is
small ($\approx10^{-8}$, using \texttt{chebfun} \cite{chebfun,Driscoll2014}).

If the microscopic flow $M$ describes a multi-particle or
high-dimensional stochastic system and is estimated using ensembles of
realizations then the computational error of the flow estimate (and,
possibly, the computation of $\lift$ and $\restrict$) is determined by
the ensemble size $N$. This error decreases asymptotically like
$1/\sqrt{N}$ for increasing $N$, unless one is able to apply variance
reduction techniques (see, for example, \cite{avitabile2014} for a
technique to reduce noise in the computations of Jacobians needed to
solve nonlinear systems). In this section we demonstrate that the
error behavior can be expected to be qualitatively the same as in
\cref{fig:plin:convergence,fig:lgaussconv}, but with stricter
limitations on $\tskip$ due to larger computational errors in $\lift$,
$\restrict$, and the flow $M$. To keep the computations simple and
comparable to the previous subsection, we perform a Monte-Carlo
simulation directly for the SDE~\cref{eq:sde_rescaled}.

\begin{figure}[t]
  \centering
    \subfigure[]{\label{fig:sde_error_a}\includegraphics[width=0.48\textwidth]{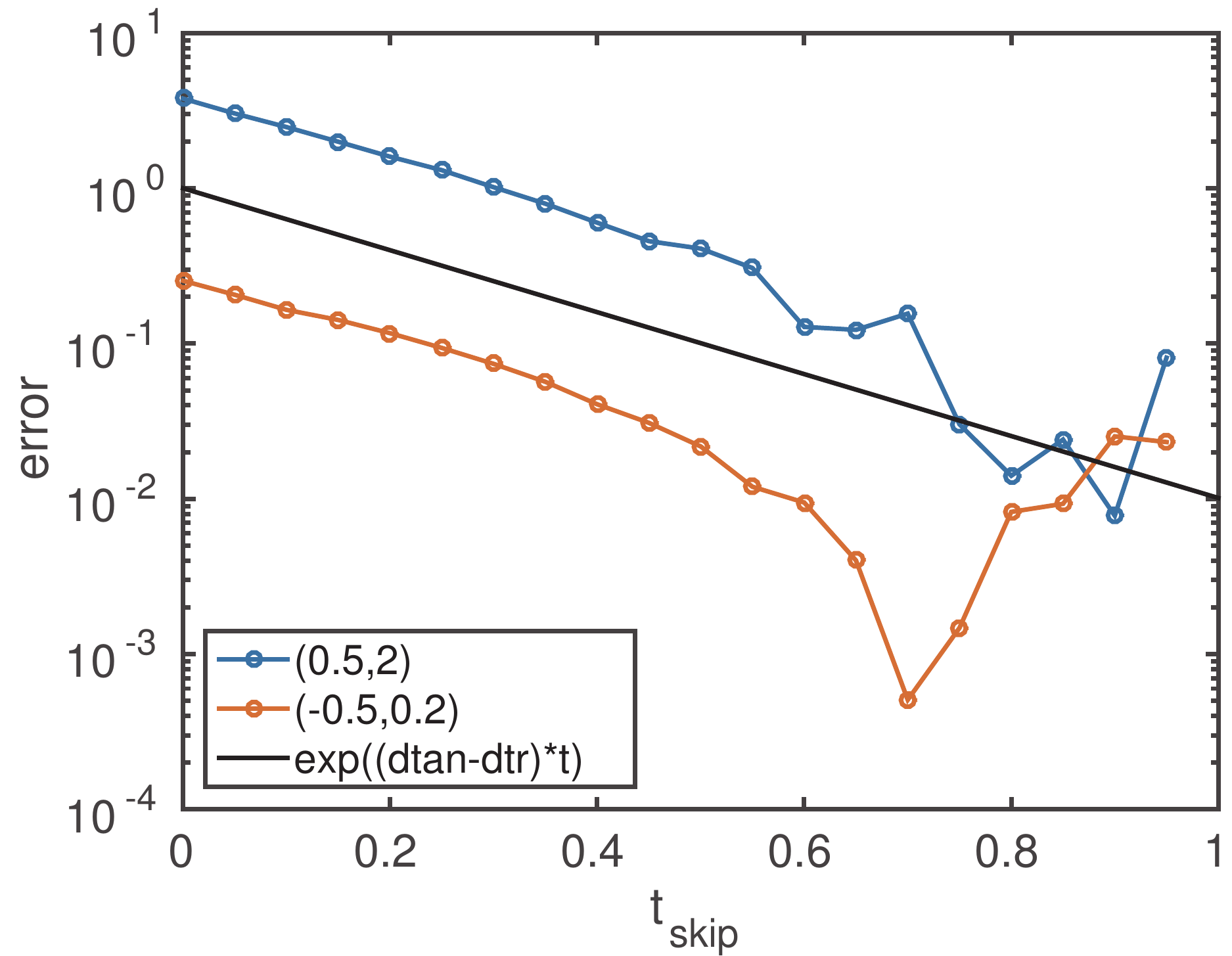}}\hfill
 \subfigure[]{\label{fig:fun_eval_error_b}\includegraphics[width = 0.48\textwidth]{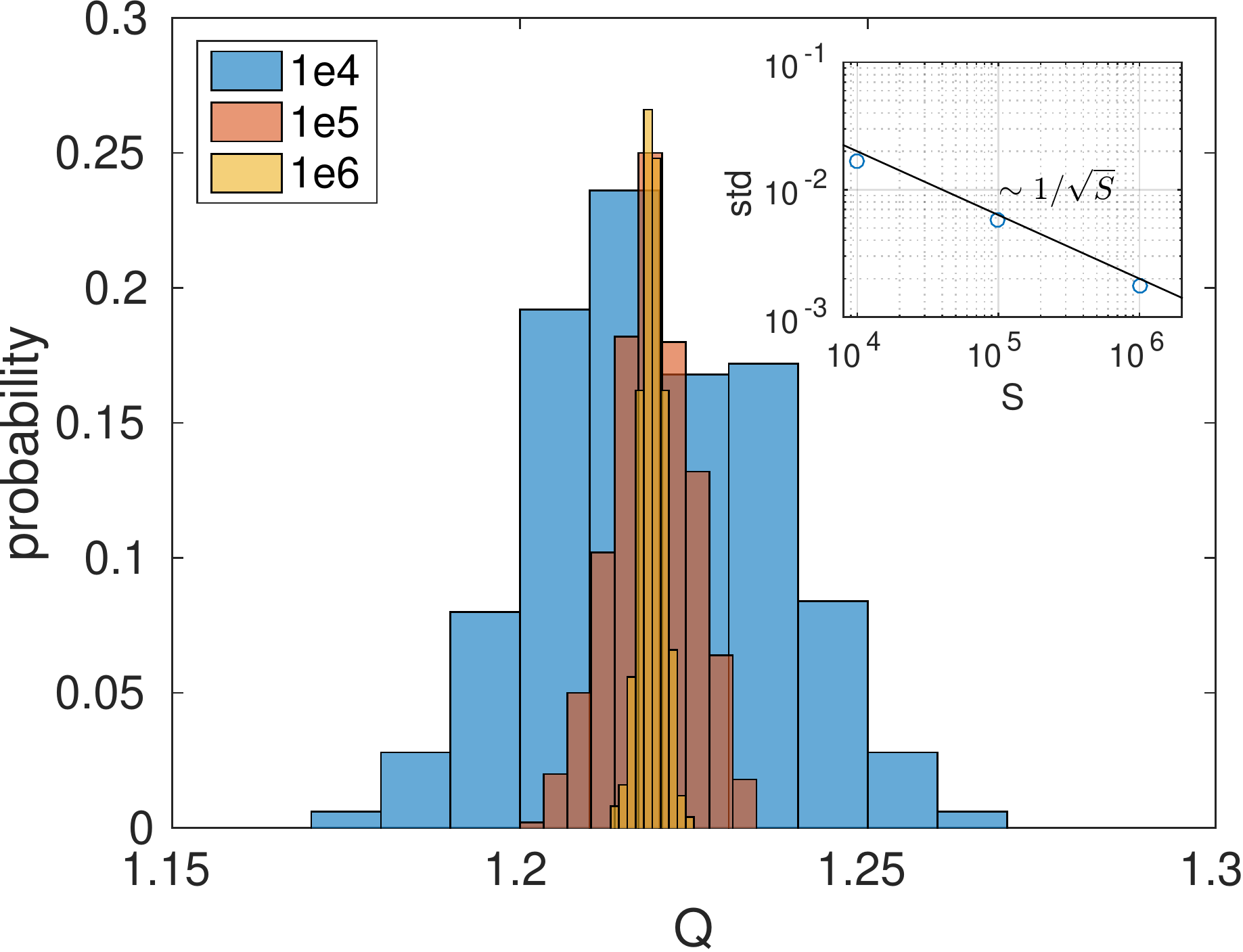}}
    \caption{\textup{(a)} Error analysis for $N=10^7, h = 10^{-2},
    \mu=6,\nu=0.3,\sigma = 1$, $\deltat=0.1$, similar to
    \cref{fig:mm_dyn_e_rot} for two different starting values $(N,\bar Q,
    \var Q)=(10^7,0.5,2)$ and $(10^7,-0.5,0.2)$.  The initial error is larger
    for the value with a larger variance. With increasing $\tskip$,
    the error shrinks with the same exponential rate for both initial
    conditions. \textup{(b)} Distribution of function
   evaluations of the lift-evolve-restrict map
   $\restrict(M(t;\lift (\bar Q, \var Q)))_2/N$ for $t=1$ and various
   ensemble sizes $N$ and $(\bar Q, \var Q)=(-0.5,0.2)$. Inset: the
   uncertainty of the function evaluation scales as $\thicksim
   1/\sqrt{N}$. 
 }
  \label{fig:sde_error}
\end{figure}

\Cref{fig:sde_error_a} shows the overall behaviour of the error
when performing computations based on random ensembles of finite size
$N$, using the lifting operator $\lgauss$, based on Gaussians. For an
ensemble size $N$, mean $\bar Q$ and we create a random set of initial
conditions
\begin{equation}
  \label{eq:sde_ensemble_lift}
  [\lift(N,\bar Q,\var Q)]_n = \bar Q + \sqrt{\var Q} \eta, \qquad n=1,\ldots,N
\end{equation}
where $\eta \thicksim \mathcal{N}(0,1)$ is a random variable drawn
from a standard normal distribution for each $n$. An ensemble of
$N$ realizations at positions $Q_n$ is restricted according to 
\begin{equation}
  \label{eq:sde_ensemble_macro}
  \restrict\left((Q_n)_{n=1}^N\right)_k= \sum_{n=1}^N Q_n^{k-1}\mbox{\quad ($k=1,2,3$).}
\end{equation}
Similar, to the definitions \cref{eq:sde_restrict} and
\cref{eq:sde_lgauss}, the first component of the argument to $\lift$
and of the output of $\restrict$ is the number of realizations, which
is preserved. In order to solve \cref{eq:sde_rescaled} numerically, we
use the Euler-Maruyama scheme
\begin{equation}
  \label{eq:euler-maruyama}
  Q_n(t+h) = Q_n(t) + f(Q_n(t)) h + \sqrt{h} \sigma \xi_n(t) 
  \mbox{\quad ($n=1,\ldots,N$),}
\end{equation}
where $h=0.01$ is the step size, $f=-\partial_Q V$ and
$\xi_n\thicksim\mathcal{N}(0,1)$ is standard normal random noise that is
uncorrelated, i.~e.\  $\langle\xi_n(t) \xi_n(t')\rangle = \delta(t-t')$.  

The error for each $\tskip$ in \cref{fig:sde_error_a} was estimated by
comparing the value of $\Phi_{\gauss,\tskip}(\delta;x)$ to the value
$\Phi_{\gauss,t_{\max}}(\delta;x)$ for the largest
$\tskip$ (called $t_{\max}$, equalling $1$). 
Thus, the
value of $\tskip$ at which the error starts to grow and the growth
rate may not have been captured accurately. However, we observe an
exponential decay with increasing $\tskip$ over approximately two
orders of magnitude and the more stringent limitation on $\tskip$, as
the error stops decreasing at $\tskip\approx0.6$.


Two problems limit the computational accuracy of function evaluations.
\begin{enumerate}
\item In Monte Carlo simulations with ensemble size $N$ the evaluation
  of the macroscopic lift-evolve-restrict map $P(t;\cdot)$ of the
  dynamics is noisy in $(\bar Q, \var Q)$.  This is due to the
  inherent noise in \cref{eq:sde_rescaled} and due to the noise in the
  lifting procedure \cref{eq:sde_ensemble_lift}.  Hence the evaluation
  of $P$ with the same input parameters might yield different
  outputs. The result of $P$ is a random variable with an
  ensemble-dependent distribution (see \cref{fig:fun_eval_error_b},
  where the distribution of the second component (the mean) of
  $P(1;(N,-0.5, 0.2))$ is shown for a range of $N$). The standard
  deviation of $P$ decreases with the ensemble size like $\thicksim
  1/\sqrt{N}$.
\item Function evaluations for large $\var Q$ become computationally
  difficult since a large $\var Q$ implies sampling of trajectories
  far away from the minima of the potential. Since the potential is
  steep away from the minima, the drift forces $V'$ become large,
  which results in stability problems of the numerical scheme
  \cref{eq:euler-maruyama} for a fixed step size $h$.
\end{enumerate}
When solving $P(\tskip;y)=P(\tskip+\deltat;x)$ for $y$ in the analysis
in \cref{fig:sde_error_a} we use a Newton iteration with damping
$\gamma = 0.5$ on the macroscopic level with tolerance $\texttt{tol} =
5\cdot 10^{-2}$ where Jacobians are computed by a central
finite-difference scheme with $\Delta \bar Q = \Delta \var Q =5\cdot
10^{-2}$. The ensemble size is $N=10^7$. 
The level of the minimal error is limited by
the finite ensemble size $N$ and the accuracy of function evaluations
and approximations of the Jacobian in the Newton iterations (see
\cite{avitabile2014} how the accuracy of the Jacobians can be
improved).

\section{Discussion}
\label{sec:discussion}

\subsection{General estimate for the influence of evaluation errors}
While the theoretical convergence result in \cref{thm:conv}
appears to suggest that a larger $\tskip$ always leads to a smaller
error, the demonstrations for the Michaelis-Menten kinetics model in
\cref{sec:exampl-mich-ment} and the SDE in
\cref{sec:sde-dynam} illustrate that there is a trade-off and,
hence, an optimal value for $\tskip$ in practice.  One source for the
difference between the estimates of \cref{thm:conv} and
numerical observations are numerical errors in the evaluation of
lifting $\lift$, evolution $M(t;\cdot)$ and restriction $\restrict$.
The effect of these errors grow along trajectories inside the slow
manifold ${\cal C}$ if the vector field tangent to ${\cal C}$ has
non-zero expansion rates forward or backward in time. This becomes
clear when looking at the arguments in the proof of
\cref{thm:conv}. The approximate solution $y_\tskip$ is the
fixed point of the map (see Equation~\cref{eq:app:Ndef})
\begin{multline}\label{eq:Ndef}
y\mapsto
  P_*(\tskip;\cdot)^{-1} \Bigl(P_*(\tskip;y_*)+\\
    \left[P(\tskip;y)-P_*(\tskip;y)\right]+
    \left[P(\tskip+\deltat;x)-P_*(\tskip+\deltat;x)\right]\Bigr)\mbox{.}
\end{multline}
According to \cref{thm:conv},
$y_\tskip-y_*\sim\exp((\dtan-\dtr)\tskip)$, where
$\dtan$ is defined as $\max\{\dtan^+,\dtan^-\}$, the maximum of the forward
($\dtan^+$) and backward ($\dtan^-$) expansion rate of the flow
$M\vert_{\cal C}$ tangential to ${\cal C}$, and $\dtr$ is the rate of
attraction transversal to ${\cal C}$. However, if we take into account
evaluation errors, we have to distinguish between the exact and
approximate operators. That is, $P_*(t;\cdot)$ equals $\restrict\circ
M(t;\cdot)\vert_{\cal C}\circ g\circ\lift$ (recall that $g$ is the
stable fiber projection) and $P(t;\cdot)$ equals
$\restrict_\Delta\circ M_\Delta(t;\cdot)\circ\lift_\Delta$ where we
use the subscript $\Delta$ to indicate that the operator is affected
by small errors. For $\lift_\Delta$ and $\restrict_\Delta$ this means
simply that they are perturbations of $\lift$ and $\restrict$ of size
$\Delta$. The evaluation error in $M$ along trajectories in ${\cal C}$
causes errors of size
\begin{align*}
  \|M(t;\cdot)\vert_{\cal C}-M_\Delta(t;\cdot)\vert_{\cal C}\|&\sim \Delta\exp(\dtan^+t)
  &&\mbox{for $t\geq0$,}\\  
  \|M(t;\cdot)\vert_{\cal C}-M_\Delta(t;\cdot)\vert_{\cal C}\|&\sim \Delta\exp(-\dtan^-t)
  &&\mbox{for $t<0$.}  
\end{align*}
These errors in $\lift_\Delta$, $\restrict_\Delta$ and
$M_\Delta(t;\cdot)$ are all part of the term $P(\tskip+\deltat;x)$ in
\cref{eq:Ndef}) such that the error grows for increasing
$\tskip$ at the rate
\begin{align*}
  \|P(\tskip+\deltat;x)-P_*(\tskip+\deltat;x)\|&\sim
  \Delta\exp(\dtan^+\tskip) \mbox{,}
\end{align*}
which gets then amplified by the expansion rate of
$M(-\tskip;\cdot)\vert_{\cal C}$ when applying
$P_*(\tskip;\cdot)^{-1}$. Thus, there will be an error between the
exact fixed point of the map \cref{eq:Ndef} and the fixed point with
evaluation errors. This error is of order
$\Delta\exp((\dtan^++\dtan^-)\tskip)$, which is growing exponentially
in $\tskip$. This is visible in all computational results:
\begin{itemize}
\item In the Michaelis-Menten kinetics model in
  \cref{sec:exampl-mich-ment} the error $\Delta$ is of order
  $10^{-10}$ and $\dtan$ is of order $\epsilon$ (which is $10^{-2}$)
  such that the growth of the error with $\tskip$ is not visible in
  the range of $\tskip$ between $0$ and $30$ in
  \cref{fig:mm_dyn_e,fig:mm_dyn_e_rot}.
\item For the stochastic differential equation in
  \cref{sec:sde-dynam}, $\dtan^+$ is zero and
  $\dtan^-=-\lambda_3\approx5.71$. For
  \cref{fig:plin:convergence,fig:lgaussconv} we computed the evolution
  of densities directly using the Fokker-Plank equation and
  \texttt{chebfun} such that the evaluation error $\Delta$ is of the
  order $10^{-8}$ (visible as the lower bound on the residuals
  $\|P_\lin(\tskip)-P_{\lin,*}(\tskip)\|$ in
  \cref{fig:plin:convergence} and in the residuals after healing in
  \cref{fig:lgaussconv}). Thus, the overall influence of the
  evaluation error is of order $\Delta\exp(\tskip\dtan^-)$. The
  amplification factor reaches $\sim10^5$ for $\tskip=2$. In
  \cref{fig:plin:convergence} evaluation errors dominate only from
  $\tskip\approx3$, while in \cref{fig:lgaussconv} they dominate from
  $\tskip\approx 2.5$.
\item In \cref{fig:sde_error_a} the growth rate of the evaluation error
  is the same as in \cref{fig:lgaussconv}, but the basic evaluation
  error of a single time step of $M_\Delta(t,\cdot)$ and the lifting
  $\lift_\Delta$ is larger (as they are generated from ensembles):
  $\Delta\sim10^{-3.5}$ for ensemble size $N=10^7$. Thus, the effects of
  evaluation error start to dominate already for
  $\tskip\approx0.7$. With smaller, more realistic, ensemble sizes the
  restriction on $\tskip$ posed by evaluation errors will be even more
  severe. Since the necessary length of $\tskip$ to reduce the
  projection error $y_\tskip-y_*$ (from \cref{thm:conv}) is dictated
  by $\dtan^--\dtr$, we have a general approximate optimal healing
  time for positive evaluation errors $\Delta$ of the order
  \begin{displaymath}
    \tskip\sim \frac{-\log\Delta}{\dtan^++\dtr}\mbox{,}
  \end{displaymath}
resulting in an optimal error of the order
\begin{displaymath}
\max\left\{  \Delta\e^{(\dtan^++\dtan^-)\tskip},\e^{(\dtan^--\dtr)\tskip}\right\}\sim
\Delta^p\mbox{\quad with\quad}
p=\frac{\dtr-\dtan^-}{\dtr+\dtan^+}\mbox{.}
\end{displaymath}
In the limit of large time scale separation ($\dtan^\pm/\dtr\to0$) the
power $p$ of the error reaches $1$ and the optimal $\tskip$ is of
order $-\log \Delta/\dtr$.
\end{itemize}

\subsection{Consequences for equation-free analysis of stochastic
  systems}
\label{sec:SDE:highdim}
The lift-evolve-restrict map $P_\gauss(t;\cdot)$ in
\cref{sec:sde-dynam} reduced the SDE $\d Q=-V'(Q)+\sigma\d W_t$
(or, more precisely, its Fokker-Planck equation) to the slow manifold
(a linear subspace) spanned by its first $3$ modes. Barkley \emph{et
  al} \cite{Barkley2006} observed that the map $P_\gauss(t;\cdot)$
(called \emph{moment map} in \cite{Barkley2006}) is nonlinear and,
hence, suspected that the nonlinearity of $P_\gauss$ may be the object
of interest for nonlinear analysis (such as finding multiple
equilibria, bifurcations under parameter changes, etc). However, as
equation~\cref{eq:phigauss} shows, the exact flow map
$\Phi_{\gauss,*}(\deltat;\cdot)$ of the low-order moments is still a
nonlinear transformation (by $T_\gauss$) of a linear map such that
there is no nonlinear dynamic behaviour present. More precisely, the
phase portrait of the exact flow map $\Phi_{\gauss,*}(\deltat;\cdot)$
is topologically conjugate to the phase portrait of a linear
system. Since the approximate flow
$\Phi_{\gauss,\tskip}(\deltat;\cdot)$, computed with $P_\gauss$,
converges to $\Phi_{\gauss,*}(\deltat;\cdot)$ for $\tskip\to\infty$ we
do not expect nonlinear behavior for $\Phi_{\gauss,\tskip}$ either.

This raises the question what the natural nonlinearity of the
underlying system is in the case of equation-free methods
applied to stochastic systems.

\subsubsection{Artificial nonlinearity}
\begin{figure}[ht]
  \centering
  \includegraphics[width=0.8\textwidth]{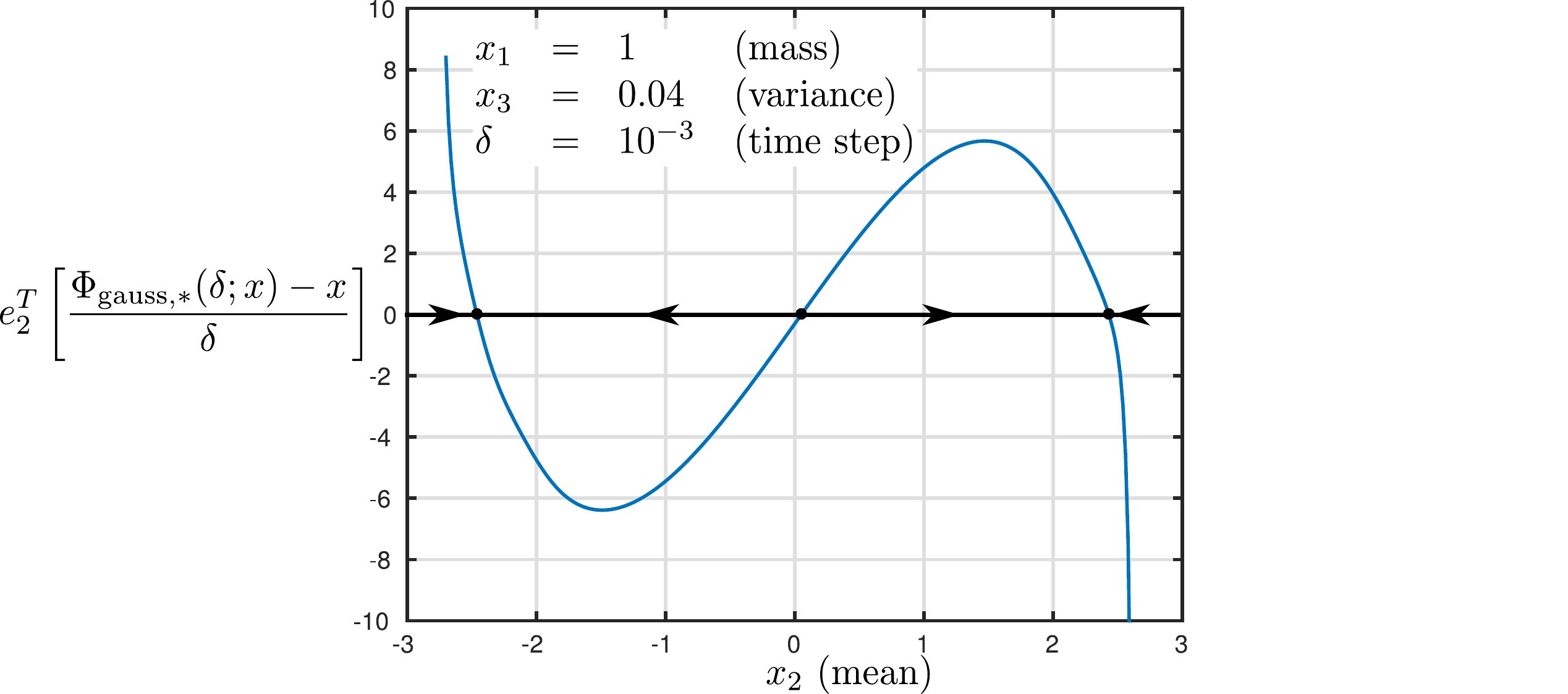}
  \caption{Apparent one-dimensional nonlinear phase portrait after
    projection of the nonlinearly transformed two-dimensional phase
    portrait for $\Phi_{\gauss,*}(\deltat;\cdot)$ in
    \cref{fig:lgauss:phaseportrait}(a) onto the horizontal line with
    variance $x_3=0.04$ \textup{(}similar to Figure 14 (left) of Barkley
    \emph{et al}. \cite{Barkley2006}\textup{)}. For healing time $\tskip$ with
    $1/\lambda_3\ll\tskip\ll1/\lambda_2(\approx10^{6})$, the
    approximate flow $\Phi_{\gauss,\tskip}$ approximates the exact
    flow $\Phi_{\gauss,*}$ accurately on the eigenspace of $L$
    corresponding to the $3$ dominant eigenvalues, but not on the
    smaller space for the $2$ dominant eigenvalues. Other
    parameters\textup{:} $\mu=6$, $\nu=0.3$, $\deltat=10^{-3}$,
    $x_1=1$, $x_3=4\times 10^{-2}$.}
  \label{fig:proj:phaseportrait}
\end{figure}
Since the Fokker-Planck equation is linear, the apparent nonlinear
dynamics arises only due to artificial projections of nonlinearly
transformed phase portraits of the linear Fokker-Planck equation when
the healing time $\tskip$ is not sufficiently large.  For example, let
us consider again the SDE with lifting to a Gaussian distribution from
\cref{sec:lgauss}. What happens if we choose a moment map for
only the zeroth and first moment but an insufficiently large $\tskip$
(which would have to be $\sim1/\lambda_2\approx10^6$ to make
\cref{thm:conv} applicable)? For illustration we choose a
lifting to near-delta Gaussian distributions, similar to
\cite{Barkley2006}. In the notation from \cref{sec:lgauss} this
means that we keep $x_1$ equal to $1$ (mass), vary $x_2$ (mean)
between $-3$ and $3$, and keep $x_3\ll1$ (variance) fixed ($x_3=0.04$
for the illustration in \cref{fig:proj:phaseportrait}). The restriction
is then the projection on the zeroth and first moment. If the healing
time $\tskip$ satisfies $1/\lambda_3\ll\tskip\ll1/\lambda_2$ (instead
of $\tskip\sim1/\lambda_2$), then we obtain for the approximate flow
$\Phi_{\gauss,\tskip}$ a projection of the phase portrait
\cref{fig:lgauss:phaseportrait}(a) onto the line with
$x_3=0.04$. \Cref{fig:proj:phaseportrait} shows this projected
phase portrait (arrows on the $x$-axis) and the associated right-hand
side (in blue). It resembles a phase portrait of a scalar ODE with two
coexisting stable fixed points, separated by an unstable fixed
point. Of course, this nonlinearity is created artificially by
projecting the accurate nonlinearly transformed two-dimensional phase
portrait of a linear system onto an arbitrarily chosen line in $\R^2$.


\subsubsection{Reduction of high-dimensional SDEs}
\label{sec:highdim:sde}
While in high-dimensional SDEs there is at first sight no obvious
nonlinearity present in the evolution of densities (see Fokker-Planck
equation~\cref{eq:fp_sde}), the reduction to low-order moments of a
multi-particle system with randomness still gives a valid dimension
reduction procedure. We give an informal outline of the argument for a
particularly simple case in which dimension reduction is in theory
possible according to Givon \emph{et al}. \cite{givon2004} (see also
textbook \cite{PS08}).

Let us assume that the simulation (say, an agent-based simulation) can
be modelled by a high-dimensional SDE (which is the microscopic model)
\begin{equation}\label{eq:SDEubig}
  \d u=F(u) \d t+\sigma_u\d W_{u,t}\mbox{,}
\end{equation}
where $u\in\R^{n_u}$ and (to keep the argument simple) $\sigma_u$ is
constant and regular, and $W_{u,t}$ are $n_u$ independent instances of
Brownian motion.  Let us also assume that there exist coordinates
$(x,y)\in\R^{n_x}\times\R^{n_y}$ ($n_x+n_y=n_u$) for $u$ such that in
these coordinates we have a time scale separation:
\begin{align}\label{eq:SDExy}
  \d x&=\epsilon f(x,y) \d t+\sqrt{\epsilon}\sigma_x\d W_{x,t}\mbox{,} &
  \d y&=g(x,y) \d t+\sigma_y\d W_{y,t}\mbox{,}  
\end{align}
and that for each $x$ the random variable $y$ converges to its
stationary density with rate of order $1$ (fast). 
Let $v_0(x,y)$ be the nullvector of the Fokker-Planck operator of the
fast subsystem of \cref{eq:SDExy}, $p\mapsto
L_0p=\partial_y[\frac{1}{2}\sigma_y^T\sigma_y \partial_yp- gp]$, with
$\int v_0(x,y)\d y=1$. Any function of the form $v_0(x,y)p_x(x)$ is
also a nullvector of $L_0$. If $(\epsilon\lambda,p)$ (with
$O(\lambda)=1$) is an eigenpair of the Fokker-Planck operator
$L_0+\epsilon L_1$ with
$L_1p=\partial_x[\frac{1}{2}\sigma_x^T\sigma_x \partial_xp- fp]$ for
the combined system \cref{eq:SDExy} in $(x,y)$ coordinates, then
$\lambda=\lambda_0+O(\epsilon)$, $p(x,y)=v_0(x,y)p_x(x)+O(\epsilon)$,
where $(\epsilon\lambda_0,p_x)$ is an eigenpair of the
of the right-hand side of the Fokker-Planck
equation for the reduced SDE
\begin{equation}\label{eq:SDEx:red}
  \d x=\epsilon \tilde f(x)\d t+\sqrt{\epsilon}\tilde \sigma_x(x)\d W_{x,t}\mbox{.}
\end{equation}
In \cref{eq:SDEx:red} $\tilde f(x)=\int f(x,y)v_0(x,y)\d y$ is the
conditional expectation with respect to $x$ of the drift in $x$ and
$\tilde \sigma_x(x)=\sigma_x\left[\int v_0(x,y)\d y/2\right]^{1/2}$ is
the standard deviation of $x$ in the stationary distribution of
$y$. 
Consequently, performing equation-free analysis on the
high-dimensional SDE~\cref{eq:SDExy} using a small number $d$ of
variables gives the same results as equation-free analysis on the
reduced system \cref{eq:SDEx:red} (up to order $\epsilon^2$).

Givon \emph{et al}. \cite{givon2004} discuss dimension reduction more
generally (independent of explicit spatial coordinates $x$ and $y$)
for Fokker-Planck operators of the form $L_0+\epsilon L_1$, assuming
that the linear operator 
$L_0$ has a non-trivial kernel (dimension greater than $1$,
implying that $\epsilon$ is a singular perturbation parameter). 
Hence, equation-free-analysis based on implicit lifting and
sufficiently large healing times can be used to perform
closure-on-demand, as described in \cite{Kevrekidis2010},
rigorously. Convergence of the approximate system created by
lift-evolve-restrict maps to the Fokker-Planck operator of the reduced
system \cref{eq:SDEx:red} occurs in the sense of classical singular
perturbation theory toward an attracting low-dimensional linear
invariant subspace of densities in the domain of definition of
$L_0+\epsilon L_1$, as ensured by \cref{thm:conv} for
sufficiently large healing times $\tskip$.

For the case that the high-dimensional SDE consists of a large number
$N$ of random variables (for example, describing agents) our analysis
in \cref{sec:sde-dynam} and the above discussion raise an
important point. Applying the equation-free procedure to initial
densities and the Fokker-Planck operator $L_0+\epsilon L_1$ does not
reduce the high-dimensional SDE to a low-dimensional SDE, but it
reduces the high-dimensional SDE to a low-dimensional linear ODE for
the coefficients of the leading modes of the Fokker-Planck
equation. Hence, increasing the number of variables $N$ (e.g., agents)
does \emph{not} increase the spectral gap or the time scale
separation. This is obvious for the simple SDE example in
\cref{sec:sde-dynam}: decreasing the noise level will let
$\lambda_2/\lambda_3$ converge to $0$ (the time scale for escape from
one well to the other), but $\lambda_3/\lambda_4$ will remain
approximately $1/2$. Hence, we need the convergence result for finite
time scale separation to prove validity of the model
reduction. Results for sufficiently large time scale separation such
as those by Zagaris \emph{et al}.
\cite{Gear2005,Zagaris2009,Zagaris2012} (using, for example,
constrained runs) and Marschler \emph{et al}. \cite{Marschler2014b}
are not applicable to equation-free methods operating on Fokker-Planck
equations, if the aim is to extract the decay rate or shape of the
dominant modes of the Fokker-Planck equation.

In summary, one possible work flow for analysing a high-dimensional SDE
with generator splittable as $L_0+\epsilon L_1$ with equation-free
methods is: (1) use the equation-free moment map to determine
properties of the leading $d$ eigenmodes $\phi_j$ and eigenvalues
$\lambda_j$ of $L_0+\epsilon L_1$; (2) if these $\phi_j$ and
$\lambda_j$ are also the leading eigenmodes and eigenvalues to an
operator $L_1$ for a Fokker-Planck equation of a low-dimensional SDE,
identify the properties of $L_1$ from the modes (for example, singular
points of the potential). 
\section{Outlook}
\label{sec:outlook}
The arguments in \cref{sec:sde-dynam}, studying the simple
scalar SDE $\d Q=-V'(Q)\d t+\d W_t$, and the discussion in
\cref{sec:SDE:highdim} treat SDEs as linear evolution equations for
densities. The sections below outline how one may have to modify the
arguments of \cref{thm:conv} for other tasks of equation-free
analysis, which are beyond the scope of this paper.

\subsection{Bifurcation analysis for the drift of the
  reduced system}
\label{sec:sderedbif}
Assume that we have access to a simulator of a system that can be
modelled by a high-dimensional SDEs of type \cref{eq:SDEubig},
\begin{equation}
\d
u=F(u)+\sigma_u\d W_{u,t}\mbox{,}\label{eq:SDEubig:out}
\end{equation}
with time-scale separation as in
\cref{eq:SDExy}. A sensible object for nonlinear equation-free
analysis is a bifurcation analysis of the deterministic part $\dot
x=\tilde f(x)$ of the reduced SDE \cref{eq:SDEx:red},
\begin{equation}
  \d x=\tilde f(x)\d t+\sigma_x\d W_{x,t}\mbox{,}\label{eq:SDEx:red:out}
\end{equation}
assuming a reduction as discussed in \cref{sec:highdim:sde} is
possible. For example, one may want to determine its phase portraits
and their parameter dependence. If one had a direct simulator of the
low-dimensional reduced SDE \cref{eq:SDEx:red:out}, one could
approximate $\tilde f$ in any given $x_0\in\R^{n_x}$ via
\begin{equation}\label{eq:driftapprox}
  \tilde  f(x_0)= \lim_{\delta\to0}\frac{1}{\delta}[E X_\delta-x_0]\mbox{,}
\end{equation}
where $X_\delta$ (a random variable in $\R^{n_x}$) is the solution of
the SDE \cref{eq:SDEx:red:out} at time $\delta$ starting from the
deterministic $x_0$, and $EX_\delta\in\R^{n_x}$ is its expectation.

Equation-free analysis based on a lift-evolve-restrict map $P$ with
healing time provides an approximation for \cref{eq:driftapprox} if
only a simulation of the high-dimensional SDE~\cref{eq:SDEubig:out} is
available. The healing time permits the fast variable $y$ to settle to
its stationary density $v_0(x,y)$ before one measures $\tilde
f$. Since the slow-fast coordinate split of $u$ into $x$ and $y$ is
unknown, one has to define a lifting $\lift$ and a restriction
$\restrict$ between $\R^{n_x}$ and the space of random variables $U$
in $\R^{n_u}$. 

Let us assume that the lift $\lift(x_L)$ of $x_L\in\R^{n_x}$ is a
random variable $U_0$ in $\R^{n_u}$ with density $p_0$ on
$\R^{n_u}$. The SDE \cref{eq:SDEubig:out} creates a Markov process
$t\mapsto U_t$ for $t\geq0$. Let us consider a restriction $\restrict$
of a random variable $U_t$ that is the expectation $E R(U_t)$ of a map
$R:\R^{n_u}\mapsto\R^{n_x}$. Thus, the lift-evolve-restrict
map $P:\R\times\R^{n_x}\mapsto\R^{n_x}$ is
$P(t;x_L)=E[R(U_t)|U_0=\lift(x_L)]$. A good approximation of the
deterministic part of the slow flow (in $x_L$ coordinates) would
\emph{not} be $(y-x_0)/\delta$ where $y$ is the solution of
$P(\tskip+\delta;x_0)=P(\tskip;y)$. Rather, a possible construction is
to define $x_R=P(\tskip;x_L)$ and then compute
\begin{equation}\label{eq:fcondexp}
  \tilde f_L(x_L)\approx
  \frac{1}{\delta}
  \left(E\left[R(U_{\tskip+\delta})\bigg\vert R(U_\tskip)=x_R\mbox{\ and\ } U_0=\lift(x_L)\right]-x_R\right)\mbox{.}
\end{equation}
This means that one first solves the SDE for the healing time
$\tskip$, then increases time to $\tskip+\delta$, and uses the
conditional expectation of $R(U_{\tskip+\delta})$, with the condition
that $R(U_\tskip)=x_R$. This conditional expectation enters the
difference quotient for $\tilde f_L(x_L)$, which is otherwise similar
to \cref{eq:driftapprox}. Constructions of the form \cref{eq:fcondexp}
do not fit into the framework of \cref{thm:conv}. Still, we conjecture
that the function $\tilde f_L$ approximates $\tilde f$ (up to a
coordinate change from $x_L$ to $x$) for sufficiently small $\delta$
and large $\tskip$. The approximation will become accurate only in the
limit of large time scale separation for a set of $n_x$ slow variables
(in contrast to \cref{thm:conv}), but we need only genericity
conditions on $\lift$ and $\restrict$.
\subsection{Averaging deterministic high-dimensional systems}
There is still another gap to applications for multi-particle systems,
which are commonly deterministic at the microscopic level. For
example, Barkley \emph{et al}. \cite{Barkley2006} used the scalar SDE
\cref{eq:sde_rescaled} as a simple model for a heat bath problem
where the position $Q$ of a heavy particle of mass $M$ and generalized
coordinates $(Q,P)$ is coupled to a heat bath of $N$ smaller particles
of masses $m_i$ and generalized coordinates $(q_i,p_i)$ for
$i=1,\ldots,N$. The full system in \cite{Barkley2006} was described by the Hamiltonian
\begin{equation}
  \label{eq:hamiltonian_hb}
  H(Q,P,q,p) = \frac{P^2}{2M} + V(Q) + \sum_{i=1}^N \frac{p_i^2}{2m_i}
  + \sum_{i=1}^N \frac{k_i}{2}(q_i - Q)^2\mbox{,}
\end{equation}
where the number $N$ of particles is large and the masses $m_i$ and
spring coupling constants $k_i$ are small (with particular
$N$-dependent distributions, see \cite{Barkley2006}, eq. (2.2)).  The
necessary assumption to enable treatment of a fast deterministic
subsystem as a stochastic system is some form of ergodicity: any
distribution of initial conditions of the fast subsystem converges
rapidly to a unique stationary distribution (conditioned on the slow
variables). This condition is hard to verify (even empirically) for
any particular system. In particular, it is not true for
\cref{eq:hamiltonian_hb} if one treats the coordinates $(Q,P)$ as the
slow variables since the small masses are only coupled through the
heavy particle. Convergence to an SDE is only guaranteed for the
system with Hamiltonian \cref{eq:hamiltonian_hb} if the initial
conditions for $q_i$ and $p_i$ are set according to the stationary
measure conditioned on $P$ and $Q$ (which was done in
\cite{Barkley2006}, see \cite{Barkley2006,KUPFERMAN2002,Pavliotis2003}
for background results). Hence, the introduction of a
healing time $\tskip$ will not have an improving effect for equation-free
analysis of the heat bath problem \cref{eq:hamiltonian_hb}.

\subsection{Approximation of stochastic slow manifolds}
\label{sec:stochslowmf}
As mentioned already in the introduction, our convergence result for
finite time scale separation relies on a result about model reduction
that is valid for finite time scale separation, namely the persistence
of normally hyperbolic invariant manifolds and their stable
fibers. While the model reduction results for stochastic systems in
\cite{givon2004,PS08} provide only statements for the limit of
infinite time scale separation, stronger results are available for
stochastic systems, if one is able to fix the noise realization (for
example, the Brownian path)
\cite{arnold1974stochastic,arnold2013random}. In this case, the
microscopic map $M$ has, for the example of an SDE of the form
$\d u=F(u) \d t+\sigma_u\d W_{u,t}$, the form $M(t;u,\omega)$, where
$\omega\in C([0,\infty);\R^D)$ is a realization of the Wiener process
$W_{u,t}$ and $M$ satisfies the invariance relation
$M(t+s;x,\omega)=M(t;M(s;x,\omega),\omega(s+\cdot))$.

Invariant stochastic manifolds ${\cal C}$ are then invariant objects
depending on the realization (one may write ${\cal C}(\omega)$).
Their persistence and attraction properties have been proven for some
cases such as finite-dimensional SDEs \cite{boxler1989,wang2013} and
SPDEs \cite{ChenRoberts2015}. For these cases, an implicit
equation-free scheme $y=\Phi_\tskip(\deltat;x,\omega)$ defined
implicitly via
\begin{equation}\label{eq:stochmf}
  \restrict M(\tskip;\lift(y),\omega)=    
  \restrict M(\tskip+\deltat;\lift(x),\omega)
\end{equation}
may converge in a similar way as claimed in
\cref{thm:conv}. However, the stochastic invariant manifold
results and the implementation of \cref{eq:stochmf} depend on the
ability to use the same realization $\omega$ throughout the
computation, as was done in \cite{kan2013} (for example, for different
arguments $y$ during a Newton iteration for \cref{eq:stochmf}). While
fixing the realization is possible for SDEs, for many of the
applications for equation-free analysis
\cite{GK08,KS09,MarschlerFaust2014,Marschler2015,Marschler2014,SGK12,TLS16}
it is not clear how to do that.

\section{Conclusion}
\label{sec:conclusion}

This paper proves convergence of equation-free methods, based on
lift-evolve-restrict maps $P(t;\cdot)=\restrict\cdot
M(t;\cdot)\circ\lift$. Our convergence proof does not assume that the
time scale separation becomes large, in contrast to previous results
\cite{Zagaris2012,Marschler2014b}. Rather, convergence is
achieved for finite time scale separation, but in the limit of large
healing time $\tskip$ and an implicit approximation of the slow flow
$\Phi_*(t;x)$: $P(\tskip;y)=P(t+\tskip;x)$ defines the approximation
$\Phi_\tskip(t;x):=y$. The original explicit equation-free framework,
as proposed by Kevrekidis \emph{et al}., corresponds to the case where
$\tskip = 0$ and $\restrict \circ \lift = \id$. The analysis is
performed for attracting slow manifolds in deterministic systems.
However, we demonstrate on a simple SDE that our result may also be
useful for stochastic systems, where the time scale separation is in
the spectrum of the Fokker-Planck equation and is often only of order
$1$. In particular, for the prototype example investigated by
\cite{Barkley2006} the implicit flow approximation $\Phi_\tskip$
converges to the true solution $\Phi_*$ of the linear Fokker-Planck
equation for large healing times $\tskip$.


\section*{Acknowledgements}
\label{sec:ack}
J.~Sieber's research was supported by funding from the European
Union's Horizon 2020 research and innovation programme under Grant
Agreement number 643073, by the EPSRC Centre for Predictive Modelling
in Healthcare (Grant Number EP/N014391/1) and by the EPSRC Fellowship
EP/N023544/1. 

C. Marschler and J. Starke would like to thank Civilingeni{\o}r Frederik
Christiansens Almennyttige Fond for financial support. J. Starke would
also like to thank the Villum Fonden (VKR-Centre of Excellence Ocean
Life), the Technical University of Denmark and Queen Mary University
of London for financial support.  

\bibliographystyle{siamplain}
\bibliography{extended_implicit}

\appendix
\section{Proof of Convergence \cref{thm:conv}}
\label{sec:app:conv}
For the proof of \cref{thm:conv} we have to analyze the two
equations (for $y$ and $y_*$ respectively)
\begin{align}\label{eq:app:perturbed}
  \restrict(M(\tskip;\lift(y)))&=
  \restrict(M(\tskip+\deltat;\lift(x)))\mbox{,}\\
    \label{eq:app:exact}
    \restrict(M(\tskip;g(\lift(y_*))))&=
    \restrict(M(\tskip+\deltat;g(\lift(x))))
    \mbox{.}
\end{align}
In both equations $x\in\R^d$ enters as a
parameter. \Cref{ass:transversality} ensures that the
solution $y_*$ of \cref{eq:app:exact} is unique and
independent of $\tskip$. For equation~\cref{eq:app:perturbed} we have
to prove the existence of a solution $y$, and prove that it is close
to $y_*$ for sufficiently large $\tskip$. Throughout this appendix we
will use the notations
\begin{displaymath}
  u(t)=O(\exp(\alpha t))\mbox{,\quad}
  v(t)=o(\exp(\alpha t))
\end{displaymath}
to describe that $\|u(t)\exp(-\alpha t)\|$ is bounded uniformly for
all $t\geq 0$, and that the function $v(t)\exp(-\alpha t)$ tends to
zero for $t\to\infty$. For the special case $\alpha=0$ we write $O(1)$
and $o(1)$. If the quantity depends also on other parameters (say,
$y\in\dom\lift$) then the expression implies uniformity (for example,
for $y$ close to $y_*$) unless stated explicitly otherwise.

Using the definitions \cref{eq:pstardef} of
$P_*(t;x)=\restrict(M(t;g(\lift(x))))$ and \cref{eq:pdef} for the map
$P(t;x)=\restrict(M(t;\lift(x)))$,
equation~\cref{eq:app:perturbed} can be written in the form (using
\cref{eq:app:exact})
\begin{align}
  \label{eq:app:ydef}
  P_*(\tskip;y)&=P_*(\tskip;y_*)+\exp(-\dtr\tskip)\left[G(\tskip;y)+H(\tskip;x)\right]
  \mbox{,\quad where}\\[1ex]
  G(t;y)&=-\exp(\dtr t)\left[P(t;y)-P_*(t;y)\right]\mbox{,}\nonumber\\
  H(t;x)&=\phantom{-}\exp(\dtr t)\left[P(t+\deltat;x)-P_*(t+\delta;x)\right]\mbox{.}\nonumber
\end{align}
The operator $P_*$ and the newly introduced $G$ and $H$ satisfy the
following conditions on their derivatives by
\cref{ass:timescales}, \cref{eq:mcebound} and
\cref{eq:gepsconv} on separation of time scales for the flow $M$:
\begin{align}
  \partial_2^jP_*(\tskip;y)&=O(\exp(\dtan \tskip))\mbox{,}\label{eq:app:pbound}\\
  \partial_2^j G(\tskip;y)&= O(1)\mbox{,}\label{eq:app:gbound}\\
  \partial_2^j H(\tskip;x)&=O(1)\label{eq:app:hbound}
\end{align}
for all $j\in\{0,\ldots,k_{\max}\}$ and all $y$ in a neighborhood of
$y_*$. In the case of $H$ the bound is also uniform for
$\deltat\in[-\delta_{\max},\delta_{\max}]$. Thus, the parameter
$\deltat$ has been dropped from the list of arguments in $H$.
Combining the separation of time scales in \cref{ass:timescales},
\cref{eq:mcebound}, with \cref{ass:transversality} on the uniform
invertibility of $\restrict\vert_{\cal C}$ and $g\circ \lift:\dom
\lift\mapsto {\cal C}$, we have a Lipschitz constant ($C$ is
independent of $y_1$, $y_2$ and $\tskip$)
\begin{align}
  \|y_1-y_2\|&\leq C\exp(\dtan
  \tskip)\|P_*(\tskip;y_1)-P_*(\tskip;y_2)\|
\mbox{,}\label{eq:app:pinvbound}
\end{align}
when inverting
$P_*(\tskip;\cdot)$ for all $y_1,y_2$ in a neighborhood of $y_*$ and
all $\tskip\geq0$. We also note that
\begin{equation}
  \label{eq:app:dystar}
  \left\|\frac{\partial^j y_*(x)}{\partial x^j}\right\|=O(1)\mbox{}
\end{equation}
Specifically, these derivatives depend only on
$\deltat\in[-\delta_{\min},\delta_{\max}]$. Thus, $\partial^j y_*(x)$
are uniformly bounded due to \cref{eq:mcebound}, and because we
required $\exp(\dtan\delta_{\max})=O(1)$.

\emph{Abbreviating notation} In the following all derivatives of the
functions $P_*$, $G$ and $H$ are with respect to their second argument
($y$ or $x$). The argument $\tskip$ enters the functions $P_*$, $G$
and $H$ as a parameter that we will drop in our notation such that we
will write, for example, $\partial^3 P_*(y_*)[\partial
y_*]^2[\partial^2y_*]$ for $\partial_2^3P_*(\tskip;y_*)[\partial y_*/\partial
x]^2[\partial^2y_*/(\partial x)^2]$. The parameter $\tskip$ enters
estimates via the bounds \cref{eq:app:pbound}--\cref{eq:app:dystar}
for $P_*$, $G$ and $H$.

The properties \cref{eq:app:pbound}--\cref{eq:app:dystar} make
Banach's contraction mapping principle applicable to
Equation~\cref{eq:app:ydef} in a sufficiently small neighborhood of
$y_*$ and for sufficiently large $\tskip$ (as shown in the paragraph
that follows). We then estimate the error of the derivatives of $y$
with respect to $x$.

\paragraph{Existence of solution $y$ and its error} We apply the
Banach contraction mapping principle to the map 
\begin{equation}\label{eq:app:Ndef}
  N:y\mapsto
  P_*^{-1} \left(
    P_*(y_*)+\exp(-\dtr\tskip)
    \left[G(y)+H(x)\right]\right)
\end{equation}
($P_*^{-1}(\cdot)$ is
the inverse of the diffeomorphism $P_*:U(y_*)\mapsto U(P_*(y_*))$).
Let $B$ be a closed ball around $y_*$ of radius $R$ in which all
estimates \cref{eq:app:pbound}--\cref{eq:app:pinvbound} on $P_*$,
$G$ and $H$ hold. Combining the estimate \cref{eq:app:pinvbound} for
the Lipschitz constant of $P_*^{-1}$ with $y_1=y$ and
$y_2=y_*$, and the bound on the derivatives for $G$ (w.r.t. $y$) gives
an estimate for the difference of $N(y)$ from $y_*$:
\begin{align*}
  \|N(y)-y_*\|&\leq C\exp((\dtan-\dtr)\tskip)
  \left[\max_B\|\partial G\|\|y\|+\|H(x)\|\right]\mbox{,}\\
  &\leq C\exp((\dtan-\dtr)\tskip) \left[\max_B\|\partial
    G\|(\|y_*\|+R)+\|H(x)\|\right]\mbox{.}
\end{align*}
Thus, choosing $\tskip$ sufficiently large, we can ensure that $N$
maps $B$ back into itself (since $\dtan<\dtr$).
Similarly, the Lipschitz constant of $N$ in $B$ can be estimated by
\begin{align*}
  \|N(y_1)-N(y_2)\|\leq C\exp((\dtan-\dtr)\tskip)\max_B\|\partial
  G\|\|y_1-y_2\|\mbox{,}
\end{align*}
where $C\exp((\dtan-\dtr)\tskip)\max_B\|\partial G\|$ is smaller than
unity for sufficiently large $\tskip$. Consequently, $N$ has a unique
fixed point $y$ in $B$, which solves the perturbed problem
\cref{eq:app:perturbed}. Moreover, the difference $y-y_*$ satisfies
\begin{equation}
  y-y_*=O(\exp((\dtan-\dtr)\tskip))\mbox{.}\label{eq:app:zbound}
\end{equation}

\paragraph{Error of derivatives}
The smoothness of the coefficients in \cref{eq:app:ydef} ensures that
$y$ is also differentiable as a function of $x$ up to order
$k_{\max}$.  We want to prove that for $\ell$ satisfying $\ell\leq
k_{\max}-1$ (where $k_{\max}$ is the order of differentiability of the
coefficients in \cref{eq:app:ydef}) and $(2\ell+1)\dtan<\dtr$ the
bound on the error is
\begin{equation}
  \label{eq:app:dyerr}
  \partial^{\ell}y-\partial^{\ell}y_*=O(\exp(((2\ell+1)\dtan-\dtr)\tskip))\mbox{.}
\end{equation}
We prove this by induction starting from $\ell=1$, which we check
first using the previous paragraph's results.

Assume that the bound \cref{eq:app:dyerr} holds for all derivatives
up to $\ell-1$. This implies, in combination with
\cref{eq:app:dystar}, that $y$, $\partial y$, \ldots,
$\partial^{\ell-1}y$ are bounded uniformly for all $\tskip\geq0$ (just
like $\partial^\ell y_*$ for $\ell=1\ldots k_{\max}$ by
\cref{eq:app:dystar}). In order to estimate the difference
$\partial^\ell y-\partial^\ell y_*$, we return to \cref{eq:app:ydef}
and differentiate each of the terms $\ell$ times with respect to $x$
(noting that $y_*$ and $y$ are also functions of $x$):
\begin{equation}\label{eq:app:ydiff}
  \frac{\partial^\ell}{\partial x^\ell}[P_*(y(x))]-\frac{\partial^\ell}{\partial x^\ell}[P_*(y_*(x))]=\exp(-\dtr\tskip)    \left[\frac{\partial^\ell}{\partial x^\ell}[G(y(x))]+\partial^\ell H(x)\right]\mbox{.}
\end{equation}
The term $\partial^\ell H(x)$ is $O(1)$ for all $\tskip\geq0$ by
\cref{eq:app:hbound}. In the term $\partial^{\ell}/(\partial
x^{\ell})[G(y)]$ we extract the highest-order derivative of $y$ by
writing it in the form
\begin{align}\label{eq:app:dghbound}\nonumber
  \frac{\partial^{\ell}}{\partial x^{\ell}}[G(y)]&= O(1)+
  \partial G(y)\partial^{\ell}y=O(1)+\partial G(y)\partial^\ell y_*+
  \partial G(y)\left[\partial^\ell y-\partial^\ell y_*\right]\\
  &=O(1)+\partial G(y)\left[\partial^\ell y-\partial^\ell y_*\right]
\end{align}
For \cref{eq:app:dghbound} the boundedness of the $O(1)$ terms
follows from the boundedness of all their parts: the derivatives of
$G$ are bounded by \cref{eq:app:gbound}, $\partial^\ell y_*$ is
bounded by \cref{eq:app:dystar}, and $y$, $\partial y$,\ldots,
$\partial^{\ell-1}y$ are bounded by induction hypothesis. The
pre-factor $\partial G(y)$ of $\partial^{\ell}y-\partial^\ell y_*$ is
also bounded uniformly for all $\tskip\geq0$.

Inserting the right-hand side of \cref{eq:app:dghbound} into the
right-hand side of \cref{eq:app:ydiff}, we obtain
\begin{equation}\label{eq:app:ydiff:lhs}
  \frac{\partial^\ell}{\partial x^\ell}[P_*(y(x))]-\frac{\partial^\ell}{\partial x^\ell}[P_*(y_*(x))]=
  \exp(-\dtr\tskip)\partial G(y)
  \left[\partial^{\ell}y- \partial^{\ell}y_*\right]+O(\exp(-\dtr\tskip))
\end{equation}
Expanding the left-hand side of the above equation using the chain
rule, we get a sequence of differences with equal powers of
derivatives of $P_*$, $y$ and $y_*$. From this sequence of differences
we extract the difference between derivatives involving $\partial^\ell
y$ and $\partial^\ell y_*$ and collect all other terms in a remainder
$r$ (which is present only for $\ell>1$ and will later turn out to be
of order $O(\exp((2\ell\dtan-\dtr)\tskip))$):
\begin{equation}\label{eq:app:dysplit}
  \frac{\partial^\ell}{\partial x^\ell}[P_*(y(x))]-\frac{\partial^\ell}{\partial x^\ell}[P_*(y_*(x))]=
  \partial P_*(y)\partial^{\ell}y-
  \partial P_*(y_*)\partial^{\ell}y_*+r\mbox{.}
\end{equation}
From the difference with the highest-order derivatives of $y$ and
$y_*$ we extract the difference $\partial^\ell y-\partial^\ell y_*$ by
adding zeroes. Using the notational convention
\begin{displaymath}
F\{x,y\}=\int_0^1F(sx+(1-s)y)\d s
\end{displaymath}
for the mean between two points of a single-argument function $F$ in
the following,
\begin{align}
\partial P_*(y)\partial^{\ell}y-
\partial P_*(y_*)\partial^{\ell}y_*
=\ &\partial P_*(y)\left[\partial^{\ell}y-\partial^{\ell}y_*\right]+\left[\partial P_*(y)-\partial P_*(y_*)\right]\partial^{\ell}y_*\nonumber\\
=\ &\partial P_*(y)\left[\partial^{\ell}y-\partial^{\ell}y_*\right]
+\partial^2P_*\{y,y_*\}
[y-y_*]\,\partial^{\ell}y_*\label{eq:app:dly:split}\\
=\ &\partial P_*(y)\left[\partial^{\ell}y-\partial^{\ell}y_*\right]+
O(\exp((2\dtan-\dtr)\tskip))\mbox{.}\label{eq:app:dly}
\end{align}
The order $O(\exp((2\dtan-\dtr)\tskip))$ of the second term follows
from the bounds on $y-y_*$ (given in \cref{eq:app:zbound}),
$\partial^2P_*$ (given in \cref{eq:app:pbound}) and the boundedness
of $\partial^\ell y_*$ (given in \cref{eq:app:dystar}). This
immediately implies the estimate for the case $\ell=1$: inserting
\cref{eq:app:dly} into \cref{eq:app:ydiff:lhs}, we have for $\ell=1$
\begin{equation}
  \label{eq:app:dy1bound}
  \partial P_*(y)\left[\partial y-\partial y_*\right]
  =\exp(-\dtr\tskip) \partial G(y)
  \left[\partial y-\partial y_*\right]
  +O(\exp((2\dtan-\dtr)\tskip))\mbox{.}
\end{equation}
In \cref{eq:app:dy1bound} we have collected the bounded terms with
pre-factors $\exp(-\dtr\tskip)$ and $\exp((2\dtan-\dtr)\tskip)$ using
the larger pre-factor $\exp((2\dtan-\dtr)\tskip)$.  Since (by
\cref{eq:app:pinvbound}) the inverse of $\partial P_*(y)$ satisfies
$\partial P_*(y)^{-1}=O(\exp(\dtan \tskip))$ we can rearrange
\cref{eq:app:dy1bound} to isolate $\partial y-\partial y_*$ for large
$\tskip$, giving the estimate (note that $\partial G(y)=O(1)$)
\begin{equation}\label{eq:app:dyerrbound}
  \partial y-\partial y_*=O(\exp((3\dtan-\dtr)\tskip))\mbox{,}
\end{equation}
which is what we had to prove for $\ell=1$.

\paragraph{Error of higher-order derivatives}
Let us assume that the assumptions of the theorem are satisfied for
all $j<\ell$ with $\ell\geq 2$. By the conditions of the theorem we
assume that $(2\ell+1)\dtan<\dtr$ and the conditions
\cref{eq:app:pbound}--\cref{eq:app:dystar} are satisfied for
$j\leq\ell$ (including existence of the corresponding derivatives).

For $\ell>1$ we have to include the remainder $r$ from
\cref{eq:app:dysplit} in our consideration. This remainder is a sum
of expressions $a_\nu$ of the form
\begin{equation}\label{eq:app:faadibruno:rem}
  a_\nu=\partial^jP_*(y)
    \left[\partial^{\nu_1} y\right]\ldots
    \left[\partial^{\nu_j} y\right]-
  \partial^jP_*(y_*)
  \left[\partial^{\nu_1} y_*\right]\ldots
  \left[\partial^{\nu_j} y_*\right]\mbox{,}
\end{equation}
where $2\leq j\leq\ell$, and $\nu$ is a $j$-tuple of integers
$\nu_i\in\{1,\ldots,\ell-1\}$ with $\sum_{i=1}^{j}\nu_i=\ell$.  All
factors $\partial^{\nu_i}y$ and $\partial^{\nu_i}y_*$ are of order
$O(1)$ with respect to $\tskip$ according to \cref{eq:app:dystar} and
induction hypothesis. The terms $\partial^jP_*(y)$ and
$\partial^jP_*(y_*)$ are of order $O(\exp(\dtan\tskip))$ according to
\cref{eq:app:pbound}. The difference in \cref{eq:app:faadibruno:rem}
can be expressed as a sum of $j+1$ differences involving
$\partial^iy-\partial^iy_*$ for some $i\in\{0\ldots,\ell-1\}$ by adding $j+1$
zeros:
\begin{align}\label{eq:app:anu:p}
  a_\nu=&\partial^{j+1}P_*\{y,y_*\}[y-y_*]     
  \left[\partial^{\nu_1} y\right]\ldots
    \left[\partial^{\nu_j} y\right]\\
  &+\sum_{i=1}^{j}\partial^jP_*(y_*)
  \left[\prod_{m<i}\partial^{\nu_m}y_*\right]
  \left[\partial^{\nu_i}y-\partial^{\nu_i}y_*\right]
  \left[\prod_{m>i}\partial^{\nu_m}y\right]
  \label{eq:app:anu:y}
\end{align}
The right-hand side in \cref{eq:app:anu:p} is of order
$O(\exp((2\dtan-\dtr)\tskip))$. The $i$th term in the sum in
\cref{eq:app:anu:y} is of order
$O(\exp((\dtan(1+(2\nu_i+1))-\dtr)\tskip))$. So, since $\nu_i\leq
\ell-1$ and $\ell>1$, all terms in the sum for $a_\nu$ are at most of
order $O(\exp((2\ell\dtan-\dtr)\tskip))$. Consequently,
\begin{equation}
  \label{eq:app:r}
  r=O(\exp((2\ell\dtan-\dtr)\tskip))\mbox{.}
\end{equation}
Inserting this estimate in combination with \cref{eq:app:dysplit} and \cref{eq:app:dly} into \cref{eq:app:ydiff:lhs}, we obtain
\begin{equation}\label{eq:app:hdyerr}
  \partial P_*(y)\left[\partial^{\ell}y-\partial^{\ell}y_*\right]=
  \exp(-\dtr\tskip)\partial G(y)
  \left[\partial^{\ell}y- \partial^{\ell}y_*\right]+
  O(\exp((2\ell\dtan-\dtr)\tskip))\mbox{.}
\end{equation}
In \cref{eq:app:hdyerr} we have included the smaller error terms
$O(\exp((2\dtan-\dtr)\tskip))$ and $O(\exp(-\dtr\tskip))$ into the
(for $\ell>1$) larger $O(\exp((2\ell\dtan-\dtr)\tskip))$. Since,
$\dtr<\dtan$, $\partial G(y)=O(1)$ and $\partial
P_*(y)=O(\exp(\dtan\tskip))$, we can isolate
$\partial^{\ell}y-\partial^{\ell}y_*$ in \cref{eq:app:hdyerr}.  This
results in the asymptotic estimate claimed in \cref{thm:conv}:
\begin{displaymath}
  \partial^{\ell}y-\partial^{\ell}y_*=
  O(\exp(((2\ell+1)\dtan-\dtr)\tskip))\mbox{.}
\end{displaymath}

\section{Brief description of supplementary material}
\label{sec:supp}
The supplementary material contains \texttt{Matlab}/\texttt{octave}
scripts and functions that reproduce \cref{fig:mm_dyn} from
\cref{sec:exampl-mich-ment}. The provided \texttt{zip} file
unpacks into folder \code{demo\_Michaelis\_Menten/}. The main script
is \code{demo\_Michaelis\_Menten.m}, which will reproduce
\cref{fig:mm_dyn}, showing phase space geometry of the Michaelis-Menten
kinetics \cref{eq:mmdyn} with explicit time scale separation as also
studied by Gear \emph{et al}. and others
\cite{OMalley1991,Gear2005,Zagaris2009,Zagaris2012}.
\begin{itemize}
\item Folder \code{demo\_Michaelis\_Menten/rotated/} contains the
  published \texttt{html} output from the script for the rotated
  coordinate system \cref{eq:rot_dyn}  in file
  \code{demo\_Michaelis\_Menten.html}.
\item Folder \texttt{demo\_Michaelis\_Menten/unrotated/} contains the
  published \texttt{html} output from the script for the coordinate
  system with explicit time scale separation \cref{eq:mmdyn} in file
  \texttt{demo\_Michaelis\_Menten.html}.
\item Folder \texttt{tools/} contains some auxiliary functions called
  in the script (a simple Newton iteration \texttt{ScSolve.m}, an
  explicit initial-value-problem solver using the Dormand-Prince
  scheme and fixed step size \texttt{ScIVP.m}, and a function for
  approximating the Jacobian with finite differences
  \texttt{ScJacobian.m}.
\end{itemize}
\end{document}